\documentclass[ECP]{ejpecp} 
\usepackage{amsmath,amssymb,
}
\usepackage{mathrsfs}
\usepackage{geometry}
\usepackage{bm}
\usepackage{mathtools}
\DeclareMathAlphabet\mathbfcal{OMS}{cmsy}{b}{n}
\DeclareMathAlphabet\mathbbcal{OMS}{cmsy}{b}{n}

\usepackage{upgreek}
\usepackage[T1]{fontenc}
\usepackage[francais,english]{babel}

\newtheorem{Definition}{Definition}[part]

\newtheorem{Assumption}{Assumption}[part]

\newtheorem{Example}{Example}[part]

\newcommand{\limsupn}{\limsup_{n\rightarrow\infty}}

\def\cal#1{\mathcal{#1}}

\def \F{I\!\!F}
\def \H{I\!\!H}

\def \N{I\!\!N}
\def \R{I\!\!R}

\def\Fc{{\cal F}}

\def\T{{\cal T}}

\def\Dzw1#1{\frac{\partial^2 #1}{\partial z \partial w_1}}

\def\Dzb1#1{\frac{\partial^2 #1}{\partial z \partial b_1}}

\newcommand{\dproof}{\noindent {Proof.} \quad}
\newcommand{\fproof}{\hfill $\square$ \bigskip}

\def\RB{\mathbb{R}}

\def\BC{\mathcal{B}}

\def\PC{\mathcal{P}}

\def\R{{\bf R}}

\def\1B{\text{1\!\!I}}

\def\tN{\tilde{N}}

\def\RB{\mathbb{R}}

\def\PC{\mathcal{P}}

\def\R{{\bf R}}

\def\1B{\text{1\!\!I}}

\def\tN{\tilde{N}}

\def\RB{\mathbb{R}}

\newcommand{\cf}{\mathcal{F}}

\newcommand{\ct}{\mathcal{T}}
\newcommand{\stopt}{\mathcal{T}_{t}}
\newcommand{\stops}{\mathcal{T}_{S,T}}
\newcommand{\stopo}{{\cal T}_0}
\newcommand{\Bwedge}{\bm{\wedge\mkern-13mu\wedge}}
\DeclareMathOperator*{\esssup}{ess\,sup}
\DeclareMathOperator*{\essinf}{ess\,inf}
\DeclareMathOperator{\e}{e}

\newcommand{\vertiii}[1]{{\left\vert\kern-0.25ex\left\vert\kern-0.25ex\left\vert #1 
    \right\vert\kern-0.25ex\right\vert\kern-0.25ex\right\vert}}

\newcommand{\vvertiii}[1]{{\vert\kern-0.25ex\vert\kern-0.25ex\vert #1 
    \vert\kern-0.25ex\vert\kern-0.25ex\vert}}

\newcommand{\triple}{\vert\kern-0.25ex\vert\kern-0.25ex\vert}

\newcommand{\Ref}{{\mathcal{R}}ef}

\DeclareMathAlphabet{\mathpzc}{OT1}{pzc}{m}{it}

\SHORTTITLE{DRBSDE\lowercase{s} and ${\cal E}^{\lowercase{f}}$-Dynkin games: beyond  right-continuity} 

\TITLE{Doubly Reflected BSDE\lowercase{s} and ${\cal E}^{\lowercase{f}}$-Dynkin games: beyond the right-continuous case} 

\AUTHORS{%
 Miryana~Grigorova\footnote{Bielefeld University, Bielefeld, Germany.
    \EMAIL{miryana_grigorova@yahoo.fr}}
  \and 
  Peter~Imkeller\footnote{Humboldt Universit\"at zu Berlin, Berlin, Germany. \EMAIL{imkeller@math.hu-berlin.de} }
  \and
Youssef~Ouknine  \footnote{Universit\'e Cadi Ayyad, Marrakech, Morocco. \EMAIL{ouknine@ucam.ac.ma} }
\and 
Marie-Claire~Quenez \footnote{Universit\'e  Paris-Diderot, Paris, France \EMAIL{quenez@lpsm.paris} }
  }


\AMSSUBJ{60G40; 93E20; 60H30} 
\AMSSUBJSECONDARY{60G07; 47N10} 

\VOLUME{}
\YEAR{}
\PAPERNUM{}
\DOI{}

\ABSTRACT{We formulate a notion of doubly reflected BSDE in the case where the barriers $\xi$ and $\zeta$ do not satisfy any regularity assumption and with general filtration.   Under a technical assumption (a Mokobodzki-type
condition), we show existence and uniqueness of the solution. 
In the case where $\xi$ is right upper-semicontinuous and $\zeta$ is right lower-semicontinuous, the solution is  characterized in terms of the value of a corresponding $\mathbfcal{E}^f$-Dynkin game, i.e. a game problem over stopping times with (non-linear) $f$-expectation, where $f$ is the driver of the doubly reflected BSDE.  In the general
case where the barriers do not satisfy any regularity assumptions, the solution of the doubly reflected BSDE is related to the value of  "an extension" of the previous  non-linear game problem over a
larger set of "stopping strategies" than the set of stopping times.
This characterization is then used to establish a comparison result and \textit{a priori} estimates with universal constants. 
}

\begin{document}
\section{Introduction}\label{Introduction}
   Backward stochastic differential equations (BSDEs) have been introduced in the case of a linear driver in \cite{Bismut2}, and then generalized to the non-linear case by Pardoux and Peng \cite{Pape90}. The theory of BSDEs provides a useful tool for the  study of financial problems such as the pricing of European options among others (cf., e.g., \cite{KPQ} and \cite{EQ96}). When the  driver $f$ is non-linear,  a BSDE induces a useful  family of non-linear operators, first introduced in \cite{EQ96} under the name of {\em non linear pricing system},
and later called $f$-{\em evaluation} (also,  $f$-{\em expectation}) and denoted by $\mathbfcal{E}^f$ (cf. \cite{Pe04}). 
  Reflected BSDEs (RBSDEs) are a variant of BSDEs in which the solution is constrained to be greater than or equal to a given process called {\em obstacle}. RBSDEs have been introduced in \cite{ElKaroui97} in the case of a Brownian filtration and a continuous obstacle, and  links with  (non-linear) optimal stopping problems with $f$-expectations  have been given in \cite{EQ96}. RBSDEs have been generalized to the case of a not necessarily continuous obstacle  
   and/or a larger filtration than the Brownian one by several authors  \cite{Ham}, \cite{CM}, \cite{HO1}, \cite{Essaky}, \cite{HO2}, \cite{17}. 
In all these works, the obstacle has been assumed to be right-continuous.  
The paper \cite{MG} is the first to study RBSDEs
 beyond the right-continuous case: there, we work under the assumption that the obstacle is only   right-uppersemicontinuous.  
 In \cite{MG3}, we address the case  where the obstacle does not satisfy any 
regularity assumption and with general filtration.  Existence and uniqueness of the solution in the irregular case is also shown in \cite{nouveau} (in the Brownian framework) by using a different approach. In \cite{MG} and  \cite{MG3}, 
links with  optimal stopping problems with $f$-expectations are also provided. 


  Doubly reflected BSDEs (DRBSDEs) have been introduced by Cvitanic and Karatzas in $\cite{CK}$ in the case of continuous barriers and a  Brownian filtration.   The solutions of such equations are constrained to stay between two adapted processes $\xi$ and $\zeta$, called {\em barriers}, with $\xi \leq \zeta$ and $\xi_T= \zeta_T$.  In the case of  non-continuous barriers  
   and/or a filtration associated with a Brownian motion and a random Poisson measure, DRBSDEs have been studied by several authors, cf. \cite{BHM}, \cite{HH}, \cite{HHd}, \cite{HL}, \cite{HHO}, \cite{CM}, \cite{EHO}, \cite{HO2}, \cite{DQS2}.  In all of the above-mentioned works on DRBSDEs, the barriers are assumed to be at least right-continuous.

 In the first part of the present paper, we formulate a notion of  doubly RBSDEs  in the case where the barriers do not satisfy any 
regularity assumption, and where the filtration is general. This allows for more flexibility in the modelling (compared to the cases of more regular payoffs and/or of particular filtrations).
We show  existence and  uniqueness of the solution of these equations. To this purpose,  we first consider the case where the driver does not depend on the solution, and is thus given by an adapted process $(f_t)$. We show that in this particular case, the solution of the DRBSDE can be written in terms of the difference of the solutions of a coupled system
of two reflected BSDEs. We  show that this system (and hence the Doubly Reflected BSDE) admits a solution if and only if the so-called \textit{Mokobodzki's condition} holds (assuming  the existence of two strong supermartingales whose difference is between $\xi$ and  $\zeta$).  
We then provide  \textit{a priori} estimates for our doubly RBSDEs, by using Gal'chouk-Lenglart's formula (cf.  Corollary A.2 in \cite{MG}). From these estimates, we derive \emph{the uniqueness of the solution} of the doubly RBSDE associated with driver process $(f_t)$. 
We then solve  the case of a general Lipschitz driver $f$ by using the \textit{a priori} estimates
 and   Banach fixed point theorem.

 In the second part of the paper, we  focus on links between the solution of the doubly  reflected BSDE with irregular  barriers from the first part  and some related two-stopper-game problems.\\
Let us first recall the "classical" Dynkin game problem which has been largely studied (cf., e.g., \cite{ALM} for general results). \\
Let  
${\cal T}_0$ denote the set of all stopping times valued in $[0,T]$, where $T>0$. 
For each pair $(\tau, \sigma) \in {\cal T}_0\times {\cal T}_0$, the terminal time of the game is given by  $\tau \wedge \sigma$ and the terminal payoff, or reward,  of the game (at time $\tau \wedge \sigma$) is given by 
\begin{equation}\label{eq0}
I(\tau, \sigma):=  \xi_{\tau}\textbf{1}_{\{ \tau \leq \sigma\}}+\zeta_{\sigma}\textbf{1}_{\{\sigma<\tau\}}.
\end{equation}
The criterion is  defined as the (linear) expectation of the pay-off, that is, $ E\, [I(\tau, \sigma)]$.
It is well-known that, if $\xi$ is right upper-semicontinuous (right u.s.c) and $\zeta$ is right lower-semicontinuous (right l.s.c) and satisfy Mokobodzki's condition, this classical Dynkin game has a (common) value, that is, the following equality holds:
\begin{equation}\label{eq0bis}
\inf_{\sigma \in {\cal T}_0}\sup_{\tau \in {\cal T}_0}   E\, [ I(\tau, \sigma)]= \sup_{\tau \in {\cal T}_0} \inf_{\sigma \in {\cal T}_0}   E\, [ I(\tau, \sigma)]. 
\end{equation}
Moreover,  under the additional assumptions that  $\xi$ is left-uppersemicontinuous (left-u.s.c), $\zeta$ is left-lowersemicontinuous (left-l.s.c), both along stopping times, and $\xi_t < \zeta_t$, $t <T$, there exists a saddle point 
(cf. \cite{ALM}, \cite{KQC}). \footnote{Actually, the strict separability condition on $\xi$ and $\zeta$ is not necessary to ensure the existence of a saddle point (cf. Remark 3.8 in \cite{DQS2} when $\xi$ and $\zeta$ are right-continuous).}

Furthermore,  when the processes $\xi$ and $\zeta$ are 
right-continuous, the (common) value of the classical  Dynkin game is equal to the solution at time $0$ of the doubly  reflected BSDE with driver equal to $0$ and barriers $(\xi,\zeta)$ (cf. \cite{CK},\cite{HL},\cite{Lepeltier-Xu}).

  


In the second part of the present paper, we consider  the following  generalization of the classical Dynkin game problem:
 For each pair $(\tau, \sigma) \in {\cal T}_0\times {\cal T}_0$, the criterion is defined by 
$
   \mathbfcal{E}^{^f}_{_{0,\tau \wedge \sigma}}[I(\tau, \sigma)],
$
 where $\mathbfcal{E}^f_{0, \tau\wedge \sigma}(\cdot)$ denotes  the $f$-expectation  at time $0$ when the  terminal time is $\tau \wedge \sigma$. We refer to this generalized game problem  as    $\mathbfcal{E}^f$-{\em Dynkin game}.
   This  non-linear game problem has been introduced in \cite{DQS2} in the case where $\xi$ and $\zeta$ are right-continuous under the name of \textit{generalized Dynkin game}, the term\textit{ generalized} referring to the presence of a (non-linear) $f$-expectation  in  place of the "classical" linear expectation. 

 In the second part of the  paper, we first generalize  the results of \cite{DQS2} beyond the right-continuity assumption on  $\xi$ and $\zeta$ 
 (and in the case of a general filtration). More precisely, by using results from the first part of the present paper, combined with some arguments from \cite{DQS2}, we show that  if $\xi$ is right-u.s.c. and $\zeta$ is right-l.s.c.\,, and if they satisfy Mokobodzki's condition, there exists a (common) value function for the $\mathbfcal{E}^f$-{\em Dynkin game}, that is 
  \begin{equation} \label{eq_intro2}
 \inf_{\sigma \in {\cal T}_0}\sup_{\tau\in {\cal T}_0} \mathbfcal{E}^{^f}_{_{0,\tau \wedge \sigma}}[I(\tau, \sigma)]= \sup_{\tau \in {\cal T}_0} \inf_{\sigma \in {\cal T}_0} \mathbfcal{E}^{^f}_{_{0,\tau \wedge \sigma}} [I(\tau, \sigma)].
  \end{equation}
and this common value is  equal to the solution at time $0$ of the doubly reflected BSDE with driver $f$ and barriers $(\xi,\zeta)$ from the first part of the paper. Moreover, under the additional assumption  that  $\xi$ is  left u.s.c. along stopping times and $\zeta$ is left l.s.c. along stopping times, we prove that there exists a saddle point for the $\mathbfcal{E}^f$-{\em Dynkin game}. 
Let us note that in the particular case when $f=0$, our results  on existence of a common value and on existence of saddle points correspond to the results from the literature on classical Dynkin games recalled above. However, even in the case when $f=0$, our characterization of the value of the classical Dynkin game \eqref{eq_intro2} with right u.s.c. payoffs via a  (linear) DRBSDE  is new; it generalizes the well-known result shown in  \cite{CK} (see also \cite{HL} and\cite{Lepeltier-Xu}) in the case of continuous (resp. right-continuous) payoffs.

We then turn to the general case where $\xi$ and $\zeta$ are completely irregular, which is technically more difficult. 
Indeed, in this case, the Dynkin game problem considered above may not have a value, that is, the equality \eqref{eq_intro2} may not hold true. This is already well known from the simpler case of a classical Dynkin game (cf. e.g. \cite{ALM})\footnote{  see Example \ref{contre-exemple} in the present paper, and also Section 3 " The defect of value" in \cite{ALM}.}.\\
An interesting question is the question of how to  interpret   the  solution of the doubly  reflected BSDE with completely irregular barriers  $(\xi,\zeta)$ in terms of a game problem.
To this aim, we formulate "an extension" of the  previous $\mathbfcal{E}^f$-Dynkin game problem over a larger set of "stopping systems" than the set of  stopping times ${\cal T}_0$. We show that this extended game has a common value which coincides with the solution of our general DRBSDE with irregular barriers, and that it admits $\varepsilon$-saddle points (of stopping systems). 
 Using these results, we prove a
comparison theorem and \textit{a priori} estimates with universal constants for DRBSDEs with completely irregular barriers.  
  

In the last section, we give an application to the pricing of game options in a complete imperfect market model beyond the case of right-continuous pay-off. 
Since Kifer's seminal work \cite{Kifer},
it is well-known that if the market model is complete and if the processes  $\xi$ and $\zeta$ are right-continuous and satisfy  Mokobodzki's condition, then  the  price of the game option (up to a discount factor)  is equal to the common value of the classical Dynkin game from equation \eqref{eq0bis}, where the expectation is taken under the unique martingale  measure of the model.
In a large class of   market models with imperfections,  European options can be priced  via an  $f$-expectation/evaluation  $\mathbfcal{E}^{^f}$, where $f$ is a \emph{nonlinear} driver in which the imperfections are encoded.
In such a framework,  the problem of  pricing of game options  has been considered in \cite{DQS3}:  when $\xi$ and $\zeta$ are right-continuous and satisfy  Mokobodzki's condition, the common value of the $\mathbfcal{E}^{f}$-{\em Dynkin game} from equation \eqref{eq_intro2}  is shown to be equal to the  price of the game option, that is the infimum of the initial wealths which allow the seller to be super-hedged (cf. \cite{DQS3}).
Using the results of the present paper, we show that the result of \cite{DQS3}  can be generalized to the case when $\xi$ is right-u.s.c. and $\zeta$ is left-l.s.c.
For example, this result can be applied to the case of cancellable  American call options with lower barrier, for which the payoff $\xi$ is not right-continuous.

 The remainder of the paper is organized as follows: In Section \ref{sec2}, we introduce the notation and some definitions.
 In Section  \ref{sec4}, we provide  first results on doubly reflected BSDEs associated with a Lipschitz driver and barriers  $(\xi,\zeta)$ which do not satisfy any regularity assumption; in particular, we 
 show existence and uniqueness of the  solution of this equation. 
   Section  \ref{sec5} is dedicated to the interpretation of  the solution in terms of a two-stopper game problem,  first in the case when $\xi$ is right u.s.c. and $\zeta$ is right l.s.c., then in the case where they 
   do not satisfy any regularity assumption. 
   In Section  \ref{sec6}, we provide a comparison theorem  and \textit{a priori} estimates with universal constants for our  doubly reflected BSDEs with irregular barriers. In Section \ref{gameoptions}, we give an application of our results to the pricing of game options with irregular payoffs in an imperfect market.
   The Appendix contains  some useful results on reflected BSDEs with an irregular obstacle and also some of  the proofs.   

\section{Preliminaries}\label{sec2}

Let $T>0$ be a fixed positive real number.  Let  $E= \R^n\setminus\{0\}, \mathscr{E}= {\cal B}({\bf R}^n\setminus\{0\})$, which we    
equip with a $\sigma$-finite positive measure $\nu$.
Let $(\Omega,  \Fc, P)$ be a probability space equipped with a right-continuous complete filtration $\F = \{\Fc_t \colon t\in[0,T]\}$. Let $W$ be
a one-dimensional $\F$-Brownian motion $W$, and let $N(dt,de)$ be an  $\F$-Poisson random measure  with compensator $dt\otimes\nu(de)$, supposed to be independent from $W$.
We denote by $\tilde N(dt,de)$ the compensated process, i.e. $\tilde N(dt,de):= N(dt,de)-dt\otimes\nu(de).$
 The notation $L^2({\cal F}_T)$  stands for the space of random variables which are  $\Fc 
_T$-measurable and square-integrable. For $t\in [0,T],$ we denote by $\stopt$ the set of stopping times $\tau$ such that $P(t \leq\tau\leq T)=1.$ More generally, for a given stopping time $S\in \stopo$, we denote by $\ct_{S,T}$ the set of stopping times $ \tau$ such that $P(S \leq\tau\leq T)=1.$   \\





 We also use the following notation: 
\begin{itemize}
\item ${\cal P}$ (resp. $\mathcal{O}$) is  the predictable (resp. optional) $\sigma$-algebra
on $ \Omega\times [0,T]$.

\item $L^2_\nu$ is the set of $(\mathscr{E}, \mathcal{B}(\R))$-measurable functions $\ell:  E \rightarrow \R$ such that  $\|\ell\|_\nu^2:= \int_{ E}  |\ell(e) |^2 \nu(de) <  \infty.$
For $\ell\in L^2_\nu$, $\mathpzc{k}\in L^2_\nu$, we define $\langle \ell, \mathpzc{k}\rangle_\nu:=\int_E \ell(e)\mathpzc{k} (e) \nu(de)$.
\item    $\H^{2}$ is the set of 
$\R$-valued predictable processes $\phi$ with
 $\| \phi\|^2_{\H^{2}} := E \left[\int_0 ^T |\phi_t| ^2 dt\right] < \infty.$ 
 \footnote{By a slight abuse of notation, we shall also write $\| \phi\|^2_{\H^{2}}$ for  $E \left[\int_0 ^T |\phi_t| ^2 dt\right]$ in the case of a progressively measurable 
real-valued process $\phi$. 
}


\item $\H_{\nu}^{2}$ is  the set of $\R$-valued processes $l: (\omega,t,e)\in(\Omega\times[0,T] \times  E)\mapsto l_t(\omega, e)$ which are {\em predictable}, that is $(\PC \otimes {\mathscr{E}},\BC(\R))$-measurable,
and such that $\| l \|^2_{\H_{\nu}^{2}} :=E\left[ \int_0 ^T \|l_t\|_{\nu}^2 \,dt   \right]< \infty.$

\item Let ${\cal M}^{2}$  be the set of square integrable  martingales $M= (M_t)_{t\in[0,T]}$ with $M_0=0$. This is a Hilbert space equipped with the scalar product 
$( M, M' ) _{{\cal M}^2}:= E[M_T M'_T]\,(= E[\, \langle M, M' \rangle_T]= E( \, [M, M' ]_T))$, for 
$M,M'$ $\in$ ${\cal M}^{2}$ (cf., e.g., \cite{Protter} IV.3). For each $M \in {\cal M}^{2}$, we set  $\| M\|^2_{{\cal M}^{2}}:= E( M_T^2)$. \\
\item Let ${\cal M}^{2,\bot}$ be the subspace of martingales $h \in {\cal M}^{2}$
satisfying  $\langle h, W \rangle_\cdot =0$, and such that,
for all predictable processes $l$ $\in$ $\H_{\nu}^{2}$, 
\begin{equation} \label{ortho}
\langle h, \int_0^\cdot \int_E l_s(e) \tilde N(ds, de) \rangle_t =0,\,\,\, 0 \leq t \leq T \quad {\rm a.s.}
\end{equation}

\end {itemize}

\begin{remark} \label{2i} Note that condition \eqref{ortho} is equivalent to the fact that the square bracket process $[ \,h\,,  \int_0^\cdot \int_E l_s(e) \tilde N(ds, de)\,]_t$ is a martingale.
\footnote{
Recall  also that the condition $\langle h, W \rangle_\cdot =0$ is equivalent to the orthogonality of $h$ (in the sense of the scalar product $( \cdot, \cdot ) _{{\cal M}^2}$) with respect to 
all stochastic integrals of the form 
$\int_0^\cdot z_s dW_s$, where 
$z \in \H^{2}$ (cf.\, e.g.\,, \cite{Protter} IV. 3 Lemma 2). 
Similarly, the condition \eqref{ortho} is equivalent to the orthogonality of $h$ with respect to 
all stochastic integrals of the form 
$\int_0^\cdot \int_E l_s(e) \tilde N(ds, de)$, where 
$l$ $\in$ $\H_{\nu}^{2}$ (cf., e.g., the Appendix 
in \cite{MG3}).}
\end{remark}

As in \cite{MG}, we  denote by  $ {\cal S}^{2}$ the vector space of $\R$-valued optional 
(not necessarily cadlag)
 processes $\phi$ such that
$\vertiii{\phi}^2_{{\cal S}^{2}} := E[\esssup_{\tau\in\T_0} |\phi_\tau |^2] <  \infty.$ 
By Proposition 2.1 in \cite{MG}, the mapping $\vertiii{\cdot}_{{\cal S}^{2}}$ is a norm on the space $ {\cal S}^{2}$,  and ${\cal S}^{2}$ endowed with this norm is a Banach space.  \\

 We recall the following orthogonal decomposition property of martingales in 
${\cal M}^2$
(cf. Lemma III.4.24 in \cite{JS}). 
\begin{lemma}[Orthogonal decomposition of martingales in 
${\cal M}^2$]\label{theoreme representation}
For each $M \in {\cal M}^2$, there exists a unique triplet 
$(Z,l,h) \in \H^{2} \times \H_{\nu}^{2} \times {\cal M}^{2, \bot}$ such that
\begin{equation}
\label{equation representation}
M_t  = \int_0^t Z_sdW_s+\int_0^t \int_{E} l_t(e) \tilde N(dt,de) +h_t\,, \quad\forall\,t\in[0,T] \quad a.s. 
\end{equation}
\end{lemma}

\begin{definition}[Driver, Lipschitz driver]\label{defd}
A function $f$ is said to be a {\em driver} if 
\begin{itemize}
\item  
$f: \Omega  \times [0,T]  \times \R^2 \times L^2_\nu \rightarrow \R $\\
$(\omega, t,y, z, \mathpzc{k}) \mapsto  f(\omega, t,y, z, \mathpzc{k})  $
  is $ {\cal P}\otimes {\cal B}(\R^2)  \otimes {\cal B}(L^2_\nu) 
- $ measurable,  
\item $E[\int_0^T f(t,0,0,0)^2dt] < + \infty$.
\end{itemize} 
A driver $f$ is called a {\em Lipschitz driver} if moreover there exists a constant $ K \geq 0$ such that $dP \otimes dt$-a.e.\,, 
for each $(y_1, z_1, \mathpzc{k}_1)\in \R^2 \times L^2_\nu$, $(y_2, z_2, \mathpzc{k}_2)\in \R^2 \times L^2_\nu$, 
$$|f(\omega, t, y_1, z_1, \mathpzc{k}_1) - f(\omega, t, y_2, z_2, \mathpzc{k}_2)| \leq 
K (|y_1 - y_2| + |z_1 - z_2| +   \|\mathpzc{k}_1 - \mathpzc{k}_2 \|_\nu).$$
\end{definition}

 We recall the definition of a non reflected BSDE in the case of a general filtration
(cf. Remark 12.1 in the Appendix of \cite{MG3} for the existence and the uniqueness of the solution). \begin{Definition}[BSDE, conditional $f$-expectation]\label{BSDE}
If   $f$ is a Lipschitz driver and if $\xi$ is  in $L^2({\cal F}_T)$, then there exists a unique solution $(X,\pi, l,h)\in{\cal S}^2 \times \H^2 \times \H^2_\nu\times {\cal M}^2$ to the following BSDE:
$-dX_t= f(t, X_t,  \pi_t, l_t)dt-  \pi_t dW_t-\int_E  l_t(e) \tilde N(dt,de)-dh_t; \quad X_T= \xi.$\\ 
For $t\in[0,T]$, the (non-linear) operator $\mathcal{E}^f_{t,T}(\cdot): L^2(\cf_T)\rightarrow L^2(\cf_t)$ which maps a given terminal condition
$\xi\in L^2(\cf_T)$ to the position $X_t$ (at time $t$) of the first component  of the solution of the above BSDE is called  \emph{conditional $f$-expectation at time $t$}.  As usual, this notion  can be extended to
the case where the (deterministic) terminal time 
$T$
is replaced by a (more general) stopping time $\tau\in\stopo$,  $t$ is replaced by a stopping time $S$ such that $S\leq \tau$ a.s. and  the domain $L^2(\cf_T)$ of the operator    is replaced  by $L^2(\cf_\tau)$.

\end{Definition}

We recall the following definition from \cite{DQS2}.
\begin{definition}\label{proba}
 Let $A=(A_t)_{0 \leq t \leq T}$ and $A'=(A_t')_{0 \leq t \leq T}$ be two real-valued  optional non-decreasing cadlag processes with $A_0=0$, $A'_0=0$ and  $E[A_T]<\infty$ and $E[A'_T]<\infty$. We say that the random measures $dA_t$ and $dA'_t$ are {\em mutually singular}, and we write $dA_t \perp dA'_t $, if there exists $D \in \mathcal{O}$ such that:
 \begin{equation}\label{eq_sing}
 E[\int_0^T \textbf{1}_{D^c} dA_t]= E[\int_0^T \textbf{1}_{D} dA'_t]=0,  
 \end{equation}
which can also be written as $\int_0^T \textbf{1}_{D^c_t} dA_t= \int_0^T \textbf{1}_{D_t} dA'_t=0 \,\, \text{ a.s.}\,,  $
where for each $t \in [0,T]$, $D_t$ is the {\em section at time $t$} of $D$, that is, $D_t:= \{\omega \in \Omega \,,\, (\omega,t) \in D \}$.
\end{definition}
For real-valued random variables $X$ and $X_n$, $n \in \N$, the notation   "$X_n\uparrow X$" stands  for "the sequence $(X_n)$ is nondecreasing and converges to $X$ a.s.".  \\
For a ladlag process $\phi$, we denote by $\phi_{t+}$ and $\phi_{t-}$ the right-hand and left-hand limit of $\phi$ at  $t$. We denote by $\Delta_+ \phi_t:=\phi_{t_+}-\phi_t$ the size of the right jump of $\phi$ at $t$, and by   $\Delta \phi_t:=\phi_t-\phi_{t-}$ the size of the left jump of $\phi$ at $t$.

\begin{definition} \label{defr} 
%
An optional process $(\phi_t)$ is said to be  left upper-semicontinuous (resp. left lower-semicontinuous) along stopping times if  for each $\tau \in \T_0$, for each nondecreasing sequence of stopping times $ (\tau_n)$ such that $\tau_n \uparrow \tau$,  a.s.\,, we have
$\phi_{\tau} \geq \limsup_{n\to \infty} \phi_{\tau_n} $ (resp. $\phi_{\tau} \leq \liminf_{n\to \infty} \phi_{\tau_n}$) a.s.
\end{definition}

\begin{remark}\label{Rmk_final_00}
If the process $( \phi_{t})$ has left limits, $( \phi_{t})$ is left upper-semicontinuous (resp. left lower-semicontinuous) along stopping times if and only if for each predictable stopping time $\tau \in \T_0$,  $\phi_{\tau-} \leq \phi_{\tau}$ (resp. $\phi_{\tau-} \geq \phi_{\tau}$) a.s.
\end{remark}

\begin{definition}[Strong supermartingale]
An optional  process $\phi_.=(\phi_t )$ 
belonging to ${\cal S}^2$ is said to be a {\em strong supermartingale} if  for all 
$\theta, \theta'$ $ \in$ $\T_0$ such that $\theta \geq \theta'$ a.s., 
$
E[\phi_{\theta} \, |\,{\cal F}_{ \theta' }] \leq \phi_{\theta'} \quad
\mbox{a.s.}   
$
\end{definition}
We recall that a  strong supermartingale in $ \mathcal{S}^2$ is necessarily right upper-semicontinuous (cf., e.g., \cite{DM2}).

  For the easing of the presentation, we define the relation  \;$\geq$\; for processes in ${\cal S}^2$  as follows: 
for $\phi, \phi' \in {\cal S}^2$,   we write $\phi \leq \phi'$, if $\phi_t \leq \phi'_t$ for all $t\in [0,T]$ a.s. 
Similarly, we define the relations \;$\leq$\; and \;$=$ \; on ${\cal S}^2$.
\section{Doubly Reflected BSDE whose obstacles are  irregular}\label{sec4}
 \subsection{Definition and first properties}

Let $T>0$ be a fixed terminal time (as before). Let $f$ be  a  driver. 
Let $\xi= (\xi_t)_{t\in[0,T]}$ and $\zeta=(\zeta_t)_{t\in[0,T]}$ be two   processes in ${\cal S}^2$ such that $\xi_t\leq \zeta_t, 0\leq t\leq T,$ a.s. and $\zeta_T=\xi_T$ a.s.
A pair of processes $(\xi, \zeta)$ satisfying the previous properties will be called a \emph{ pair of admissible  barriers}, or a \emph{ pair of admissible  obstacles}.

For each process  $\phi \in  {\cal S}^2$, the process
$(\overline{\phi}_t)$ (resp. $(\underline{\phi}_t)$) denotes the {\em left upper- (resp.  left lower-) semicontinuous envelope} of the process $\phi$, defined by  
$\overline{\phi}_t:=\limsup_{s \uparrow t, s< t} \phi_s$ (resp. $\underline{\phi}_t:=\liminf_{s \uparrow t, s< t} \phi_s$), for all $t\in]0,T]$. 
The process $(\overline{\phi}_t)$ (resp. $(\underline{\phi}_t)$) is predictable (cf. \cite[Thm. 90, page 225]{DM1}) and left upper- (resp. left lower-) semicontinuous. \\
Note that when $\phi$ is left-limited, we have $\overline{\phi}_t= \underline{\phi}_t = \phi_{t^-}$ for all $t\in]0,T]$ a.s.

\begin{definition} \label{def_solution_DRBSDE}
A process $(Y,Z,k,h, A,C,A',C')$ is said to be a solution to the doubly reflected BSDE with parameters $(f,\xi,\zeta)$, where $f$ is a driver and $(\xi,\zeta)$ is a pair of admissible obstacles, if 
\begin{align}\label{DRBSDE}
&(Y,Z,k,h,A,C,A',C')\in {\cal S}^2 \times  \H^2 \times \H^2_\nu \times {\cal M}^{2, \bot} \times ({\cal S}^2)^2\times ({\cal S}^2)^2  \text{and}  \text{ a.s.} \text{ for all } t\in[0,T] \nonumber \\
&Y_t=\xi_T+\int_{t}^T f(s,Y_s,  Z_s,k_s)ds-\int_{t}^T  Z_s dW_s-\int_{t}^T \int_E k_s(e) \tilde N(ds,de)- (h_T-h_t) \,+ \nonumber
\\&+A_T-A_t-(A'_T-A'_t)+C_{T-} -C_{t-}-(C'_{T-} -C'_{t-}),\\
&  \xi_t \leq Y_t \leq \zeta_t,   \text{ for all } t\in[0,T]  \text{ a.s.,} \label{RBSDE_inequality_barrier}\\
& A \text{ and } A' \text{ are nondecreasing right-continuous predictable processes 
with } A_0=  A'_0= 0, \nonumber\\
& \int_0^T {\bf 1}_{\{Y_{t-} > \overline{\xi}_t\}} dA_t = 0 \text{ a.s. and } \; \int_0^T {\bf 1}_{\{Y_{t-} < \underline{\zeta}_t\}} dA'_t = 0 \text{ a.s.}\label{RBSDE_A}\\   
& C \text{ and } C' \text{ are nondecreasing right-continuous adapted purely discontinuous processes with }   \nonumber\\
& C_{0-}=C'_{0-}= 0,  \nonumber \\
&(Y_{\tau}-\xi_{\tau})(C_{\tau}-C_{\tau-})=0 \text{ and } (Y_{\tau}-\zeta_{\tau})(C'_{\tau}-C'_{\tau-})=0 \text{ a.s. for all }\tau\in\T_0,\label{RBSDE_C}\\
&dA_t \perp dA'_t \;\text { and }\; dC_t \perp dC'_t. \label{RBSDE_D}
\end{align}
\end{definition}
Here  $A^c$ denotes the continuous part of the process $A$ and $A^d$  its discontinuous part.\\
Equations \eqref{RBSDE_A} and \eqref{RBSDE_C} are referred to as \emph{minimality conditions} or \emph{Skorokhod conditions}.  \\
Let us note that if  $(Y,Z,k,h,A,C,A',C')$ satisfies \eqref{DRBSDE}, then the process $Y$
has left and right limits.

\begin{remark}[Left-limited case] \label{Left-limited-case} 
We note that when $\xi$ and $\zeta$ are  left-limited processes, then $\overline{\xi}_t$ (resp. $\underline{\zeta}_t$)  can be replaced by $\xi_{t-}$ (resp. $\zeta_{t-}$) in the Skorokhod conditions \eqref{RBSDE_A}.
\end{remark}


\begin{remark} \label{sautsY}
When $A$ and $A'$ (resp. $C$ and $C'$) are not required to be mutually singular, they can simultaneously increase 
on $\{\overline{\xi}_{t}=\underline{\zeta}_{t}\}$ (resp. on $\{\xi_{t}=\zeta_{t}\}$). 
The constraints $dA_t \perp dA_t'$  and $dC_t \perp dC'_t$ will allow us to obtain the uniqueness of the nondecreasing processes $A$, $A'$, $C$ and $C'$ without the  strict separability condition $\xi < \zeta$. \\
We note also that, due to Eq.  \eqref{DRBSDE}, we have $\Delta C_t-\Delta C'_t=-(Y_{t+}-Y_t)=-\Delta_+Y_t$. This, together with the condition $dC_t \perp dC_t'$ gives  $\Delta C_t=(Y_{t+}-Y_t)^-$ for all $t$ a.s., and $\Delta C'_t=(Y_{t+}-Y_t)^+ $ for all $t$ a.s. 
\end{remark}
\begin{remark}\label{sautsYA} (Quasi-left-continuous filtration)
When the filtration is quasi-left-continuous (for example, when the filtration is the natural filtration associated $W$ and $N$ \footnote{as it is the case in the literature on RBSDEs and DRBSDEs.}), martingales have only totally inaccessible jumps. Hence, in this case, for each predictable $\tau\in\T_0$, $\Delta A_\tau-\Delta A_\tau^{'}=-\Delta Y_\tau$ (cf. Eq.  \eqref{DRBSDE}). 
  This, together with the condition  $dA_t \perp dA_t'$, ensures that for each predictable  $\tau \in \T_0$, 
$\Delta A_\tau=(\Delta Y_\tau)^- $ and $\Delta A_\tau^{'}=  (\Delta Y_\tau)^+ \text{ a.s.}$ \\
 We stress that in the case of a general filtration, this property does not necessarily hold. Indeed, by equation  \eqref{DRBSDE}, for each predictable $\tau\in\T_0$, we have $\Delta Y_\tau =-\Delta A_\tau+\Delta A_\tau^{'} + \Delta h_\tau$ a.s.\, and  $\Delta h_\tau$ is here not necessarily equal to $0$, since in this case martingales may admit jumps at some predictable stopping times. 
\end{remark}

\begin{proposition}\label{continuous} Let $f$ be a driver and $(\xi,\zeta)$ be a pair of admissible obstacles.\\
 Let $(Y,Z,k,h,A,C,A',C')$ be a solution to the doubly reflected BSDE with parameters $(f,\xi,\zeta)$.\\
(i) For each $\tau \in \T_0$, we have
 $$ Y_\tau = (Y_{\tau^+} \vee \xi_\tau)  \wedge \zeta_\tau \quad {\rm a.s.}$$
(ii)  If $\xi$ (resp. $\zeta$) is right continuous, then $C=0$ (resp. $C'=0$).\\
(iii) If $\xi$ (resp. $\zeta$) is left upper-semicontinuous (resp. left lower-semicontinuous) along stopping times, then the process $A$ (resp. $A'$) is continuous. 
\end{proposition}
Note that the assertion (iii) will be used to obtain Corollary \ref{optimal} and Proposition \ref{superhedging}.

\dproof  Let us show the first assertion. Let $\tau \in \T_0$. By the previous Remark \ref{sautsY}, we have $\Delta C_\tau=(Y_{\tau^+}-Y_\tau)^-$ and $\Delta C'_\tau=(Y_{\tau+}-Y_\tau)^+ $ a.s.\, Since $C$ and $C'$ satisfy the Skorokhod condition \eqref{RBSDE_C}, we have 
$$
(Y_{\tau^+}-Y_\tau)^-=\textbf{1}_{ \{ Y_{\tau}=\xi_{\tau}\}}(Y_{\tau^+}-Y_\tau)^- \,\,
\text{ and }\,\,(Y_{\tau^+}-Y_\tau)^+=\textbf{1}_{ \{ Y_{\tau}=\zeta_{\tau}\}}(Y_{\tau^+}-Y_\tau)^+ \,\,
\text{ a.s. }
$$
Hence, on the set $\{\xi_{\tau} < Y_{\tau} < \zeta_\tau \}$, we have $Y_\tau = Y_{\tau^+}$ a.s.\,, which implies that 
 $(Y_{\tau^+} \vee \xi_\tau)  \wedge \zeta_\tau= Y_\tau$ a.s. 
 Now, on the set $\{\xi_{\tau} < Y_{\tau} = \zeta_\tau \}$, we have $(Y_{\tau^+}-Y_\tau)^-=0$ a.s.\,, which gives 
$Y_{\tau^+} \geq Y_\tau = \zeta_\tau \geq \xi_\tau$ a.s.\,, which implies that 
$(Y_{\tau^+} \vee \xi_\tau)  \wedge \zeta_\tau$ $=$ $Y_{\tau^+} \wedge \zeta_\tau$ $=$ $\zeta_{\tau}$
$=$ $Y_{\tau}$ a.s.\, Similarly, on the set $\{\xi_{\tau} = Y_{\tau} < \zeta_\tau \}$, we have $(Y_{\tau^+} \vee \xi_\tau)  \wedge \zeta_\tau= Y_\tau$ a.s.\, The first assertion thus holds.

Let us show the second assertion. Suppose that $\xi$ is right-continuous. 
Let $\tau \in \T_0$. We show $\Delta C_\tau=0$ a.s.
As seen above, we have
$$
\Delta C_{\tau}=\textbf{1}_{ \{ Y_{\tau}=\xi_{\tau}\}}(Y_{\tau^+}-Y_\tau)^-=\textbf{1}_{ \{ Y_{\tau}=\xi_{\tau}\}}(Y_{\tau^+}-\xi_\tau)^-=\textbf{1}_{ \{ Y_{\tau}=\xi_{\tau}\}}(Y_{\tau^+}-\xi_{\tau^+})^- \text{ a.s., }
$$
where the last equality follows from the right-continuity of $\xi$. Since $  Y\geq \xi$, we derive that $\Delta C_{\tau} = 0$ a.s.\, This equality being true for all $\tau \in \T_0$, it follows that $C=0$. Similarly, it can be shown that if  $\zeta$ is right-continuous, then $C'=0$. 
Hence, the second assertion holds.


 It remains to show the third assertion.
Suppose that $\xi$ is left u.s.c.along stopping times. We show $\Delta A_\tau=0$ a.s.
 Let $\tau\in\T_0$ be a predictable stopping time. For each martingale $M$, we have 
$E[\Delta M_{\tau}/ {\cal F}_{\tau^-}] =0$ a.s.\,Moreover, since $A$ and $A'$ are predictable, we have 
$E[\Delta A_{\tau}/ {\cal F}_{\tau^-}]= \Delta A_{\tau}$ a.s. and $E[\Delta A'_{\tau}/ {\cal F}_{\tau^-}]= \Delta A'_{\tau}$ a.s. By \eqref{DRBSDE}, we derive that
\begin{equation} \label{ineqY}
E[\Delta Y_{\tau}/ {\cal F}_{\tau^-}]= -\Delta A_{\tau}+ \Delta A'_{\tau} = -\Delta A_{\tau}{\bf 1}_{\{Y_{\tau^-} = \overline{\xi}_\tau\} \cap D} + \Delta A'_{\tau} {\bf 1}_{\{Y_{\tau^-} =\underline{\zeta}_\tau\} \cap D'} 
\quad {\rm a.s.}
\end{equation}
where $D:= \{\Delta A_{\tau} >0\}$ and $D':= \{\Delta A'_{\tau} >0\}$. Note that the sets  $D $ and $D'$ belong to 
 ${\cal F}_{\tau^-}$. Since  $dA_t \perp dA_t'$, we get $D \cap D' = \emptyset$ a.s. 
 Hence, on  $\{Y_{\tau^-} = \overline{\xi}_\tau\} \cap D$, we  have 
 $$E[Y_{\tau}/ {\cal F}_{\tau^-}] - Y_{\tau^-}= E[\Delta Y_{\tau}/ {\cal F}_{\tau^-}]= -\Delta A_{\tau} \leq 0 \quad {\rm a.s.}$$
Since $\xi$ is left-u.s.c. along stopping times, we thus derive that on  $\{Y_{\tau^-} = \overline{\xi}_\tau\} \cap D$, we have
$$\overline{\xi}_\tau \leq E[\xi_{\tau}/ {\cal F}_{\tau^-}]  \leq E[Y_{\tau}/ {\cal F}_{\tau^-}] \leq Y_{\tau^-} \quad {\rm a.s.}$$
and the inequalities are even equalities (still on the set $\{Y_{\tau^-} = \overline{\xi}_\tau\} \cap D$).
Hence, $E[Y_{\tau}/ {\cal F}_{\tau^-}] = Y_{\tau^-}$ a.s. on  $\{Y_{\tau^-} = \overline{\xi}_\tau\} \cap D$. By \eqref{ineqY}, we derive that $\Delta A_{\tau}=0$ a.s. This equality being true for every predictable stopping time $\tau \in \T_0$, it follows that $A$ is continuous. Similarly, it can be shown that if  $\zeta$ is left lower-semicontinuous
 along stopping times, then $A'$ is continuous, which ends the proof.
\fproof

\begin{remark}
Note that the proof of the third assertion relies on different arguments from the proof given in the literature on DRBSDEs (cf. the proof of Theorem 3.7 (i) in \cite{DQS2}), which concerns the particular case when the filtration is the natural filtration associated $W$ and $N$, and is thus quasi-left-continuous. 
More precisely, the proof given in  \cite{DQS2} relies on Remark \ref{sautsYA}, which does not necessarily holds when the filtration is not quasi-left-continuous. 
\end{remark}

\begin{remark}[Right-continuous left-limited case] It follows from the second assertion in the above proposition that if $\xi$ and $\zeta$ are right-continuous, then $C=C'=0$.  If moreover $\xi$ and $\zeta$  are left-limited, by Remark \ref{Left-limited-case}, our Definition \ref{def_solution_DRBSDE}  thus corresponds to the one  given in the literature on DRBSDEs (cf. e.g. \cite{DQS2}).
\end{remark}

We now provide a necessary condition for the existence of a solution of the doubly reflected BSDE from Definition \ref{def_solution_DRBSDE}.

Let $(Y, Z,k, h, A, C,A',C')\in {\cal S}^2 \times \H^2 \times \H^2_\nu \times {\cal M}^{2, \bot} \times({\cal S}^2)^2\times ({\cal S}^2)^2$ be  a  solution to the DRBSDE  associated with  driver $f$ and with a pair of admissible  barriers $(\xi,\zeta)$. 
By taking the conditional expectation with respect to  ${\cal F}_t$ in the equality \eqref{DRBSDE}, we derive that 
$
Y=H- H'$, 
where $H$ and $H'$ are the two nonnegative strong supermartingales given by 
$$H_t := E[\xi_T^+ +\int_{t}^T f^+(s,Y_s,  Z_s,k_s)ds+A_T-A_t +C_{T-} -C_{t-} \,\vert \, {\cal F}_t ];$$
$$H'_t := E[\xi_T^- +\int_{t}^T f^-(s,Y_s,  Z_s,k_s)ds+A'_T-A'_t +C'_{T-} -C'_{t-} \,\vert \, {\cal F}_t ].$$
Since $Y=H-H'$ and $\xi \leq Y \leq \zeta$, we get $\xi \leq H-H'\leq \zeta$, which ensures that the  following  condition holds:
\begin{definition}[Mokobodzki's condition]
Let $(\xi,\zeta) \in{\cal S}^2\times{\cal S}^2 $  be a \emph{pair of admissible barriers}.  
 We say that the pair $(\xi,\zeta)$ satisfies \emph{Mokobodzki's condition}  if there exist two nonnegative strong supermartingales $H$ and $H'$ in $ \mathcal{S}^2$ such that:
\begin{equation}\label{Moki}
\xi_t  \leq H_t -H'_t \leq \zeta_t  \quad 0\leq t \leq T \quad {\rm a.s.}
\end{equation}
\end{definition}

The above reasoning gives us the following property. 
\begin{lemma}\label{Mok_is_necessary}
Mokobodzki's condition   is a necessary condition for the existence of a solution to the DRBSDE \eqref{def_solution_DRBSDE}.
\end{lemma}
We will see in Theorem \ref{exiuni} that it is also  a sufficient condition for the existence of a solution.



%
\subsection{The case when $f$ does not depend on the solution}\label{subsec3}
Let us now  investigate the question of existence and uniqueness of the solution to the DRBSDE defined above in the case where the driver $f$ does not depend on $y$, $z$, and $\mathpzc{k}$, that is, $f=(f_t)$, where $(f_t)$ is a process belonging to $\H^2$.

In this section, to simplify the notation, we suppose  that  the processes $\xi$ and $\zeta$ are left-limited. In this case, we can replace the process  
$(\overline{\xi}_t)$ (resp. $(\underline{\zeta}_t)_{t\in]0,T]}$) by $(\xi_{t-})_{t\in]0,T]}$ 
(resp. $(\zeta_{t-})_{t\in]0,T]}$) 
in the Skorokhod conditions \eqref{RBSDE_A} from Definition \ref{def_solution_DRBSDE}.

We stress that all the results of this section still hold true in the case where $\xi$ and $\zeta$ do not have left limits, provided we replace the process  $(\xi_{t-})_{t\in]0,T]}$ by the  process  $(\overline{\xi}_s)_{t\in]0,T]}$ and the process 
$(\zeta_{t-})_{t\in]0,T]}$ by the process $(\underline{\zeta}_t)_{t\in]0,T]}$.

\subsubsection{Equivalent formulation via  two (coupled) reflected BSDEs}
We will first assume that there exists a solution of the DRBSDE  associated with  driver $f(\omega,t)$.
We will show that (up to the process $(E[\xi_T+\int_t^Tf_sds | \cal{F}_t])$)  the first component  of this solution  can be written as the difference of two reflected BSDEs.

Suppose that $(Y, Z,k, h, A, C,A',C')\in {\cal S}^2 \times \H^2 \times \H^2_\nu \times {\cal M}^{2, \bot} \times({\cal S}^2)^2\times ({\cal S}^2)^2$ is  a  solution to the DRBSDE  associated with  driver $f(\omega,t)$ and with a pair of admissible  barriers  
$(\xi,\zeta)$.\\
Let $\tilde Y_t:= Y_t - E[\xi_T +\int_{t}^T f_sds \,\vert \, {\cal F}_t ]$, for all $t\in[0,T]$. 
From this definition, 
together with Eq. \eqref{DRBSDE}, we have 
\begin{equation} \label{dif}
-d\tilde Y_t= - Z_t dW_t-\int_E k_t(e) \tilde N(dt,de)-  dh_t +dA_t   -dA'_t + dC_{t-}- dC'_{t-},
\end{equation}
where the processes $A$, $C$, $A'$, $C'$, satisfy the Skorokhod conditions \eqref{RBSDE_A} and \eqref{RBSDE_C}. Formally, 
the process $\tilde Y_t= Y_t - E[\xi_T +\int_{t}^T f_sds \,\vert \, {\cal F}_t ]$ thus naturally appears as the  difference of the first coordinates of the solutions of two  reflected BSDEs with driver $0$, the first (resp. second) one admitting $A$ and $C$ (resp. $A'$ and $C'$) as non decreasing processes \footnote{The definition  of a solution of a  reflected BSDE (with driver $0$) is recalled in Proposition \ref{rexiuni} in the Appendix.}. 
We now precise this assertion, and specify the obstacles of these two reflected BSDEs in terms of the processes $A$, $C$, $A'$, $C'$, $\xi$, $\zeta$ and $f$. 
By \eqref{dif}, we get
\begin{equation}\label{tildey}
{\tilde Y}_t = X^{f}_t - X^{' f}_t \text{ for all } t\in[0,T] \text{ a.s., }
\end{equation}
where the processes $X^{f}$ and $X^{' f}$ are defined by 
\begin{equation}\label{JJprime}
X^{f}_t  := E[A_T-A_t +C_{T-} -C_{t-} \,\vert \, {\cal F}_t ] \text{ and }
 X^{' f}_t := E[A'_T-A'_t +C'_{T-} -C'_{t-} \,\vert \, {\cal F}_t ], \text{ for all } t\in[0,T].
 \end{equation}

\begin{remark} \label{X^{f}} Note that $X^{f}$ and $ X^{' f}$ are two nonnegative (right-u.s.c.) strong supermartingales in $\mathcal{S}^2$ such that  $X^{f}_T= X^{' f}_T =0$ a.s.
\end{remark}
By the   orthogonal decomposition property of martingales in 
${\cal M}^2$
(recalled in Lemma \ref{theoreme representation}), there exist $(\pi, l,h^1)$ , $(\pi', l',h^2)$ $\in \H^2 \times \H^2_\nu \times {\cal M}^{2, \bot}$ such that 
\begin{equation}\label{un}
X^{f}_t= -\int_{t}^T  \pi_s dW_s-\int_{t}^T \int_E l_s(e) \tilde N(ds,de)- (h^1_T-h^1_t)+A_T-A_t+C_{T-} -C_{t-};
\end{equation}
\begin{equation}\label{deux}
X^{' f}_t= -\int_{t}^T  \pi'_s dW_s-\int_{t}^T \int_E l'_s(e) \tilde N(ds,de)- (h^2_T-h^2_t)+A'_T-A'_t+C'_{T-} -C'_{t-}.
\end{equation}

We introduce the following optional processes:
\begin{equation}\label{defz}
\Tilde{\xi}_t^{f}:=\xi_t-E[\xi_T+\int_t^Tf_sds|\cal{F}_t], \;
\quad \Tilde{\zeta}_t^{f}:=\zeta_t-E[\zeta_T+\int_t^Tf_sds|\cal{F}_t], \quad 0 \leq t \leq T.
\end{equation}
\begin{remark}\label{tildexi}  
Note that $\Tilde{\xi}$ and $\Tilde{\zeta}$  satisfy  
$\Tilde{\xi}_T^{f} = \Tilde{\zeta}_T^{f} = 0 \; \text{ a.s. }$   
 We also have   $\Tilde{\xi}^{f}$ $\in$ ${\cal S}^2$ and $\Tilde{\zeta}^{f}$ $\in$ ${\cal S}^2$. Indeed,
$\vert \Tilde{\xi}_t^{f} \vert $ $\leq $ $\vert \xi_t \vert + E[U|\cal{F}_t]$, where 
$U:=\vert \xi_T \vert +\int_0^T \vert f_s \vert ds$.
Now, since $\xi$ $\in$ ${\cal S}^2$ and $f \in {\mathbb H}^2$, we have 
$U$ $\in$ $L^2$. Thus, by Doob's martingale inequalities in $L^2$, 
the martingale $(E[U \,|\, \cal{F}_t])$ belongs to ${\cal S}^2$, which implies that 
$\Tilde{\xi}^{f}$ $\in$ ${\cal S}^2$. Similarly, it can be shown that $ \Tilde{\zeta}^{f}$ $\in$ ${\cal S}^2$.
\end{remark}

 From $ \xi\leq Y \leq \zeta$ and  the definitions of $\tilde Y$, $\Tilde{\xi}_t^{f}$, $\Tilde{\zeta}_t^{f}$, we derive 
 $ \Tilde{\xi}^{f} \leq \tilde Y \leq \Tilde{\zeta}^{f}$; since ${\tilde Y}_t = X^{f}_t - X^{' f}_t$, we have 
 $X^{f}_t \geq X^{' f}_t + \Tilde{\xi}_t^{f}$ and $X^{' f}_t \geq X^{ f}_t - \Tilde{\zeta}_t^{f}$. \\
 Note that $Y- \xi = \tilde Y- \Tilde{\xi}^{f}=  X^{f} - X^{' f}- \Tilde{\xi}^{f}$.
 The  Skorokhod condition \eqref{RBSDE_C} satisfied by $C$ can thus be written:
$
\Delta C_{\tau}( X^{f}_\tau - X^{' f}_\tau - \Tilde{\xi}_\tau^{f})=0$ a.s. We also have 
$\{Y_{t-} > \xi_{t-}\}= \{  X^{f}_{t-} >X^{' f}_{t-} + \Tilde{\xi}_{t-}^{f}\}$.
Hence, the  Skorokhod condition \eqref{RBSDE_A} satisfied by $A$ can  be written:  $\int_0^T {\bf 1}_{ \{  X^{f}_{t-} >X^{' f}_{t-} + \Tilde{\xi}_{t-}^{f}\}} dA_t = 0$ a.s. 
It follows that $(X^{f}, \pi, l, h^1, A,C)$ is the solution of the reflected BSDE
 associated with driver $0$ and obstacle 
$(X^{' f} + \Tilde{\xi}^{f})\mathbb{I}_{[0,T)}$.\footnote{We note that this obstacle process is equal to $E[A'_T-A'_t +C'_{T-} -C'_{t-} \,\vert \, {\cal F}_t ] + \Tilde{\xi}_t^{f}$ if  $t<T$, and $0$ if $t=T$. Moreover, this process belongs to $\mathcal{S}^2$ (due to Remarks \ref{X^{f}} and \ref{tildexi}), and thus,  
  is an admissible obstacle for  RBSDEs.}\\
 By similar arguments we get that $(X^{' f}, \pi', l',h^2, A',C')$ is the solution of the reflected BSDE 
associated with driver $0$ and obstacle 
$(X^{ f} - \Tilde{\zeta}^{f})\mathbb{I}_{[0,T)}$.\footnote{We note that this obstacle process is equal to $E[A_T-A_t +C_{T-} -C_{t-} \,\vert \, {\cal F}_t ] - \Tilde{\zeta}_t^{f}$ if  $t<T$, and $0$ if $t=T$.}\\
 Hence, by equality \eqref{tildey}, we get that the process $(Y_t - E[\xi_T +\int_{t}^T f_sds \,\vert \, {\cal F}_t ])$ can be written as the difference of the solutions of two coupled reflected BSDEs. More precisely, the following result holds.
 \begin{lemma}
 Let $Y$ be the first component of a solution of the DRBSDE  with parameters $(f,\xi,\zeta)$ (where $f$ is a driver process). We then have
 $$
 Y_t= X^f_t-X^{'f }_t+E[\xi_T+\int_t^Tf_sds | \cal{F}_t], \; 0 \leq t \leq T, \text{ a.s. }\,,
$$
where the processes $X^f$ and $X^{' f} $ satisfy the following coupled system of reflected BSDEs:
\begin{equation}\label{systemf}
X^{f} = \, \Ref [ (X^{' f} + \Tilde\xi^{f})\mathbb{I}_{[0,T)}]; \quad
X^{' f} =\, \Ref[(X^{f}-\Tilde \zeta^{f})\mathbb{I}_{[0,T)}],
\end{equation}
where $\Ref $ is the operator induced by the RBSDE with driver $0$ (cf. Definition \ref{definitionRef} in the Appendix). 
\end{lemma}
By this lemma, we derive that the existence of a solution to the DRBSDE with parameters $(f,\xi,\zeta)$ (where $f$ is a driver process) implies the existence of a solution to the coupled system of RBSDEs \eqref{systemf}. 
We will see in the following proposition that the converse statement also holds true.  
\begin{proposition}\label{critereexistence}(Equivalent formulation)
The DRBSDE  associated with  driver process $f=(f_t) \in\H^{2}$ and with a pair of admissible barriers $(\xi,\zeta)$ has a solution  if and only if there exist two  processes  $X_{\cdot}\in{\cal S}^2$  and $X^{' }_{\cdot}\in{\cal S}^2$   
  satisfying the coupled system of RBSDEs: 
  \begin{equation}\label{system}
X = \, \Ref [ (X^{' } + \Tilde\xi^{f})\mathbb{I}_{[0,T)}]; \quad
X^{' } =\, \Ref[(X-\Tilde \zeta^{f})\mathbb{I}_{[0,T)}].
\end{equation}
In this case, the optional process $Y$ defined by
\begin{equation}\label{oY}
 Y_t:= X_t-X^{' }_t+E[\xi_T+\int_t^Tf_sds | \cal{F}_t], \; 0 \leq t \leq T, \text{ a.s. }
\end{equation} 
gives the first component of a solution to the DRBSDE.\\
 
 
 %
\end{proposition}

\dproof
The "only if part" of the first assertion has been proved above. Let us prove the "if part" of the first statement, together with the second statement. Let $X_{\cdot}\in{\cal S}^2$  and $X^{' }_{\cdot}\in{\cal S}^2$ be two processes   satisfying the coupled system \eqref{system}. 
Let $(\pi, l, h^1,A,C)$ (resp. $(\pi',l',h^2,A',C')$) be the vector of the remaining components of the solution to the RBSDE whose first component is $X$ (resp. whose first component is $X^{' }$). 
We note that equations \eqref{un} and \eqref{deux} hold for $X$ and $X^{'}$ (in place of $X^f$ and $X^{'f}$).  
We define the optional process  $Y$ as in   \eqref{oY}. 

 
 Since by assumption $X$ and $X^{' }$ belong to ${\cal S}^2$, it follows that $X$ and $X^{' }$ are real-valued, which implies that the process $Y$ is well- defined.
From \eqref{oY} and the property  $X_T=  X^{' }_T=0$ a.s., we get   $ Y_T = \xi_T$ a.s.\, 
 From the system \eqref{system} we get $X _t \geq X^{' }_t + \Tilde{\xi}^{f}_t$ and 
$X^{' }_t \geq X_t -  \Tilde{\zeta}^{f}_t$ for all $t \in [0,T]$ a.s.\,
 By using the definitions of $\Tilde{\xi}^{f}$, $\Tilde{\zeta}^{f}$ and $Y$, we derive that 
 $\xi_t \leq  Y_t \leq \zeta_t$  for all $t \in [0,T]$ a.s.


Moreover, the processes $A,C$ (resp. $A',C'$) satisfy the Skorokhod conditions for RBSDEs. More precisely,  for $A$ and $C$ we have: for all $\tau\in\T_0$, $\Delta C_{\tau}= {\bf 1}_{\{X_{\tau} = X^{' }_{\tau} + \Tilde{\xi}^{f}_\tau \}}\Delta C_{\tau}$ a.s.; for all predictable $\tau\in\T_0$, $\Delta A_{\tau}$ $=$ $ {\bf 1}_{\{X_{\tau^-}=X^{' }_{\tau^-}+\Tilde \xi^{f}_{\tau^-}\}}\Delta A_{\tau}$ a.s.; and $\int_0^T {\bf 1}_{\{X_{t}>X^{' }_{t}+\Tilde \xi^{f}_{t}\}} dA^c_t = 0$ a.s. Similar conditions hold for $A'$ and $C'$. 


Now, by using the definitions of $\Tilde{\xi}^{f}$ and $Y$, we get  $\{X_{\tau}=X^{' }_{\tau}+\Tilde \xi^{f}_{\tau}\}=\{ Y_{\tau}= \xi_{\tau}\}$, $\{X_{\tau-}=X^{' }_{\tau-}+\Tilde \xi^{f}_{\tau-}\}=\{ Y_{\tau-}= \xi_{\tau-}\}$ and $\{X_{t}>X^{' }_{t}+\Tilde \xi^{f}_{t}\}$ $=$ 
$\{Y_t > \xi_t\}$. Combining this with the previous observation gives 
$\Delta C_{\tau}= {\bf 1}_{\{ Y_\tau= \xi_\tau \} }\Delta C_{\tau}$ a.s. for all $\tau\in\T_0$ and $\int_0^T {\bf 1}_{\{Y_{t-} > \xi_{t-}\}} dA_t = 0$ a.s.
By applying the same arguments to $A'$ and $C'$, we get 
$\Delta C'_{\tau}= {\bf 1}_{\{Y_\tau= \zeta_\tau \} }\Delta C'_{\tau}$ a.s.\,for all $\tau \in \T_0$ and  
$\int_0^T {\bf 1}_{\{ Y_{t-} < \zeta_{t-}\}} dA'_t = 0$ a.s.\,\\
 We now note that the process $(E[\xi_T+\int_t^T f_sds | \cal{F}_t])_{t\in[0,T]}$ (which appears in the definition of $Y$) corresponds to the first component of the solution to the (non-reflected) BSDE with terminal condition $\xi_T$ and driver $f$. Hence,  there exist $\overline{\pi} \in \mathbb{H}^{2}$, $\overline{l} \in \mathbb{H}_{\nu}^{2}$  and $\overline{h}\in {\cal M}^{2, \bot}$ such that
$
E[\xi_T+\int_t^T f_sds | \cal{F}_t]=\xi_T+\int_t^T f_sds- \int_t^T \overline{\pi}  dW_s -\int_t^T \int_{ E}  \overline{l}_s(e) \tilde{N}(ds,de)- (\overline{h}_T-\overline{h}_t). 
$
 From this, together with the definition of $Y$ and equations \eqref{un} and \eqref{deux} for $X$ and $X^{' }$,  we obtain 
 $$  Y_t = \xi_T+ \int_t^T f_sds - \int_t^T Z_s  dW_s -\int_t^T \int_{ E}  k_s(e) \tilde{N}(ds,de)
 - (h_T-h_t)+\alpha_T-\alpha_t +\gamma_{T-}-\gamma_{t-},$$
where $Z:= \pi-\pi' + \overline{\pi} $, $k:= l-l' + \overline{l} $, 
$ h:= h^1-h^2 + \overline{h}$,
 $\alpha:= A-A^{'}$ and $\gamma:= C- C'$.\\
 If $dA_t \perp dA'_t$ and  $dC_t \perp dC'_t$, 
 then  
 $(Y,Z,k,h,A,C,A',C')$ is a solution to the doubly reflected BSDE with parameters $(f,\xi,\zeta)$, which gives the desired result.\\ 
Otherwise,
by the canonical decomposition of  RCLL predictable (resp. optional) processes with integrable variation 
(cf. Proposition A.7 in  \cite{DQS2}), there exist two nondecreasing right-continuous predictable (resp. optional) processes 
$B$ and $B'$ (resp. $D$ and $D'$) belonging to ${\cal S}^2$ such that $\alpha=B-B'$ (resp. $\gamma=D-D'$)
with $dB_t \perp dB'_t$ (resp. $dD_t \perp dD'_t$). Moreover, $dB_t \!<\!\!<\! dA_t$, $dB'_t \!<\!\!<\! dA'_t$,
$dD_t \!<\!\!<\! dC_t$ and $dD'_t \!<\!\!<\! dC'_t$. \\
Hence, 
since $\int_0^T \textbf{1}_{ \{Y_{t^-} > \xi_{t^-}\}}dA_t = 0$ a.s.\,, we get
 $ \int_0^T \textbf{1}_{\{Y_{t^-} > \xi_{t^-}\}}dB_t=0$ a.s.\, Similarly, we obtain $\int_0^T \textbf{1}_{\{ Y_{t^-} < \zeta_{t^-}\}}dB'_t=0$ a.s. 
 Moreover, since $dD_t \!<\!\!<\! dC_t$, the process $D$ is purely discontinuous and 
 $\Delta D_{\tau}= {\bf 1}_{\{  Y_\tau= \xi_\tau \} }\Delta D_{\tau}$ a.s.\,\,for all $\tau \in \T_0$. Similarly, $D'$ is purely discontinuous and $\Delta D'_{\tau}= {\bf 1}_{\{  Y_\tau= \zeta_\tau \} }\Delta D'_{\tau}$ a.s.\,for all $\tau \in \T_0$.
 The nondecreasing RCLL processes   $D,D'$ are thus purely discontinuous and  satisfy the Skorokhod condition  \eqref{RBSDE_C}. The nondecreasing RCLL processes   $B,B'$ satisfy the Skorokhod condition  \eqref{RBSDE_A}. 
The process
 $(Y,Z,k,h,B,D,B',D')$ is thus a solution to the DRBSDE with parameters $(f,\xi,\zeta)$.
\fproof


%
%

In the next section, we show that, under Mokobodzki's condition, there exists a solution of the coupled system of reflected BSDEs \eqref{system}, which, by  Proposition \ref{critereexistence}, will imply the existence of a solution of the doubly reflected BSDE associated with  driver process 
$f=(f_t)\in\H^2$.

\subsubsection{Existence of a solution of the coupled system of  RBSDEs. Existence of \\
 a solution of the DRBSDE with driver process  $(f_t)$}
Let $f=(f_t)\in\H^2$ be a driver process (as above). We  show the existence of a solution to the system \eqref{system} under Mokobodzki's condition. To do that, we use Picard's iterations. \\
We set ${\cal X}^0=0$ and ${\cal X}^{' 0}=0$,  and we define recursively, for each $n \in \mathbb{N}$, the   processes:
\begin{equation}\label{sys1}
{\cal X}^{ n+1}:=\Ref[({\cal X}^{' n}+\Tilde{\xi}^{f}){\bf 1}_{[0,T)}] \quad ; \quad {\cal X}^{' n+1}:=\Ref[({\cal X}^{ n}-\Tilde{\zeta}^{f}){\bf 1}_{[0,T)}]
\end{equation}
We see, by induction, that the processes ${\cal X}^{n}$ and ${\cal X}^{' n}$ are well-defined; moreover,  ${\cal X}^{n}$ and ${\cal X}^{' n}$ are strong supermartingales  in $ {\cal S}^2$. For the sake of simplicity, we have omitted the dependence on  $f$ in the notation  for  ${\cal X}^{ n}$ and ${\cal X}^{' n}$.



\begin{proposition}\label{seq}(Existence of a solution of the coupled system of RBSDEs)
Assume that the admissible pair $(\xi,\zeta)$ satisfies Mokobodzki's condition. 
The sequences of optional processes $({\cal X}_{\cdot}^{n})_{n \in \mathbb{N}}$ and 
$({\cal X}^{' n}_{\cdot})_{n \in \mathbb{N}}$ defined above are  nondecreasing.
The limit processes 
\begin{equation}\label{calX}
{\cal X}^{f} _{\cdot} := \lim_{n \rightarrow + \infty} {\cal X}_{\cdot}^{n} \,\,{\rm and} \,\,
{\cal X}^{' f} _{\cdot} := \lim_{n \rightarrow + \infty} {\cal X}_{\cdot}^{' n}
\end{equation}
satisfy the system \eqref{system} of coupled RBSDEs.\\  
Moreover, ${\cal X}^{f}_{\cdot}, {\cal X}^{' f}_{\cdot}$ are  the {\em smallest} processes in $\mathcal{S}^2$ satisfying  the system \eqref{system}.\\
The processes ${\cal X}^{f}, {\cal X}^{' f}$ are  also characterized as the  {\em minimal} nonnegative strong supermartingales  in  $\mathcal{S}^2$ satisfying the inequalities
$ \tilde \xi^{f} \leq {\cal X}^{f} -{\cal X}^{'f} \leq \tilde \zeta^{f}.$ 
\end{proposition}
The proof is given in the Appendix.

From this result together with Proposition \ref{critereexistence}, we derive the existence of a solution of the doubly reflected BSDE \eqref{DRBSDE} with driver process $(f_t)$.
\begin{corollary}(Existence of a solution of the DRBSDE)
Let $(\xi,\zeta)$ be an  admissible pair  satisfying Mokobodzki's condition. Then,  there exists a solution of the DRBSDE \eqref{DRBSDE} associated with driver process $f=(f_t)$, whose first component is given by
$Y_t= {\cal X}^{f}_t-{\cal X}^{'f}_t+E[\xi_T+\int_t^Tf_sds | \cal{F}_t], \; 0 \leq t \leq T, \text{ a.s. }$ (where ${\cal X}^{f}$, ${\cal X}^{'f}$ are the processes defined in \eqref{calX}).
\end{corollary}

In the following theorem we summarize some of the properties established so far.  
\begin{theorem}\label{Mokoko}
Let $f=(f_t) \in \H^2$ be a driver process.   Let $(\xi,\zeta)$ be a  pair of admissible barriers.  The following assertions are equivalent:
\begin{itemize}
\item[(i)] The pair $(\xi,\zeta)$ satisfies  Mokobodzki's condition.
\item[(ii)]  The system \eqref{system} of coupled RBSDEs  admits a solution.
\item[(iii)] The  DRBSDE \eqref{DRBSDE} with driver process $f=(f_t)$ has a  solution. 
\end{itemize}
\end{theorem}

\dproof The implication $(i)\Rightarrow (ii)$ follows from Proposition \ref{seq}. The equivalence between $(ii)$ and $(iii)$ has been established in Proposition \ref{critereexistence}. 
By Lemma \ref{Mok_is_necessary}, the implication $(iii)\Rightarrow (i)$ holds.
\fproof 

\begin{remark}We note that Mokobodzki's condition is satisfied if $(\xi, \zeta)$ is an admissible pair such that $\xi$ and/or $\zeta$ is an optional semimartingale of the form given in Lemma \ref{Moki_lem} in the Appendix.
\end{remark}

\subsubsection{Uniqueness of the solution of the DRBSDE with driver process $(f_t)$}
Let us now  investigate the question of uniqueness of the solution to the DRBSDE with  driver process  
$(f_t) \in\H^2$. 
To this purpose, we first state a lemma which will be used in the sequel.

Let $\beta>0$. For $\phi\in\H^2$, $\|\phi\|_\beta^2:= E[\int_0^T \e^{\beta s} \phi_s^2 ds ].$  For $l\in\H^2_\nu$, $\|l\|_{\nu,\beta}^2:=E[\int_0^T \e^{\beta s} \|l_s\|_\nu^2  ds].$ 
For $\phi\in{\cal S}^2$, we define $\vertiii{\phi}^2_\beta:=E[\esssup_{\tau \in \T_0}\e^{\beta \tau}\phi_\tau^2]$. We note that $\vertiii{\cdot}_\beta$ is a norm on ${\cal S}^{2}$  equivalent to the norm $\vertiii{\cdot}_{{\cal S}^{2}}$.   
For $M \in {\cal M}^{2}$,
$\| M\|^2_{\beta, {\cal M}^{2}}:=  E(\int_{]0,T]} \e^{\beta s} d[M] _s)$.

\begin{lemma}[A priori estimates]\label{Lemma_estimate}
Let $(Y, Z,k, h, A, C,A',C')\in {\cal S}^2 \times \H^2 \times \H^2_\nu \times {\cal M}^{2, \bot} \times({\cal S}^2)^2\times ({\cal S}^2)^2$ (resp. $(\bar Y, \bar Z, \bar k, \bar h, \bar A, \bar C, \bar A',\bar C')\in {\cal S}^2 \times \H^2 \times \H^2_\nu \times {\cal M}^{2, \bot} \times({\cal S}^2)^2\times ({\cal S}^2)^2$) be  a  solution to the DRBSDE  associated with  driver process  $f=(f_t) \in\H^2$ (resp. $\bar f= (\bar f_t) \in\H^2$) and with a pair of admissible  obstacles $(\xi,\zeta)$.  Then, there exists $c>0$ such that for all $\varepsilon>0$,  for all $\beta\geq\frac 1 {\varepsilon^2}$ we have 
\begin{eqnarray}\label{eq_initial_Lemma_estimate}
\|k-\bar k\|^2_{\nu,\beta} \, \leq  \, \varepsilon^2  \|f-\bar f\|^2_\beta\,;  & &\! \! \! \! \! \! \! \! \! \! 
\|Z-\bar Z\|^2_\beta \, \leq \, \varepsilon^2  \|f-\bar f\|^2_\beta \,;  \,\,\,\,  \|h-\bar h\|^2_{\beta,{\cal M}^2} \, \leq  \, \varepsilon^2  \|f-\bar f\|^2_\beta ; \nonumber\\
  \quad \vvertiii{Y-\bar Y}^2_\beta  &\leq& 4\varepsilon^2(1+6c^2)  \|f-\bar f\|^2_\beta.
\end{eqnarray}

\end{lemma}

The proof, which relies on Gal'chouk-Lenglart's formula (cf.  Corollary A.2 in \cite{MG}), is given in the Appendix.

We prove below \textit{the uniqueness} of
the solution to the DRBSDE associated with the driver process $(f_t)$ and with the admissible pair of  barriers $(\xi,\zeta)$ satisfying Mokobodzki's condition. 

\begin{theorem}\label{OY} 
Let $(\xi, \zeta)$ be an admissible pair of barriers satisfying Mokobodzki's condition. Let $f=(f_t) \in\H^2$ be a driver process. 
 There exists  a unique solution  to the DRBSDE \eqref{DRBSDE} 
 associated with  parameters  $(\xi,\zeta,f)$. 
 \end{theorem}

 \dproof
Theorem \ref{Mokoko} yields the existence of a solution. It remains to show the uniqueness. Let $(Y,Z,k,h,A,C,A',C')$ be a solution of the DRBSDE associated with the driver process $(f_t)$ and the barriers $\xi$ and $\zeta$. 
By the \emph{a priori estimates} (cf. Lemma \ref{Lemma_estimate}), we derive the uniqueness of $(Y,Z,k,h)$. By Remark \ref{sautsY}, we have $\Delta C_t=(Y_{t+}-Y_t)^-$ for all $t$ a.s. and $\Delta C'_t=(Y_{t+}-Y_t)^+ $ for all $t$ a.s.\,, which implies the uniqueness of the purely discontinuous processes $C$ and $C'$. Moreover, since $(Y,Z,k,h,A,C,A',C')$ satisfies the equation \eqref{DRBSDE}, it follows that the process $A-A'$ can be expressed in terms  of $Y,C,C'$, the integral of the driver process $(f_t)$ with respect to the Lebesgue measure, the martingale $h$, and the stochastic integrals of $Z$ and  $k$ with respect to $W$ and $\tilde N$, respectively, which yields the uniqueness of the finite variation process $A- A'$. Now, since $dA_t \perp dA'_t$, the nondecreasing processes $A$ and $A'$ correspond to the (unique) canonical decomposition of this finite variation process, which ends the proof. 
\fproof

Using the minimality property of  $({\cal X}^{f}, {\cal X}^{'f})$ (cf. Proposition \ref{seq}), together with the uniqueness property of the solution of the DRBSDE \eqref{DRBSDE} with driver process $f=(f_t)$ and Proposition \ref{critereexistence}, 
we show that ${\cal X}^{f}= X^{f}$ and $ {\cal X}^{'f}=X^{'f}$, where  the processes $X^{f}$ and $X^{'f}$ are defined by \eqref{JJprime} (in terms of the  solution of the DRBSDE). More precisely, we have the following result.


\begin{proposition}[Identification of ${\cal X}^{f}$ and $ {\cal X}^{'f}$] \label{Rk_identification} 
Let $(Y, Z,k, h, A, C,A',C')$  be  the  solution to the DRBSDE  associated with  driver process  $f=(f_t) \in\H^2$. Let ${\cal X}^{f}$ and $ {\cal X}^{'f}$ be the  strong supermartingales defined by \eqref{calX}.  We have a.s.
$${\cal X}^{f}_t  = E[A_T-A_t +C_{T-} -C_{t-} \,\vert \, {\cal F}_t ]\,\, \text{and}\,\,
 {\cal X}^{'f}_t = E[A'_T-A'_t +C'_{T-} -C'_{t-} \,\vert \, {\cal F}_t ], \text{ for all } t\in[0,T],$$
and
 $Y_t= {\cal X}^{f}_t-{\cal X}^{'f}_t+E[\xi_T+\int_t^Tf_sds | \cal{F}_t], \; 0 \leq t \leq T, \text{ a.s. }$ 
 \end{proposition}
 The proof is given in the Appendix.


\subsection{The case of a general Lipschitz driver $f(t,y,z,\mathpzc{k} )$}\label{subsec4}

We now prove the existence and uniqueness of the solution to the DRBSDE from  Definition  \ref{def_solution_DRBSDE} in the case of a general Lipschitz driver.

\begin{theorem}[Existence and uniqueness of the solution]\label{exiuni}
Let  $(\xi, \zeta)$  be  a pair of admissible barriers satisfying Mokobodzki's condition and let $f$ be a  Lipschitz driver. 
The DRBSDE with parameters $(f,\xi,\zeta)$ from Definition \ref{def_solution_DRBSDE} admits a unique solution $(Y,Z,k,h,A,C,A',C')\in \mathcal{S}^2  \times \H^2 \times \H^2_\nu \times {\cal M}^{2, \bot} \times(\mathcal{S}^2)^2\times (\mathcal{S}^2)^2.$

\end{theorem}
The proof, which relies on the estimates provided in Lemma \ref{Lemma_estimate} and a fixed point theorem, is given in the Appendix.

\section{Doubly reflected BSDEs with irregular barriers and $\mathbfcal{E}^{^{f}}$-Dynkin games with irregular rewards}\label{sec5}

The purpose of this section is to connect our DRBSDE  with irregular barriers  to a zero-sum  game problem between  two "stoppers" whose pay-offs are irregular and are assessed by  non-linear $f$-expectations. 

  In the "classical" case where $f\equiv 0$ (or, more generally, where $f$ is a given process $(f_t)$ $\in \mathbb{H}^2$), this topic has been first studied in \cite{CK} in the case of continuous barriers, and in \cite{Ham} and \cite{H} in the case of  right-continuous barriers. The case of right-continuous barriers and a general Lipschitz driver $f$  has been studied in \cite{DQS2}.

The following assumption holds in the sequel. 
 
\begin{Assumption}\label{Royer} 
Assume that  $dP \otimes dt$-a.s\, for each $(y,z, \mathpzc{k}_1,\mathpzc{k}_2)$ $\in$ $ \RB^2 \times (L^2_{\nu})^2$,
$$f( t,y,z, \mathpzc{k}_1)- f(t,y,z, \mathpzc{k}_2) \geq \langle \gamma_t^{y,z, \mathpzc{k}_1,\mathpzc{k}_2}  \,,\,\mathpzc{k}_1 - \mathpzc{k}_2 \rangle_\nu,$$ 
\begin{equation*}
\text{with } \;\; \gamma:  [0,T]  \times \Omega\times \RB^2 \times  (L^2_{\nu})^2  \rightarrow  L^2_{\nu}\,; \, (\omega, t, y,z, \mathpzc{k}_1, \mathpzc{k}_2) \mapsto 
\gamma_t^{y,z,\mathpzc{k}_1,\mathpzc{k}_2}(\omega,.)
\end{equation*}
 ${\cal P } Â·\otimes {\cal B}(\R^2) \otimes  {\cal B}( (L^2_{\nu})^2 )$-measurable and satisfying the inequalities
  \begin{equation}\label{condi}
\gamma_t^{y,z, \mathpzc{k}_1, \mathpzc{k}_2} (e)\geq -1\,\,\,  \;\; \text{ and }
\,\,  \;\;\|\gamma_t^{y,z, \mathpzc{k}_1, \mathpzc{k}_2}\|_{\nu}  \leq C,
\end{equation}
for each $(y,z, \mathpzc{k}_1, \mathpzc{k}_2)$ $\in$ $\RB^2 \times (L^2_{\nu})^2$, respectively $ dP\otimes dt \otimes d\nu(e)$-a.s.\, and 
$ dP\otimes dt$-a.s.\, (where $C$ is a positive constant).
\end{Assumption}
 Assumption \ref{Royer} ensures the non decreasing property of $\mathbfcal{E}^{f}$ by the comparison theorem for BSDEs with jumps (cf. Theorem 4.2 in \cite{16}). 

\subsection{The case where $\xi$ and $- \zeta$ are right upper-semicontinuous}\label{subsect_right}





 In this subsection we focus on the case where $\xi$ is \emph{right upper-semicontinuous} (right u.s.c.) and $\zeta$ is \emph{right lower-semicontinuous} (right l.s.c.). We  interpret the solution of our Doubly Reflected BSDE  in terms of the value process of a suitably defined zero-sum game problem on \emph{stopping times} with (non-linear) $f$-expectations. 
  This result generalizes the one shown in \cite{DQS3} in the case of RCLL payoffs and a filtration associated with the Brownian motion and a Poisson random measure.

 Let $\xi\in\mathcal{S}^2$ and $\zeta\in\mathcal{S}^2$. 
 We suppose that $\xi \leq \zeta$. 
We consider a game problem with two players where  each of the players' strategy is a stopping time in  $\T_0$ and the players payoffs are defined in terms of the given processes $\xi$ and $\zeta$. More precisely, if the first agent chooses $\tau\in\T_0$ as his/her strategy and the second agent chooses $\sigma\in\T_0$, then, at time $\tau\wedge\sigma$ (when the game ends),  the pay-off (or reward)  is $I(\tau,\sigma)$,   where 
\begin{equation} \label{terminal}
I(\tau,\sigma):=\xi_{\tau} {\bf 1} _{\tau \leq \sigma}+ \zeta_{\sigma} {\bf 1} _{\sigma < \tau}.
\end{equation}
The associated criterion (from time $0$ perspective) is defined as the $f$-evaluation of the pay-off, that is, 
by $\mathbfcal{E}^{^{f}} _{_{0,\tau \wedge \sigma}}[I(\tau, \sigma)]$. 
The first agent aims at choosing a stopping time $\tau\in\T_0$ which maximizes the criterion. 
The second  agent has the antagonistic objective of choosing a strategy $\sigma \in\T_0$ which minimizes the criterion.

As is usual in stochastic control, we embed the above (game) problem in a dynamic setting, by considering the game from time $\theta$ onwards, where $\theta$ runs through $\T_0$. From \emph{the   perspective of time $\theta$} (where $\theta\in\T_0$ is given), the first agent aims at choosing a strategy $\tau\in\T_{\theta} $ such that  $\mathbfcal{E}^{^{f}} _{_{\theta,\tau \wedge \sigma}}[I(\tau, \sigma)]$ be maximal. The second  agent has the antagonistic objective of choosing $\sigma\in\T_{\theta}$ such that $\mathbfcal{E}^{^{f}} _{_{\theta,\tau \wedge \sigma}}[I(\tau, \sigma)]$ be  minimal. 


The following notions will be used in the sequel: 
\begin{definition} \label{Def_values}
Let  $\theta \in {\cal T}_0$.
\begin{itemize}
\item   {\em The upper value} $\overline{V}(\theta)$ and {\em the  lower value} $\underline{V}(\theta)$ of the game at time $\theta$ are the random variables  defined respectively by 
 \begin{equation} \label{dessus}
\overline{V}(\theta):=\,\essinf_{\sigma \in \T_\theta }\,\esssup_{\tau \in \T_\theta} \, \mathbfcal{E}^{^{f}}_{_{\theta,\tau \wedge \sigma}}[I(\tau, \sigma)]; \quad 
\underline{V}(\theta):=\,\esssup_{\tau \in \T_{\theta}} \, \essinf_{\sigma \in \T_\theta} \, \mathbfcal{E}^{^{f}}_{_{\theta,\tau \wedge \sigma}}[I(\tau, \sigma)].
\end{equation}
\item 
We say that  {\em there exists a value} for the game 
 at time $\theta$  if $\overline{V}(\theta)=\underline{V}(\theta)$ a.s.
\item 
A pair $(\hat \tau, \hat \sigma) \in \T_{\theta}^2$ is called a \emph{ saddle point  at time $\theta$} for the  game if for all $(\tau, \sigma) \in \T_{\theta}^2$ we have
$$  \mathbfcal{E}^{^{f}}_{_{\theta, \tau \wedge \hat \sigma}}[I(\tau, \hat \sigma)] \leq  \mathbfcal{E}^{^{f}}_{_{\theta, \hat \tau \wedge \hat \sigma}}[I(\hat \tau , \hat \sigma)] \leq  \mathbfcal{E}^{^{f}}_{_{\theta, \hat \tau  \wedge \sigma}}[I(\hat \tau , \sigma)]  \quad \text{ a.s. }$$
\item Let $\varepsilon>0$.  A pair $(\hat \tau, \hat \sigma) \in \T_{\theta}^2$ is called an \emph{ $\varepsilon$-saddle point  at time $\theta$} for the  game if for all $(\tau, \sigma) \in \T_{\theta}^2$ we have
$$  \mathbfcal{E}^{^{f}}_{_{\theta, \tau \wedge \hat \sigma}}[I(\tau, \hat \sigma)]-\varepsilon \leq  \mathbfcal{E}^{^{f}}_{_{\theta, \hat \tau \wedge \hat \sigma}}[I(\hat \tau , \hat \sigma)] \leq  \mathbfcal{E}^{^{f}}_{_{\theta, \hat \tau  \wedge \sigma}}[I(\hat \tau , \sigma)]+\varepsilon  \quad \text{ a.s. }$$
\end{itemize}
\end{definition}
The inequality  $\underline{V}(\theta) \leq \overline{V}(\theta)$ a.s. is trivially true. 
As mentioned in the introduction, in the case where the processes $\xi$ and $\zeta$ are RCLL,  we recover a game problem which appears in \cite{DQS2} under the name of \textit{generalized Dynkin game}. 
In the case $f=0$, we have $\mathbfcal{E}^{^0}_{_{\theta,\tau \wedge \sigma}}[I(\tau, \sigma)] = 
 E\, [I(\tau, \sigma)\, \vert \, {\cal F}_{\theta}],$ and, in this case,  our game problem corresponds to  the classical Dynkin game (cf., e.g., \cite{ALM}).

We also recall the following definition:
\begin{definition}\label{defmart}
Let $Y$ 
 $ \in \cal S^2$. The process $Y$ is said to be a strong $\mathbfcal{E}^{^{f}}$-supermartingale (resp $\mathbfcal{E}^{^{f}}$-submartingale), if $\mathbfcal{E}^{^{f}}_{_{\sigma ,\tau}}[Y_{\tau}] \leq Y_{\sigma}$ (resp. $\mathbfcal{E}^{^{f}}_{\sigma ,\tau}[Y_{\tau}]\geq Y_{\sigma})$ a.s.\, on $\sigma \leq \tau$,  for all $ \sigma, \tau \in  \T_0$. 
\end{definition}
\begin{remark} Recall that $Y$ is right u.s.c.(cf. e.g. Lemma 5.1 in \cite{MG}).
\end{remark}
Let $Y$ be the first component of the solution to the DRBSDE with parameters $(f,\xi,\zeta)$ from Definition \ref{def_solution_DRBSDE}. 
For each $\theta$ $\in$ $ \T_0$ and each $\varepsilon >0$, we define the stopping times 
$\tau^{\varepsilon}_{\theta}$ and $\sigma^{\varepsilon}_{\theta}$ by
\begin{equation}\label{tauepsilon}
\tau^{\varepsilon}_{\theta} := \inf \{ t \geq \theta,\,\, Y_t \leq \xi_t + \varepsilon\}; \quad \sigma^{\varepsilon}_{\theta} := \inf \{ t \geq \theta,\,\, Y_t \geq \zeta_t - \varepsilon\}.
\end{equation}
\begin{lemma}\label{lalaun} 
The process
$(Y_t, \,\theta \leq t \leq
\tau^{\varepsilon}_{\theta} )$ is a strong $\mathbfcal{E}^{^{f}}$-submartingale and the process $(Y_t, \,\theta \leq t \leq \sigma^{\varepsilon}_{\theta} )$ is a strong $\mathbfcal{E}^{^{f}}$-supermartingale.
\end{lemma}

%

\dproof
Let us first prove that the process
$(Y_t, \,\theta \leq t \leq
\tau^{\varepsilon}_{\theta} )$ is a strong $\mathbfcal{E}^{^{f}}$-submartingale. By definition of $\tau^{\varepsilon}_{\theta}$, we have $Y_t> \xi_t  + \varepsilon$ on $[\theta, \tau^{\varepsilon}_{\theta}[$ a.s.; hence, $Y_{t-}\geq \overline{\xi}_t+\varepsilon$ on $[\theta, \tau^{\varepsilon}_{\theta}[$ a.s. 
Therefore, $A^c$ is constant on $[\theta, \tau^{\varepsilon}_{\theta}[$ a.s.\, (cf. Skorokhod conditions); by continuity of  the process $A^c$, $A^c$ is constant on the closed interval $[\theta, 
\tau^{\varepsilon}_{\theta}]$, a.s. Also, the process  $A^d$ is constant on $[\theta, \tau^{\varepsilon}_{\theta}[, $ a.s. (cf. Skorokhod conditions). 
Moreover, $ Y_{ (\tau^{\varepsilon}_{\theta})^-  } \geq \overline{\xi}_{ \tau^{\varepsilon}_{\theta}  }+ \varepsilon\,$ a.s.\,, which implies that $\Delta A^d _ {\tau^{\varepsilon}_{\theta}  } =0$ a.s.\, Finally, for a.e. $\omega\in\Omega $, for all $t\in[\theta(\omega), \tau^{\varepsilon}_{\theta}(\omega)[$, $\Delta C_t(\omega)=C_t(\omega)-C_{t-}(\omega)=0$; we deduce that for a.e. $\omega\in\Omega $, $(C_{t-}(\omega))$ is constant on $[\theta(\omega), \tau^{\varepsilon}_{\theta}(\omega)[$, and even on the closed interval $[\theta(\omega), \tau^{\varepsilon}_{\theta}(\omega)]$, since the trajectories of $(C_{t-})$ are left-continuous.  Thus,  the process $(A_t + C_{t-})$ is constant on $[\theta, \tau^{\varepsilon}_{\theta}]$ a.s.\\
Hence, $Y$ satisfies on $[\theta, 
\tau^{\varepsilon}_{\theta}]$ the following dynamics:
\begin{equation*} 
-d Y_t= f(t,Y_t, Z_t, k_t)dt - Z_t dW_t-\int_E k_t(e) \tilde N(dt,de)-  dh_t  -dA'_t - dC'_{t-}.
\end{equation*}
 By Lemma 12.2 in \cite{MG3}, the process
$(Y_t, \,\theta \leq t \leq
\tau^{\varepsilon}_{\theta} )$ is thus a strong $\mathbfcal{E}^{^{f}}$-submartingale. By similar arguments, we  can show that $(Y_t, \,\theta \leq t \leq \sigma^{\varepsilon}_{\theta} )$ is a strong $\mathbfcal{E}^{^{f}}$-supermartingale, which ends the proof of the lemma. 
\fproof

\noindent We now prove the following result under additional regularity assumptions on $\xi$ and $\zeta$. 

\begin{lemma}\label{laladeux} 
Suppose that $\xi$ is right upper-semicontinuous (resp. $\zeta$ is right lower-semicontinuous). We then have
\begin{equation}\label{tau}
 Y_{\tau^{\varepsilon}_{\theta}} \, \leq \, \xi_{\tau^{\varepsilon}_{\theta}} + \varepsilon  \quad ({\rm resp.} \quad Y_{\sigma^{\varepsilon}_{\theta}}\, \geq \, \zeta_{\sigma^{\varepsilon}_{\theta}} - \varepsilon) \quad \mbox{a.s.}
\end{equation}
\end{lemma}

\begin{remark} 
We stress that, contrary to Lemma  \ref{lalaun}, the result of the above Lemma \ref{laladeux} does not generally hold when $\xi$ and $\zeta$ are completely irregular.
\end{remark}

\dproof Suppose that $\xi$ is right u.s.c.
Let us  prove that $Y_{\tau^{\varepsilon}_{\theta}} \, \leq \, \xi_{\tau^{\varepsilon}_{\theta}} + \varepsilon$ a.s. By way of contradiction, we suppose $P(Y_{\tau^{\varepsilon}_{\theta}}> \xi_{\tau^{\varepsilon}_{\theta}} + \varepsilon)>0$. By the Skorokhod conditions, we have  $\Delta C_{\tau^{\varepsilon}_{\theta}}=C_{\tau^{\varepsilon}_{\theta}}- C_{(\tau^{\varepsilon}_{\theta})-}=0$ on the set $\{Y_{\tau^{\varepsilon}_{\theta}}> \xi_{\tau^{\varepsilon}_{\theta}} + \varepsilon\}$. On the other hand, due to Remark \ref{sautsY}, $\Delta C_{\tau^{\varepsilon}_{\theta}}=Y_{\tau^{\varepsilon}_{\theta}}-Y_{(\tau^{\varepsilon}_{\theta})+}.$ Thus, $Y_{\tau^{\varepsilon}_{\theta}}=Y_{(\tau^{\varepsilon}_{\theta})+}$ on the set $\{Y_{\tau^{\varepsilon}_{\theta}}> \xi_{\tau^{\varepsilon}_{\theta}} + \varepsilon\}.$ Hence,
\begin{equation}\label{eq_contradiction}
Y_{(\tau^{\varepsilon}_{\theta})+}> \xi_{\tau^{\varepsilon}_{\theta}} + \varepsilon \text{ on the set } \{Y_{\tau^{\varepsilon}_{\theta}}> \xi_{\tau^{\varepsilon}_{\theta}} + \varepsilon\}.
\end{equation}
 We will obtain a contradiction with this statement.
Let us fix $\omega\in\Omega$. By definition of $\tau^{\varepsilon}_{\theta}(\omega)$, there exists a non-increasing sequence $(t_n)=(t_n(\omega))\downarrow \tau^{\varepsilon}_{\theta}(\omega) $ such that $Y_{t_n}(\omega)\leq \xi_{t_n}(\omega) + \varepsilon,$ for all $n\in\N$. Hence,
$\limsupn Y_{t_n}(\omega)\leq \limsupn \xi_{t_n}(\omega)+\varepsilon.$ As, by assumption, the process $\xi$ is right-u.s.c.\,, we have $\limsupn \xi_{t_n}(\omega)\leq \xi_{\tau^{\varepsilon}_{\theta}}(\omega)$. On the other hand, as $(t_n(\omega))\downarrow \tau^{\varepsilon}_{\theta}(\omega)$, we have $\limsupn Y_{t_n}(\omega)=Y_{(\tau^{\varepsilon}_{\theta})+}(\omega).$ Thus, $Y_{(\tau^{\varepsilon}_{\theta})+}(\omega)\leq \xi_{\tau^{\varepsilon}_{\theta}}(\omega)+\varepsilon,$ which is in contradiction with \eqref{eq_contradiction}. We conclude that $Y_{\tau^{\varepsilon}_{\theta}}\leq  \xi_{\tau^{\varepsilon}_{\theta}} + \varepsilon \text{ a.s.}$ By similar arguments, one can show that if $\zeta$ is right-l.s.c.\,, then $Y_{\sigma^{\varepsilon}_{\theta}}\, \geq \, \zeta_{\sigma^{\varepsilon}_{\theta}} - \varepsilon$ a.s. 
The proof 
of the lemma is  thus complete.
\fproof

Using the above lemmas, we show that, under Mokobodzki's condition, the  game problem defined above has a value. Moreover, we characterize  the value of the game in terms of  the (first component of the) solution of the DRBSDE (\ref{DRBSDE}), and we also show  the existence of 
$\varepsilon$-saddle points.

\begin{theorem}[Existence and characterization of the value function]\label{caracterisation}

Let $f$ be a Lipschitz  driver satisfying Assumption $\eqref{Royer}$. Let $(\xi, \zeta)$ be an admissible  pair of barriers satisfying Mokobodzki's condition, and such that $\xi$ is right u.s.c.and $\zeta$ is right l.s.c.
Let $(Y,Z,k,h,A,A', C,C')$ be the solution of the DRBSDE (\ref{DRBSDE}). 
There exists a common value function for the $\mathbfcal{E}^f$-Dynkin game \eqref{dessus}, and for each stopping time $\theta$ $\in$ $\T_0$, we have
 \begin{equation}\label{car}
Y_{\theta} = \overline{V}(\theta)=\underline{V}(\theta) \quad \mbox{a.s.}
\end{equation}

Let $\theta$ $\in$ $\T_0$ and let $\varepsilon >0$. For each $(\tau, \sigma) \in \T_{\theta}^2$, the stopping times $\tau^{\varepsilon}_{\theta}$ and $\sigma^{\varepsilon}_{\theta}$, defined by \eqref{tauepsilon}, satisfy the  inequalities:
\begin{equation}\label{fifi}
\mathbfcal{E}^{^{f}}_{_{\theta, \tau \wedge \sigma^{\varepsilon}_{\theta}}} [I(\tau,  \sigma^{\varepsilon}_{\theta})]  - L\varepsilon  \,\,\leq \,\,  Y_{\theta} \,\,\leq \,\,  \mathbfcal{E}^{^{f}}_{_{\theta, \tau^{\varepsilon}_{\theta} \wedge \sigma}}[I(\tau^{\varepsilon}_{\theta},  \sigma)]  + L\varepsilon \quad \mbox{a.s.}\,,
\end{equation}
where $L$ is a positive constant which only depends on  the Lipschitz constant $K$ of $f$ and on the terminal time $T$.
In other terms, the pair $(\tau^{\varepsilon}_{\theta}, \sigma^{\varepsilon}_{\theta})$ is an $L \varepsilon$-saddle point at time $\theta$ for the $\mathbfcal{E}^f$-Dynkin game \eqref{dessus}.

%
\end{theorem}

\dproof 
 Let $\theta$ $\in$ $\T_0$ and let $\varepsilon >0$. Let us show that 
$(\tau^{\varepsilon}_{\theta},\sigma^{\varepsilon}_{\theta})$ satisfies the inequalities \eqref{fifi}.
By Lemma \ref{lalaun}, the process $(Y_t, \,\theta \leq t \leq \tau^{\varepsilon}_{\theta} )$ is a strong $\mathbfcal{E}^{^{f}}$-submartingale. We thus get
\begin{equation}\label{lala1}
Y_{\theta} \, \leq \,  \mathbfcal{E}^{^{f}}_{_{\theta, \tau^{\varepsilon}_{\theta} \wedge \sigma}}[Y_{\tau^{\varepsilon}_{\theta} \wedge \sigma}] \quad \mbox{a.s.}\
\end{equation}
Now, by assumption, $\xi$ is right-u.s.c.\, By Lemma \ref{laladeux}, we thus get  $ Y_{\tau^{\varepsilon}_{\theta}} \, \leq \, \xi_{\tau^{\varepsilon}_{\theta}} + \varepsilon$ a.s. Using the inequality $Y \leq \zeta$, we derive  
\begin{equation*}
Y_{\tau^{\varepsilon}_{\theta} \wedge \sigma} \, \leq \, (\xi_{\tau^{\varepsilon}_{\theta}} + \varepsilon) {\bf 1} _{\tau^{\varepsilon}_{\theta}  \leq \sigma}+ \zeta_\sigma {\bf 1} _{\sigma < \tau^{\varepsilon}_{\theta} } \, \leq \,  I(\tau^{\varepsilon}_{\theta} , \sigma) +\varepsilon \quad {\rm a.s.}\,,
\end{equation*}
where the last inequality follows from the definition of $I(\tau^{\varepsilon}_{\theta} , \sigma)$.
By using the inequality $\eqref{lala1}$ and the monotonicity of $\mathbfcal{E}^{^{f}}$, we  get 
\begin{equation}\label{22}
Y_{\theta} \, \leq \, \mathbfcal{E}^{^{f}}_{_{\theta, \tau^{\varepsilon}_{\theta} \wedge \sigma}}[I(\tau^{\varepsilon}_{\theta} , \sigma) +\varepsilon] \leq \mathbfcal{E}^{^{f}}_{_{\theta, \tau^{\varepsilon}_{\theta} \wedge \sigma}}[I(\tau^{\varepsilon}_{\theta} , \sigma)] + L \varepsilon\,\quad \mbox{a.s.}\,,
\end{equation}
where the last inequality follows from an estimate on BSDEs (cf. Proposition A.4 in \cite{16}).
%
By Lemma \ref{lalaun}, the process $(Y_t, \,\theta \leq t \leq \sigma^{\varepsilon}_{\theta} )$ is a strong $\mathbfcal{E}^{^{f}}$-supermartingale. We thus get
\begin{equation}\label{ba}
 Y_{\theta}\,\, \geq \,\,\mathbfcal{E}^{^{f}}_{_{\theta, \tau \wedge \sigma^{\varepsilon}_{\theta}}}[ Y_{\tau \wedge \sigma^{\varepsilon}_{\theta}}]
 \quad
 {\rm a.s.}
 \end{equation}
 Now, by assumption, $\zeta$ is right lower-semicontinuous. Hence, by Lemma \ref{laladeux}, we have
 $Y_{\sigma^{\varepsilon}_{\theta}}\, \geq \, \zeta_{\sigma^{\varepsilon}_{\theta}} - \varepsilon$ a.s. 
 Using similar arguments as above, we derive that 
$
 Y_{\theta}\,\, \geq \,\,\mathbfcal{E}^{^{f}}_{_{\theta, \tau \wedge \sigma^{\varepsilon}_{\theta}}}[I(\tau, \sigma^{\varepsilon}_{\theta})]  - L\varepsilon$ a.s\,, which, together with \eqref{22}, leads to the desired inequalities \eqref{fifi}.

 Now, since inequality \eqref{22} holds for all $\sigma \in \T_{\theta}$, it follows that 
 $$\,Y_{\theta} \leq \, \essinf_{\sigma \in \T_{\theta}} \mathbfcal{E}^{^{f}}_{_{\theta, \tau^{\varepsilon}_{\theta} \wedge \sigma}}[I(\tau^{\varepsilon}_{\theta},  \sigma)]+L\varepsilon \leq \, \,\esssup_{\tau \in \T_{\theta}} \, \essinf_{\sigma \in \T_\theta} \, \mathbfcal{E}^{^{f}}_{_{\theta,\tau \wedge \sigma}}[I(\tau, \sigma)] +L\varepsilon 
 \,\,\text{ a.s.}\,$$
 From this, together with  the definition of  $\underline{V}(\theta)$ (cf. 
\eqref{dessus}), we obtain  $Y_{\theta} \,\,\leq \,\, \underline{V}(\theta)+ L \varepsilon$ a.s. Similarly, we  show that $\overline{V}(\theta)- L \varepsilon \,\, \leq \,\, Y_{\theta}$ a.s. for all $\varepsilon>0$. We thus get
$\overline{V}(\theta) \,\, \leq \,\, Y_{\theta} \,\,\leq \,\, \underline{V}(\theta)$ a.s.\,  This, together with the inequality $\underline{V}(\theta) \leq \overline{V}(\theta)$ a.s.\,,  yields 
$\underline{V}(\theta)= Y_{\theta} = \overline{V}(\theta) \text{ a.s.}$
\fproof

We will now show the existence of saddle points under an additional regularity assumption on the barriers. 
Let $(Y,Z,k,h,A,A', C,C')$ be the solution of the DRBSDE (\ref{DRBSDE}). 
For each $\theta$ $\in$ $\T_0$, we introduce the following stopping times:
\begin{equation}\label{tauetoile}
\tau^*_{\theta}:= \inf \{ t \geq \theta,\,\, Y_t = \xi_t\} ;  \quad \, \sigma^*_{\theta}:= \inf \{ t  \geq \theta,\,\, Y_t = \zeta_t\},
\end{equation}
and
\begin{equation}\label{taubar}
\overline\tau_{\theta}:= \inf \{ t  \geq \theta,\,\, A_t > A_{\theta} \,\,{\rm or}\,\, C_{t^-}>C_{\theta^-}\} ;  \quad \overline\sigma_{\theta}:= \inf \{ t  \geq \theta,\,\, A'_t > A'_{\theta}\,\,{\rm or}\,\, C'_{t^-}>C'_{\theta^-}\}.
\end{equation}

\begin{remark}\label{etoilebar} 
Note that $\tau^*_{\theta} \leq \overline\tau_{\theta}$ a.s. and $\sigma^*_{\theta} \leq \overline\sigma_{\theta}$ a.s. Indeed, by definition of $\tau^*_{\theta}$, we have $Y_t >\xi_t$ on  $[\theta, \tau^*_{\theta}[$ a.s.\,; hence, by the Skorokhod condition satisfied by $A$ (resp. $C$), the process $A$ (resp. $C_{t^-}$) is constant on $[\theta, \tau^*_{\theta}[$. By definition of $\overline\tau_{\theta}$, it thus follows that $\tau^*_{\theta} \leq \overline\tau_{\theta}$ a.s. By similar arguments, we show that $\sigma^*_{\theta} \leq \overline\sigma_{\theta}$ a.s.
%
\end{remark}

We now prove the following result, which will be used to obtain the existence of saddle points under additional regularity assumptions on $\xi$ and $\zeta$ (cf. Corollary \ref{optimal}), as well as to study the pricing and (super)hedging of game options with irregular payoffs in Section \ref{gameoptions} (cf. Proposition \ref{superhedging}).

  \begin{proposition}\label{lalabis} Let $f$ be a driver satisfying Assumption $\eqref{Royer}$. Let $(\xi, \zeta)$ be an admissible  pair of barriers  satisfying Mokobodzki's condition and such that $\xi$ is right-u.s.c and and $\zeta$ is right l.s.c. Let $(Y,Z,k,h,A,A', C,C')$ be the solution of the DRBSDE (\ref{DRBSDE}). 
 For each $\theta \in \T_0$, the following assertions hold:
 \begin{enumerate}
 \item
Assume that  $A'$  is \emph{continuous} (which holds if, for example, $\zeta$ is left l.s.c.along stopping times \footnote{see Proposition \ref{sautsYA} (iii).}). Then, the process 
$(Y_t, \,\theta \leq t \leq \overline\sigma_{\theta} )$ is a strong $\mathbfcal{E}^{^{f}}$-supermartingale.\footnote{Note that since $\sigma^*_{\theta} \leq \overline\sigma_{\theta}$ a.s.\,, this implies that
$(Y_t, \,\theta \leq t \leq
\sigma^*_{\theta})$ is a strong $\mathbfcal{E}^{^{f}}$-supermartingale. } Moreover, we have
\begin{equation}\label{eg2}
  Y_{\sigma^*_{\theta}}\, = \, \zeta_{\sigma^*_{\theta}}   \quad  {\rm and} \quad   Y_{\overline\sigma_{\theta}}\, = \, \zeta_{\overline\sigma_{\theta}} \quad \mbox{a.s.}
\end{equation}


\item 
Assume that  $A$  is \emph{continuous} (which holds if, for example, $ \xi $ is left u.s.c. along stopping times). Then, the process $(Y_t, \,\theta \leq t \leq
\overline\tau_{\theta})$ is a strong $\mathbfcal{E}^{^{f}}$-submartingale.
Moreover, we have
\begin{equation}\label{eg1}
 Y_{\tau^*_{\theta}} \, = \, \xi_{\tau^*_{\theta}} \quad  {\rm and} \quad Y_{\overline\tau_{\theta}} \, = \, \xi_{\overline\tau_{\theta}}\quad  \mbox{a.s.}
\end{equation}
\end{enumerate}


\end{proposition}

%

 \dproof 
We suppose $A'$ is continuous.  
By definition of $\overline{\sigma}_{\theta}$,  we have $A'_{\overline{\sigma}_{\theta}}=A'_{\theta}$ a.s. and 
 $C'_{\overline{\sigma}_{\theta^-}}=C'_{\theta^-}$  a.s. 
 because $(A'_t)$ and $(C'_{t^-})$ are left-continuous.\\ 
Hence, $Y$ satisfies on  $[\theta, 
\overline{\sigma}_{\theta}]$ the following dynamics:
\begin{equation*} 
-d Y_t= f(t,Y_t, Z_t, k_t)dt - Z_t dW_t-\int_E k_t(e) \tilde N(dt,de)-  dh_t  +dA_t +dC_{t-}.
\end{equation*}
 By Lemma 12.2 in \cite{MG3},
  the process 
$(Y_t, \,\theta \leq t \leq
\overline\sigma_{\theta} )$ is thus a strong $\mathbfcal{E}^{^{f}}$-supermartingale. 

 Moreover, by definition of $\overline{\sigma}_{\theta}$, since the continuous process $A'$ increases only on $\{Y_t=\zeta_t\}$ and $\Delta C'_t= 
 {\bf 1}_{\{ Y_t = \zeta_t\}}\Delta C'_t$, we get $Y_{\overline{\sigma}_{\theta}}=\zeta_{\overline{\sigma}_{\theta}}$ a.s.\\
 It remains to show the equality $ Y_{\sigma^*_{\theta}} = \zeta_{\sigma^*_{\theta}}$ a.s. Note first that $ Y_{\sigma^*_{\theta}} \leq \zeta_{\sigma^*_{\theta}}$ a.s., since $Y$ is (the first component of) the solution to the DRBSDE with barriers $\xi$ and $\zeta$. We show that $ Y_{\sigma^*_{\theta}} \geq \zeta_{\sigma^*_{\theta}}$ a.s.  by using the assumption of right-lower semicontinuity on the process $\zeta$; the  arguments are similar to those used in the proof of Lemma \ref{laladeux} and are left to the reader.     
The proof of the first assertion is thus complete.\\
The case where $A$ is continuous can be treated by similar arguments.   
\fproof

Using the above proposition, Proposition \ref{sautsYA} (iii) and Theorem \ref{caracterisation}, we derive the existence of saddle points when $ \xi $ is left u.s.c.and $\zeta$ is left l.s.c.along stopping times.

\begin{corollary}\label{optimal}
Let the assumptions of  
Theorem \ref{caracterisation} hold.  We assume moreover that 
 $ \xi $ is left u.s.c. along stopping times and $\zeta$ is left l.s.c. along stopping times.
Then, for each $\theta \in \T_0$, the pairs of stopping times $(\tau_{\theta}^*, \sigma_{\theta}^*)$ and $(\overline{\tau}_{\theta}, \overline{\sigma}_{\theta})$, defined by \eqref{tauetoile} and  \eqref{taubar},
are saddle points at time  $\theta$ for the $\mathbfcal{E}^f$-Dynkin game.
%
\end{corollary}

\dproof
The proof of this corollary is given in the Appendix. 
\fproof


\paragraph{Classical Dynkin game with irregular rewards}
In this paragraph, we consider the particular case where $f\equiv 0$, that is, the case where the $f$-expectation  reduces to the classical linear expectation.
Let $(\xi, \zeta)$ be an admissible  pair of barriers.   
Let $\theta \in  \T_0$. For  $\tau \in \T_{\theta}$  and $\sigma \in \T_{\theta} $, it holds $\mathbfcal{E}^{^0} _{_{\theta,\tau \wedge \sigma}}[I(\tau, \sigma)] = 
 E\, [I(\tau, \sigma)\, \vert \, {\cal F}_{\theta}].$ 
 The {\em upper} and {\em lower values} at time $\theta$ are then given by 
 \begin{equation} \label{dessusclassical}
\overline{V}(\theta):=\,\essinf_{\sigma \in \T_\theta }\,\esssup_{\tau \in \T_\theta} \,  E\,[I(\tau, \sigma)\, \vert \, {\cal F}_{\theta}]; \quad 
\underline{V}(\theta):=\,\esssup_{\tau \in \T_{\theta}} \, \essinf_{\sigma \in \T_\theta} \,  E\,[I(\tau, \sigma)\, \vert \, {\cal F}_{\theta}],
\end{equation}
We thus  recover the classical Dynkin game on stopping times (with linear expectations) recalled in the introduction (cf., e.g., \cite{Bismuth} and \cite{ALM}). 
Recall that in \cite{ALM}, it has been shown that, if $(\xi, \zeta)$ satisfies Mokobodzki's condition, $\xi$ is right u.s.c.and $\zeta$ is right l.s.c.\,, then the  Dynkin game \eqref{dessusclassical} has a value (cf. Section 4.1 in \cite{ALM}).  
Thanks to our Theorem \ref{caracterisation}, we recover this result.

Moreover, we provide an  infinitesimal characterization of the common value of this classical Dynkin game via a (linear) doubly reflected BSDE, which generalizes the one shown in the literature  (cf. \cite{CK},\cite{HL},\cite{Lepeltier-Xu}) to the case of  
non right-continuous payoffs and general filtration.
More precisely, Theorem \ref{caracterisation} yields the following result.

\begin{corollary} (Existence and characterization of the  value of the classical Dynkin game) Suppose that the admissible pair of barriers $(\xi, \zeta)$ satisfies Mokobodzki's condition, that $\xi$ is right u.s.c.and $\zeta$ is right l.s.c.
 Then, for all $\theta \in \T_0$, 
$$ \overline{V}(\theta)=\underline{V}(\theta)= Y_\theta \quad {\rm a.s.},$$
where 
$Y$ is equal to the first component of the 
solution  of the DRBSDE (\ref{DRBSDE}) associated with driver $f=0$ and with  barriers $\xi$ and $\zeta$. 
\end{corollary}
It has also  been proven in the literature on classical Dynkin games (cf. \cite{ALM}) that if, moreover, $ \xi $ and $-\zeta$ are left u.s.c. along stopping times, and if $\xi_t < \zeta_t$, $t <T$ then
for each $\theta \in \T_0$, there exists a
saddle point at time $\theta$ for the  classical Dynkin game \eqref{dessusclassical} (cf. Theorem 4.3 in \cite{ALM}).
Note  that this result (even without the assumption $\xi < \zeta$) can be derived  by applying Corollary \ref{optimal} with $f=0$. 

\subsection{The general irregular case}

In this subsection $(\xi,\zeta)$ is an admissible pair of barriers satisfying Mokobodzki's condition. Contrary to the previous subsection, here we do not make   any regularity assumptions on the pair $(\xi,\zeta)$. 
Our   $\mathbfcal{E}^{^{f}}$-Dynkin game
might not admit a value (cf. the discussion in the introduction). 
In this general case, we will interpret   the  DRBSDE with  a pair of obstacles  $(\xi,\zeta)$  in terms of the value of "an extension" of the  zero-sum game of the previous subsection  over a larger set of "stopping strategies" than the set of  stopping times ${\cal T}_0$. 
To this purpose we  introduce the following notion of {\em stopping system}.

\begin{Definition}\label{weak}
Let $\tau \in {\cal T}_0$ be a  stopping time (in the usual sense). Let $H$ be a set in $\cf_\tau$. Let $H^c$ denote its complement in $\Omega$.  
The  pair  $\rho=(\tau,H)$ is called a \emph{stopping system} if $H^c\cap \{\tau=T\}=\varnothing.$
\end{Definition}  
By taking $H=\Omega$  in the above definition, we see that the notion of a stopping system generalizes that of a stopping time (in the usual sense). 
\begin{remark}
A {\em stopping system} is an example of {\em divided stopping time} (from the French "temps d'arr\^et divis\'e") in the sense of \cite{EK} Section 2.37
(see also Definition 3.1 in \cite{ALM} and the recent work by \cite{BB}).\\
Note that the notion of {\em stopping system} is simpler than the notion of {\em divided stopping time} (see 
Remark \ref{divise} for some additional comments on the usage of {\em divided stopping times} for classical linear Dynkin games).
\end{remark}


We denote by $\mathcal{S}_{0}$ the set of all stopping systems; for a stopping time $\theta\in{\cal T}_0$, we denote by $\mathcal{S}_{\theta}$ the set of stopping systems $\rho=(\tau,H)$ such that such that $\theta\leq \tau$.  \\

\noindent For an optional {\em right-limited} process $\phi$ and a stopping system $\rho=(\tau,H)$, we define $\phi_{\rho}$ by
$$\phi_{\rho}:=\phi_\tau \textbf{1}_H+{\phi}_{\tau^+} \textbf{1}_{H^c}.$$
In the particular case where $\rho=(\tau,\Omega)$, we have $\phi_\rho=\phi_\tau$, so the notation is consistent. \\

\noindent For  an optional  (not necessarily right-limited) process $\phi$ and for a  stopping system $\rho=(\tau,H)$, we set
$$\phi^{^u}_{\rho}:=\phi_\tau \textbf{1}_H+\hat{\phi}_\tau \textbf{1}_{H^c}\,\,{\rm and}\,\,
 \phi^{^l}_{\rho}:=\phi_\tau \textbf{1}_H+\check{\phi}_\tau \textbf{1}_{H^c},$$
 where $(\hat{\phi}_t)$ (resp. $(\check{\phi}_t)$) denotes the right upper- (resp. right lower-) semicontinuous envelope of the process $\phi$, defined by  $\hat{\phi}_t:=\limsup_{s \downarrow t, s> t} \phi_s$ (resp. $\check{\phi}_t:=\liminf_{s \downarrow t, s> t} \phi_s$), for all $t\in[0,T[$ (cf., e.g., \cite[page 133]{EK}). 
The process $\bar{\phi}$ (resp. $(\check{\phi}_t)$) is progressive and right upper- (resp. right lower-) semicontinuous. \\
Note that when $\phi$ is right-limited, we have $\phi^{^u}_{\rho}= \phi^{^l}_{\rho}= \phi_{\rho}$.\\ Moreover, in the particular case where $\rho=(\tau,\Omega)$, we have $\phi^{^u}_{\rho}= \phi^{^l}_{\rho}=\phi_\tau$, so the notation is consistent. \\
With the help of the above definitions and notation we formulate an extension of the  zero-sum game problem from Subsection \ref{subsect_right} where  the set of "stopping strategies" of the agents is the set of stopping systems. More precisely, for two stopping systems  $\rho=(\tau,H)\in\mathcal{S}_{0}$ and $\delta=(\sigma,G)\in\mathcal{S}_{0}$, we define the pay-off $I(\rho,\delta)$ by 
\begin{equation}\label{eq_new_payoff}
I(\rho,\delta):=\xi^{^u}_{\rho} {\bf 1} _{\tau \leq \sigma}+ \zeta^{^l}_{\delta} {\bf 1} _{\sigma < \tau}.
\end{equation}
We note that, by definition, $I(\rho,\delta)$ is an $\cf_{\tau\wedge \sigma}$-measurable random variable. 
 As in the previous subsection, the pay-off is assessed by an $f$-expectation, where $f$ is a Lipschitz driver. Let $\theta\in{\cal T}_0$ be a stopping time. 
The upper and lower value of the game at time $\theta$ are defined by:
\begin{equation}\label{dessus_2}
\bm{\overline{V}}(\theta):=\,\essinf_{\delta=(\sigma,G) \in \mathcal{S}_{\theta} }\,\esssup_{\rho=(\tau,H) \in \mathcal{S}_{\theta}} \, \mathbfcal{E}^{^{f}}_{_{\theta,\tau \wedge \sigma}}[I(\rho, \delta)]; \quad 
\bm{\underline{V}}(\theta):=\,\esssup_{\rho=(\tau,H) \in \mathcal{S}_{\theta}} \, \essinf_{\delta=(\sigma,G) \in \mathcal{S}_{\theta}} \, \mathbfcal{E}^{^{f}}_{_{\theta,\tau \wedge \sigma}}[I(\rho, \delta)].
\end{equation}

The other definitions from Definition \ref{Def_values} are generalized to the above framework in a similar manner, by replacing the set of stopping times $\T_\theta$ by the  set of stopping systems  $\mathcal{S}_{\theta}$. 
We will refer to this game problem as "extended" $\mathbfcal{E}^{f}$-Dynkin game (over the set of stopping systems). We will show that, for any $\theta\in\stopo$, the "extended" $\mathbfcal{E}^{f}$-Dynkin game  defined above has a value $\bm{V}(\theta)$, that  is, we have  $\bm{V}(\theta)=\bm{\overline{V}}(\theta)=\bm{\underline{V}}(\theta)$ a.s.,  and that this (common) value coincides with the first component of the solution (at time $\theta$) to the DRBSDE with driver $f$ and obstacles $(\xi,\zeta)$; we also show the existence of \emph{$\varepsilon$-optimal} stopping systems.  

Let $(Y,Z,k,h,A,A', C,C')$ be the solution of the DRBSDE (\ref{DRBSDE}). Let us give some definitions. For each $\theta$ $\in$ $ \T_0$ and each $\varepsilon >0$, we define the sets 
 $$A^\varepsilon:=\{(\omega,t)\in\Omega\times[0,T]: Y_t\leq \xi_t+\varepsilon\} \quad  B^\varepsilon:=\{(\omega,t)\in\Omega\times[0,T]: Y_t\geq \zeta_t-\varepsilon\}.$$   
 We recall that  the  stopping times $\tau^{\varepsilon}_{\theta}$ and $\sigma^{\varepsilon}_{\theta}$ have been defined  
as the \emph{débuts} after $\theta$ of the sets $A^\varepsilon$ and $B^\varepsilon$ (cf. Eq. \eqref{tauepsilon}). We now set 
$$H^\varepsilon:=\{\omega\in\Omega: (\omega,\tau^{\varepsilon}_{\theta}(\omega))\in A^\varepsilon\} \quad G^\varepsilon:=\{\omega\in\Omega: (\omega,\sigma^{\varepsilon}_{\theta}(\omega))\in B^\varepsilon\}$$ and we define the stopping systems 
\begin{equation}\label{rhoepsilon}
\rho^{\varepsilon}_{\theta} := (\tau^{\varepsilon}_{\theta}, H^\varepsilon) \text{ and }  \delta^{\varepsilon}_{\theta} := (\sigma^{\varepsilon}_{\theta}, G^\varepsilon).
\end{equation}
The following lemma uses an additional piece of notation.\\
For an optional \textit{right-limited} process $\phi$, and for two stopping systems  $\rho=(\tau,H)\in\mathcal{S}_{0}$ and $\delta=(\sigma,G)\in\mathcal{S}_{0}$, we set 
$$\phi_{\rho \Bwedge\delta}:=\phi_{\rho} {\bf 1} _{\tau \leq \sigma}+ \phi_{\delta} {\bf 1} _{\sigma < \tau}.
$$
\begin{remark}
For general stopping systems, the above notation is not symmetric (i.e. the equality   $\phi_{\rho\Bwedge\delta}= \phi_{\delta\Bwedge\rho}$ is not necessarily true). In the particular case where $\rho=(\tau,\Omega)$ and $\delta=(\sigma,\Omega)$ (i.e. the particular case of stopping times), we have 
$\phi_{\rho\Bwedge\delta}= \phi_{\tau\wedge\sigma}$, where    $\tau\wedge\sigma$ is the usual notation for the minimum of the two stopping times $\tau$ and $\sigma$, and we have the equality $\phi_{\rho\Bwedge\delta}= \phi_{\tau\wedge\sigma}=\phi_{\sigma\wedge\tau}=\phi_{\delta\Bwedge\rho}.$ 
\end{remark}

The following lemma is to be compared with Lemmas \ref{lalaun} and \ref{laladeux}. In the general irregular  case where $\xi$ is not necessarily right upper-semicontinuous and $\zeta$ is not necessarily left lower-semicontinuous, the inequalities \eqref{tau} from Lemma \ref{laladeux} do not necessarily hold true. In this case, working with the "regularized" processes $\xi^u$ and $\zeta^l$ and with stopping systems (instead of stopping times) allows us to have  inequalities which are analogous to those of Lemma \ref{laladeux}, as well as some properties which are, in a certain sense (cf. Remark \ref{remarqueST} below), analogous to those of Lemma \ref{lalaun}. 
\begin{lemma}\label{lala2} Let $(\xi, \zeta)$ be an admissible  pair of barriers satisfying Mokobodzki's condition.
Let $(Y,Z,k,h,A,A', C,C')$ be the solution of the DRBSDE (\ref{DRBSDE}). The following assertions hold:
\begin{enumerate}
\item We have
\begin{equation}\label{tau_2}
 Y_{\rho^{\varepsilon}_{\theta}} \, \leq \, \xi^{^u}_{\rho^{\varepsilon}_{\theta}} + \varepsilon  \quad {\rm and} \quad Y_{\delta^{\varepsilon}_{\theta}}\, \geq \,  \zeta^{^l}_{\delta^{\varepsilon}_{\theta}} - \varepsilon \quad \mbox{a.s.}
\end{equation}
\item  For all stopping systems $\rho=(\tau,H)$ and $\delta=(\sigma,G)$, we have 
\begin{equation}\label{lala11}
 \mathbfcal{E}^{^{f}}_{\theta,\tau_\theta^\varepsilon\wedge \sigma}[Y_{\rho_\theta^\varepsilon\Bwedge \delta}]\geq  Y _\theta \quad
\text{ and } \quad \mathbfcal{E}^{^{f}}_{\theta,\tau \wedge \sigma_\theta^\varepsilon}[Y_{\rho \Bwedge \delta_\theta^\varepsilon}]\leq  Y_\theta \quad \mbox{a.s.}
\end{equation}
\end{enumerate}
\end{lemma}

\begin{remark}\label{remarqueST}
Note that the inequalities \eqref{lala11} are the analogue, for the stopping systems $\rho_\theta^\varepsilon$, $\delta$, $\delta_\theta^\varepsilon$ and $\rho$,  of the inequalities \eqref{lala1} and \eqref{ba} satisfied by the stopping times $\tau_\theta^\varepsilon, $ $\sigma$, $\sigma_\theta^\varepsilon$ and $\tau$.
\end{remark}

\dproof 
Let us prove the first point. On the set $H^\varepsilon$, we have 
$Y_{\rho^{\varepsilon}_{\theta}}=Y_{\tau^{\varepsilon}_{\theta}}\leq \xi_{\tau^{\varepsilon}_{\theta}} + \varepsilon=\xi^{^u}_{\rho^{\varepsilon}_{\theta}} + \varepsilon$, where we have used  the definitions of $\rho^{\varepsilon}_{\theta}$,  $Y_{\rho^{\varepsilon}_{\theta}}$, $\xi^{^u}_{\rho^{\varepsilon}_{\theta}}$ and   $H^\varepsilon$.   
On the complement $H^{\varepsilon,c}$,  we have:  
\begin{equation}\label{eq_varepsilon}
Y_{\rho^{\varepsilon}_{\theta}}= Y_{\tau^{\varepsilon}_{\theta}+} \text{ and }\xi^{^u}_{\rho^{\varepsilon}_{\theta}}=\hat \xi_{\tau^{\varepsilon}_{\theta}}.
\end{equation}  On the other hand, by definitions of $\tau^{\varepsilon}_{\theta}$ and of $H^{\varepsilon,c}$, for a.e. $\omega\in\Omega$, there exists a decreasing sequence $(t_n):=(t_n(\omega))$ such that   $t_n(\omega)\downarrow\downarrow \tau^{\varepsilon}_{\theta}(\omega)$ and $Y_{t_n}\leq \xi_{t_n}+\varepsilon$, for all $n\in\N$. Hence,
$\limsupn Y_{t_n}(\omega)\leq \limsupn \xi_{t_n}(\omega)+\varepsilon.$ Now, be definition of $\hat \xi$,  we have $\limsupn \xi_{t_n}(\omega)\leq \hat\xi_{\tau^{\varepsilon}_\theta}(\omega)$. On the other hand, 
we have $\limsupn Y_{t_n}(\omega)=Y_{\tau^{\varepsilon}_\theta+}(\omega).$ Hence, 
$Y_{\tau^{\varepsilon}_\theta+}(\omega)\leq \hat \xi_{\tau^{\varepsilon}_\theta}(\omega)+\varepsilon.$  This inequality, together with \eqref{eq_varepsilon} gives that 
 $Y_{\rho^{\varepsilon}_{\theta}} \, \leq \, \xi^{^u}_{\rho^{\varepsilon}_{\theta}} + \varepsilon$ a.s. on $H^{\varepsilon,c}$. We thus derive the desired result, namely $Y_{\rho^{\varepsilon}_{\theta}}\leq \xi^{^u}_{\rho^{\varepsilon}_{\theta}}+\varepsilon$ a.s. on $\Omega$. 

Let us prove the second inequality.  
On the set $G^\varepsilon$, we have 
$Y_{\delta^{\varepsilon}_{\theta}}=Y_{\sigma^{\varepsilon}_{\theta}}\geq \zeta_{\sigma^{\varepsilon}_{\theta}} - \varepsilon$ $= \zeta^{^l}_{\delta^{\varepsilon}_{\theta}} - \varepsilon$, where we have used  the definitions of $\delta^{\varepsilon}_{\theta}$,  $Y_{\delta^{\varepsilon}_{\theta}}$, $\zeta^{^l}_{\delta^{\varepsilon}_{\theta}}$  and   $G^\varepsilon$.    
On the complement $G^{\varepsilon,c}$, we have 
 \begin{equation}\label{eq_zeta}
 Y_{\delta^{\varepsilon}_{\theta}}=Y_{\sigma^{\varepsilon}_\theta+} \,\,{\rm and}\,\, \zeta^{^l}_{\delta^{\varepsilon}_{\theta}}= \,\,\check{\zeta }_{\sigma^{\varepsilon}_\theta}.
 \end{equation}
 
 Now, for a.e. $\omega \in \Omega$,
  there exists  a decreasing sequence $(t_n):=(t_n(\omega))$ such that   $t_n(\omega)\downarrow\downarrow \sigma^{\varepsilon}_{\theta}(\omega)$ and $Y_{t_n}(\omega)\geq \zeta _{t_n}(\omega)-\varepsilon$, for all $n\in\N$. Hence,
  \begin{equation}\label{moinszeta}
\liminf_{n \rightarrow \infty} Y_{t_n}(\omega)\geq  \liminf_{n \rightarrow \infty} \zeta _{t_n}(\omega)-\varepsilon.
\end{equation}

Now,  $\liminf_{n \rightarrow \infty} Y_{t_n}(\omega)=Y_{\sigma^{\varepsilon}_\theta+}(\omega).$ Moreover, by definition of $\check{\zeta }$,  we have\\
 $ \liminf_{n \rightarrow \infty} \zeta _{t_n}(\omega)\geq \check{\zeta }_{\sigma^{\varepsilon}_\theta}(\omega)$. Hence, by \eqref{moinszeta}, we get $Y_{\sigma^{\varepsilon}_\theta+}(\omega)
\geq \check{\zeta }_{\sigma^{\varepsilon}_\theta}(\omega)-\varepsilon$. 
Using \eqref{eq_zeta}, we derive that on $G^{\varepsilon,c}$, $Y_{\delta^{\varepsilon}_{\theta}}\, \geq \, \zeta^{^l}_{\delta^{\varepsilon}_{\theta}} - \varepsilon $ a.s.\, We have thus shown that $Y_{\delta^{\varepsilon}_{\theta}}\, \geq \, \zeta^{^l}_{\delta^{\varepsilon}_{\theta}} - \varepsilon$ a.s.\,on $\Omega$.

 %
%
%
%
%
%
%
 
 Let us prove now the first inequality of \eqref{lala11}. We have $$Y_{\rho_\theta^\varepsilon\Bwedge \delta}=
 Y_{\rho_\theta^\varepsilon} {\bf 1} _{\tau_\theta^\varepsilon \leq \sigma}+ Y_{\delta} {\bf 1} _{\sigma < \tau_\theta^\varepsilon}.$$
 For \emph{the first term} of the second member of the equality, we have 
 $Y_{\rho_\theta^\varepsilon}=Y_{\tau_\theta^\varepsilon} {\bf 1} _{H^\varepsilon}+ Y_{\tau_\theta^\varepsilon+} {\bf 1} _{H^{\varepsilon,c}}.$ Now, on $H^{\varepsilon,c}$, we have $ Y_{\tau_\theta^\varepsilon}> \xi_{\tau_\theta^\varepsilon}+\varepsilon.$ The Skorokhod condition thus gives $\Delta C_{\tau_\theta^\varepsilon}=0$. This, together with Remark \ref{sautsY}, gives $(Y_{\tau_\theta^\varepsilon+}-Y_{\tau_\theta^\varepsilon})^-=0$. Hence, $Y_{\tau_\theta^\varepsilon+}\geq Y_{\tau_\theta^\varepsilon}$ on $H^{\varepsilon,c}$. Hence, $Y_{\rho_\theta^\varepsilon}\geq Y_{\tau_\theta^\varepsilon}$.   
 For \emph{the second term}, we have  $Y_{\delta} {\bf 1} _{\sigma < \tau_\theta^\varepsilon}= (Y_{\sigma} {\bf 1} _{H}+ Y_{\sigma+} {\bf 1} _{H^c}){\bf 1} _{\sigma < \tau_\theta^\varepsilon}$. By using the fact that $Y$ is a strong $\mathbfcal{E}^{f}$-submartingale on $[\theta, \tau_\theta^\varepsilon]$ (cf. Lemma \ref{lalaun}), we have $Y_{\sigma+}\geq Y_{\sigma}$ on   $\{\sigma < \tau_\theta^\varepsilon\}$. 
 Hence, $(Y_{\sigma} {\bf 1} _{H}+ Y_{\sigma+} {\bf 1} _{H^c}){\bf 1} _{\sigma < \tau_\theta^\varepsilon}\geq Y_\sigma {\bf 1} _{\sigma < \tau_\theta^\varepsilon}$. By combining the two terms, we get  $ Y_{\rho_\theta^\varepsilon\Bwedge \delta}\geq Y_{\tau_\theta^\varepsilon} {\bf 1} _{\tau_\theta^\varepsilon \leq \sigma}+ Y_{\sigma} {\bf 1} _{\sigma < \tau_\theta^\varepsilon}= Y_{\tau_\theta^\varepsilon\wedge \sigma}.$
 Using this and the nondecreasingness of $\mathbfcal{E}^{^{f}}_{\theta,\tau_\theta^\varepsilon\wedge \sigma}[\cdot]$, we obtain $\mathbfcal{E}^{^{f}}_{\theta,\tau_\theta^\varepsilon\wedge \sigma}[ Y_{\rho_\theta^\varepsilon\Bwedge \delta}]\geq \mathbfcal{E}^{^{f}}_{\theta,\tau_\theta^\varepsilon\wedge \sigma}[Y_{\tau_\theta^\varepsilon\wedge \sigma}].$  As $Y$ is  a strong $\mathbfcal{E}^{f}$-submartingale on $[\theta, \tau_\theta^\varepsilon]$ (cf. Lemma \ref{lalaun}), we get  
 $\mathbfcal{E}^{^{f}}_{\theta,\tau_\theta^\varepsilon\wedge \sigma}[ Y_{\tau_\theta^\varepsilon\wedge \sigma}]\geq Y_\theta$, from which we conclude that $ \mathbfcal{E}^{^{f}}_{\theta,\tau_\theta^\varepsilon\wedge \sigma}[Y_{\rho_\theta^\varepsilon\Bwedge \delta}]\geq  Y _\theta$. The proof of the second inequality of \eqref{lala11} is  similar.
\fproof

With the help of the previous lemma, we establish the following inequalities which are to be compared with the inequalities  \eqref{fifi} from Theorem \ref{caracterisation}.  
\begin{lemma}\label{Lemma_inequalities} 
The following inequalities hold:
\begin{equation}\label{fifibis}
\mathbfcal{E}^{^{f}}_{_{\theta, \tau \wedge \sigma^{\varepsilon}_{\theta}}}[I(\rho, \delta^{\varepsilon}_{\theta})]  - L\varepsilon  \,\,\leq \,\,  Y_{\theta} \,\,\leq \,\,  \mathbfcal{E}^{^{f}}_{_{\theta, \tau^{\varepsilon}_{\theta} \wedge \sigma}}[ I(\rho^{\varepsilon}_{\theta} , \delta) ]  + L\varepsilon \quad \mbox{a.s.}\,,
\end{equation}
where $L$ is a positive constant which only depends on  the Lipschitz constant $K$ of $f$ and on the terminal time $T$.
\end{lemma}
\dproof 
 Let $\theta$ $\in$ $\T_0$ and let $\varepsilon >0$. We first show the right-hand inequality. 
By Lemma \ref{lala2}, 
\begin{equation}\label{lala1bis}
Y_{\theta} \, \leq \,  \mathbfcal{E}^{^{f}}_{_{\theta, \tau^{\varepsilon}_{\theta} \wedge \sigma}}[Y_{\rho^{\varepsilon}_{\theta} \Bwedge \delta}] \quad \mbox{a.s.}\
\end{equation}
By definition of $Y_{\rho^{\varepsilon}_{\theta} \Bwedge \delta}$, we have
\begin{equation*}
Y_{\rho^{\varepsilon}_{\theta} \Bwedge \delta} \,=  \, Y_{\rho^{\varepsilon}_{\theta}}  {\bf 1} _{\tau^{\varepsilon}_{\theta}  \leq \sigma}+ Y_{\delta}  {\bf 1} _{\sigma < \tau^{\varepsilon}_{\theta} }\quad {\rm a.s.}
\end{equation*}

Now, $ Y_{\rho^{\varepsilon}_{\theta}} \, \leq \, \xi^{^u}_{\rho^{\varepsilon}_{\theta}} + \varepsilon$ (cf. Lemma \ref{lala2}). Moreover, since $ Y \leq \zeta$ and since $Y$ is right-limited,  we have 
$Y_{\delta}= Y^{^l}_{\delta} \leq  \zeta^{^l}_{\delta}$. We thus get
\begin{equation*}
Y_{\rho^{\varepsilon}_{\theta} \Bwedge \delta} \, \leq \, (\xi^{^u}_{\rho^{\varepsilon}_{\theta}} + \varepsilon) {\bf 1} _{\tau^{\varepsilon}_{\theta}  \leq \sigma}+ \zeta ^{^l}_\delta {\bf 1} _{\sigma < \tau^{\varepsilon}_{\theta} } \, \leq \,  I(\rho^{\varepsilon}_{\theta} , \delta) +\varepsilon \quad {\rm a.s.}
\end{equation*}
where the last inequality follows from the definition of $ I(\rho^{\varepsilon}_{\theta} , \delta) $.
By using the inequality $\eqref{lala1bis}$ and the nondecreasingness of $\mathbfcal{E}^{^{f}}$, we  derive 
\begin{equation}\label{22bis}
Y_{\theta} \, \leq \, \mathbfcal{E}^{^{f}}_{_{\theta, \tau^{\varepsilon}_{\theta} \wedge \sigma}}[ I(\rho^{\varepsilon}_{\theta} , \delta)  +\varepsilon] \leq \mathbfcal{E}^{^{f}}_{_{\theta, \tau^{\varepsilon}_{\theta} \wedge \sigma}}[ I(\rho^{\varepsilon}_{\theta} , \delta) ] + L \varepsilon\,\quad \mbox{a.s.}\,,
\end{equation}
where the last inequality follows from an estimate on BSDEs (cf. Proposition A.4 in \cite{16}).
Using similar arguments, it can be shown that 
$
 Y_{\theta}\,\, \geq \,\,\mathbfcal{E}^{^{f}}_{_{\theta, \tau \wedge \sigma^{\varepsilon}_{\theta}}}[I(\rho, \delta^{\varepsilon}_{\theta})]  - L\varepsilon$ a.s\,, which, together with \eqref{22bis}, leads to the desired inequalities \eqref{fifibis}.
 \fproof
 
 In the following theorem we show that the "extended" $\mathbfcal{E}^f$-Dynkin game has a value which coincides with the first component of the DRBSDE with irregular barriers. 
 \begin{theorem} [Existence of a value and  characterization]\label{Thm_extended_game}
 Let $f$ be a Lipschitz  driver satisfying Assumption $\eqref{Royer}$. Let $(\xi, \zeta)$ be an admissible  pair of barriers satisfying Mokobodzki's condition.
Let $(Y,Z,k,h,A,A', C,C')$ be the solution of the DRBSDE (\ref{DRBSDE}). 
There exists a common value  for the "extended" $\mathbfcal{E}^{f}$-Dynkin game, and for each stopping time $\theta$ $\in$ $\T_0$, we have
$$\bm{\underline{V}}(\theta)= Y_{\theta} = \bm{\overline{V}}(\theta) \text{ a.s.}$$

Moreover, for each $\theta$ $\in$ $\T_0$ and each $\varepsilon >0$, the pair of stopping systems $(\rho^{\varepsilon}_{\theta}, \delta^{\varepsilon}_{\theta})$, defined by \eqref{rhoepsilon}, is an $L \varepsilon$-saddle point at time $\theta$ for the "extended" $\mathbfcal{E}^f$-Dynkin game, that is satisfies 
the  inequalities \eqref{fifibis}.
\end{theorem}

 \dproof
 The proof relies on  Lemma \ref{Lemma_inequalities}. Since the right-hand  inequality in \eqref{fifibis} from Lemma \ref{Lemma_inequalities} holds for all $\delta=(\sigma,G) \in \mathcal{S}_{\theta}$, we have 
 $$\,Y_{\theta} \leq \, \essinf_{\delta=(\sigma,G) \in \mathcal{S}_{\theta}} \mathbfcal{E}^{^{f}}_{_{\theta, \tau^{\varepsilon}_{\theta} \wedge \sigma}}[I(\rho^{\varepsilon}_{\theta} , \delta)]+L\varepsilon \leq \, \,\esssup_{\rho=(\tau,H) \in {\cal S}_\theta} \, \essinf_{\delta=(\sigma,G) \in \mathcal{S}_{\theta}} \, \mathbfcal{E}^{^{f}}_{_{\theta,\tau \wedge \sigma}}[I(\rho, \delta)] +L\varepsilon 
 \,\,\text{ a.s.}\,$$
 From this, together with  the definition of  $\bm{\underline{V}}(\theta)$ (cf. 
\eqref{dessus}), we obtain  $Y_{\theta} \,\,\leq \,\, \bm{\underline{V}}(\theta)+ L \varepsilon$ a.s. Similarly, we  show that $\bm{\overline{V}}(\theta)- L \varepsilon \,\, \leq \,\, Y_{\theta}$ a.s. for all $\varepsilon>0$. We thus get
$\bm{\overline{V}}(\theta) \,\, \leq \,\, Y_{\theta} \,\,\leq \,\, \bm{\underline{V}}(\theta)$ a.s.\,  This, together with the inequality $\bm{\underline{V}}(\theta) \leq \bm{\overline{V}}(\theta)$ a.s.\,,  yields 
$\bm{\underline{V}}(\theta)= Y_{\theta} = \bm{\overline{V}}(\theta) \text{ a.s.}$ The proof is thus complete.
%
\fproof


\begin{remark}\label{divise}
Recall that in the case of the classical Dynkin game (that is when $f=0$), in \cite{ALM}, the authors state that in the case of complete irregular payoffs (satisfying Mokobodzki's condition), the Dynkin game might does not admit a value (cf. Section 3 in \cite{ALM}).\footnote{see also Example \ref{contre-exemple} in the present paper.}
This leads them to  extend the game so that there exists a value for the extended game. More precisely, they consider the extended game to the set of {\em divided stopping times} for which they show that there exists a value  (cf. Theorem 3.3 in \cite{ALM}).\\
One might think to generalize this property to the case of (non linear) expectation ${\cal E}^f$, but it appears that the {\em divided stopping times} are not appropriate to our non linear case.\footnote{They do not allow us to obtain the inequalities of type \eqref{lala11} (which are the analogous of Lemma \ref{lalaun} in the case of right-u.s.c. barriers).}
\end{remark}

\section{Two useful corollaries}\label{sec6}
Using the characterization of the solution of the  nonlinear DRBSDE as the value function of the  "extended" $\mathbfcal{E}^{f}$-Dynkin game (over the set of stopping systems)
from Theorem \ref{Thm_extended_game}, we now establish a  comparison theorem and  \emph{a priori} estimates with universal constants (i.e.  depending only on the terminal time $T$ and the 
common Lipschitz constant $K$) for DRBSDEs with completely irregular barriers.


\begin{corollary}[Comparison theorem for DRBSDEs.]\label{thmcomprbsde}
  Let $(\xi^{1},  \zeta^{1})$ and  $(\xi^{2}, \zeta^{2})$ be two admissible pairs of barriers satisfying Mokobodzki's condition.
Let $f^1, f^2$ be Lipschitz drivers satisfying Assumption \ref{Royer}. For $i=1,2$, let $(Y^{i}, Z^{i}, k^i, A^i, A^{' i}, C^i, C^{' i})$  be the solution of the DRBSDE associated with driver $f^i$ and barriers $\xi^i$, $\zeta^i$. \\
Assume  that $\xi^{2}\le \xi^{1}$ and $\zeta^{2}\le \zeta ^{1}$ and 
$f^{2}(t,Y_t^2, Z_t^2, k_t^2) \le f^{1}(t,Y_t^2, Z_t^2, k_t^2)$ $dP\otimes dt$-a.s.\,\\
Then, we have
$Y^{2}\le
Y^{1}$.
\end{corollary}


\dproof 
{\bf Step 1}: Let us first assume that $\xi^2 \le \xi^1$, $\zeta^{2}\le \zeta ^{1}$, and that 
$f^{2}(t,y,z,\mathpzc{k}) \le f^{1}(t,y,z,\mathpzc{k})$ for all  $(y,z,\mathpzc{k}) \in \R^2 \times L^2_\nu,$ $dP\otimes dt$-a.s.\, 
Let $\theta\in\mathcal{S}_0$. For $i=1,2$ and for all stopping systems $\rho=(\tau,H)\in {\cal S}_\theta, \delta=(\sigma,G) \in {\cal S}_\theta$, let $\mathbfcal{E}^{i}_{_{\cdot, \tau \wedge \sigma}}[I^i(\rho,\delta)]$ be the first coordinate of the solution 
 of
the BSDE associated with driver $f^{i}$, terminal time $\tau \wedge \sigma$ and terminal condition $I^i(\rho, \delta)=(\xi^i)^u_{\rho} \textbf{1}_{ \tau \leq \sigma}+(\zeta^i)^l_{\delta}\textbf{1}_{\sigma<\tau}$. Since $\xi^{2}\le \xi^{1}$ and $\zeta^{2}\le \zeta ^{1}$, we have $I^2(\rho, \delta) \leq I^1(\rho, \delta)$ a.s. Since, moreover $f^2 \leq f^1$, the comparison theorem 
for BSDEs gives: for all stopping systems $\rho=(\tau,H)\in{\cal S}_\theta, \delta=(\sigma,G) \in {\cal S}_\theta$, 
$\mathbfcal{E}^{2}_{_{\theta, \tau \wedge \sigma}}[I^2(\rho,\delta)]\, \leq \, \mathbfcal{E}^{1}_{_{\theta, \tau \wedge \sigma}}[I^1(\rho,\delta)]$ a.s.\,
Taking the essential supremum over $\rho$ in ${\cal S}_\theta$ and the essential infimum over $\delta$ in ${\cal S}_\theta$ in this inequality, and using  the characterization of the solution of the  DRBSDE with  obstacles  $(\xi,\zeta)$ as the value function of the  "extended" $\mathbfcal{E}^{f}$-Dynkin game (cf. Theorem \ref{Thm_extended_game}), we obtain: 
$$ Y^{2}_{\theta}= \, \essinf_{\delta=(\sigma,G) \in {\cal S}_\theta} \,\esssup_{\rho=(\tau,H) \in {\cal S}_\theta}\, \mathbfcal{E}^{2}_{_{\theta, \tau \wedge \sigma}}[I^2(\rho,\delta)]\, \leq \,  \essinf_{\delta=(\sigma,G) \in {\cal S}_\theta} \,\esssup_{\rho=(\tau,H) \in {\cal S}_\theta}\, \mathbfcal{E}^{1}_{_{\theta, \tau \wedge \sigma}}[I^1(\rho,\delta)]= Y^{1}_{\theta} \quad {\rm a.s.}\,$$
Since this inequality holds 
for each $\theta \in  {\cal T}_0$, we get $Y^{2}\le
Y^{1}$.\\
 {\bf Step 2}: We now place ourselves under  the assumptions of the theorem (which are weaker than those made in Step 1). Let $\tilde f$ be the process defined by $\tilde f_t:= f^{2}(t,Y_t^2, Z_t^2, k_t^2)-f^{1}(t,Y_t^2, Z_t^2, k_t^2) $, which, by assumption, is non positive. Note that $(Y^2, Z^2, k^2)$ is the solution of the DRBSDE associated with barriers $\xi ^{2},\zeta ^{2}$ and driver $f^1(t,y,z,\mathpzc{k}) +  \tilde f_t$. We have $f^1(t,y,z,\mathpzc{k})+ \tilde f_t $ $\leq$ $f^1(t,y,z,\mathpzc{k}) $  for all $(y,z,\mathpzc{k})$. By Step 1 applied to the driver $f^1$ and the driver 
$f^1(t,y,z,\mathpzc{k}) + \tilde f_t$ (instead of $f^2$), we get  $Y^2 \leq Y^1$. 
\fproof

Using  the results provided in Theorem \ref{Thm_extended_game}, in particular
the existence of $\varepsilon$-saddle points (of stopping systems) for the "extended" $\mathbfcal{E}^f$-Dynkin game,
 we  now prove the following estimates for the spread of the solutions of two DRBSDEs with completely irregular barriers.

\begin{corollary}[A priori estimates for DRBSDEs] \label{oubli} Let $(\xi^{1},  \zeta^{1})$ and  $(\xi^{2}, \zeta^{2})$ be two admissible pairs of barriers satisfying Mokobodzki's condition.
Let $f^1, f^2$ be Lipschitz drivers satisfying Assumption \ref{Royer} with common Lipschitz constant $K>0$. For $i=1,2$, let $Y^i$ be the solution of the DRBSDE associated with driver $f^i$ and barriers $\xi^i$, $\zeta^i$.  \\
Let $\tilde{Y}:=Y^1-Y^2$, $\tilde{\xi}:=\xi^1-\xi^2$, $\tilde{\zeta}:=\zeta^1-\zeta^2$.  Let $\eta, \beta >0$ with $\beta\geq \dfrac{3}{\eta}+2C$ and $\eta \leq \dfrac{1}{C^2}$. Setting $\delta f_s:= f^{2}(t,Y_s^2, Z_s^2, k_s^2)-f^{1}(t,Y_s^2, Z_s^2, k_s^2)$, $0 \leq s \leq T$, for each $\theta \in \T_0$, we have 
\begin{equation}\label{eqA.1}
(\tilde{Y_\theta}) ^2  \leq e^{\beta (T-\theta)}E[\esssup_{\tau \in {\cal T}_\theta} \tilde{\xi_\tau}^2+\esssup_{\tau \in {\cal T}_\theta}  \tilde{\zeta_\tau}^2| \mathcal{F}_\theta]+ 
 \eta E[\int_\theta^Te^{\beta (s-\theta)}(\delta f_s)^2 ds|\mathcal{F}_\theta] \quad {\rm a.s.}\,
\end{equation}
\end{corollary}

%

\dproof 
The proof is divided into two steps.\\
 {\bf Step 1}:
For $i=1,2$ and for all {\em stopping systems} $\rho=(\tau,H), \delta=(\sigma,G) \in {\cal S}_\theta$, let $(X^{i,\rho, \delta}$, $\pi^{i,\rho, \delta}, l^{i, \rho, \delta})$  be the solution of the BSDE associated with driver $f^i$, terminal time $\tau \wedge \sigma$ and terminal condition $I^i(\rho, \delta)$, where $I^i(\rho, \delta)=(\xi^i)^u_{\rho} \textbf{1}_{ \tau \leq \sigma}+(\zeta^i)^l_{\delta}\textbf{1}_{\sigma<\tau}$.
Set  $\tilde{X}^{\rho, \delta}:=X^{1,\rho, \delta}-X^{2,\rho, \delta}$ and $\tilde{ I}(\rho, \delta):= I^1(\rho, \delta)- I^2(\rho, \delta)=((\xi^1)^u_{\rho}-  (\xi^2)^u_{\rho})\textbf{1}_{\tau \leq \sigma}+((\zeta^1)^l _{\delta}-(\zeta^2)^l _{\delta})\textbf{1}_{\sigma<\tau}$.\\
By an estimate on BSDEs (see Proposition $A.4$ in \cite{17}), for each $\theta \in \T_0$, we have a.s.:
\begin{equation*}\label{A.2}
 (\tilde{X}_{\theta}^{\tau, \delta})^2  \leq e^{\beta (T-\theta)} E[\tilde{ I}(\rho, \delta)^2 \mid \mathcal{F}_{\theta}]+ \eta  E[\int_{\theta}^T e^{\beta (s-\theta)}{[(f^1-f^2)(s, X_s^{2,\rho, \delta},\pi_s^{2,\rho, \delta}, l_s^{2,\rho, \delta}) ]^2ds} \mid \mathcal{F}_{\theta}].
\end{equation*}
From this, together with the definitions of $(\xi^i)^u_{\rho}$ and $(\zeta^i)^l _{\delta}$,  we derive 
\begin{equation}\label{A.3}
 (\tilde{X}_\theta^{\rho, \delta})^2  \leq e^{\beta (T-\theta)}E[\esssup_{\tau \in {\cal T}_\theta} \tilde{\xi_\tau}^2+\esssup_{\tau \in {\cal T}_\theta}  \tilde{\zeta_\tau}^2| \mathcal{F}_\theta]+ 
 \eta E[\int_\theta^Te^{\beta (s-\theta)}(\tilde f_s)^2 ds|\mathcal{F}_\theta] \quad {\rm a.s.}\,
\end{equation}
where $\tilde{f}_s:= \sup_{y,z,\mathpzc{k}}|f^1(s,y,z,\mathpzc{k})-f^2(s,y,z,\mathpzc{k})|$.\\
For each $\varepsilon >0$, let $\rho_{\theta}^{1,\varepsilon}$ (resp. $\delta_{\theta}^{2,\varepsilon}$) 
be the {\em stopping system} $\rho_{\theta}^{\varepsilon}$ (resp. $\delta_{\theta}^{\varepsilon}$) associated with $(Y^1, \xi^1)$ (resp. $(Y^2, \zeta^2)$) defined by \eqref{rhoepsilon}. 
By using inequality \eqref{fifibis} in Lemma \ref{Lemma_inequalities}, we obtain that for all $\varepsilon >0$ and for all {\em stopping systems} $\rho$, $\delta \in {\cal S}_\theta$,
\begin{equation*}
Y_\theta^1-Y_\theta^2 \leq X_\theta^{1, \rho_{\theta}^{1,\varepsilon}, \delta}-X_\theta^{2, \rho, \delta_{\theta}^{2,\varepsilon}}+2L\varepsilon \quad {\rm a.s.}
\end{equation*}
Applying this inequality to the {\em stopping systems} $\rho=\rho_{\theta}^{1,\varepsilon}$ and $\delta= \delta_{\theta}^{2,\varepsilon}$, we get
\begin{equation*}\label{ineq2}
Y_\theta^1-Y_\theta^2 \leq X_\theta^{1,\rho_{\theta}^{1,\varepsilon}, \delta_{\theta}^{2,\varepsilon}}-X_\theta^{2, \rho_{\theta}^{1,\varepsilon}, \delta_{\theta}^{2,\varepsilon}}+2L\varepsilon
\leq |X_\theta^{1, \rho_{\theta}^{1,\varepsilon}, \delta_{\theta}^{2,\varepsilon}}-X_\theta^{2, \rho_{\theta}^{1,\varepsilon}, \delta_{\theta}^{2,\varepsilon}}|+2L\varepsilon \quad {\rm a.s.}
\end{equation*}
This inequality together with $\eqref{A.3}$ gives
$$
Y_\theta^1-Y_\theta^2 \leq \sqrt{ e^{\beta (T-\theta)}E[\esssup_{\tau \in {\cal T}_\theta} \tilde{\xi_\tau}^2+\esssup_{\tau \in {\cal T}_\theta}  \tilde{\zeta_\tau}^2| \mathcal{F}_\theta]+ 
 \eta E[\int_\theta^Te^{\beta (s-\theta)}(\tilde f_s)^2 ds|\mathcal{F}_\theta]}+2L\varepsilon \quad {\rm a.s.}
$$
By symmetry, the last inequality is also verified by $Y_\theta^2-Y_\theta^1$. Since this holds for all $\varepsilon>0$, we derive that
\begin{equation*}\label{eqA.1b}
 (\tilde{Y_\theta})^2  \leq e^{\beta (T-\theta)}E[\esssup_{\tau \in {\cal T}_\theta} \tilde{\xi_\tau}^2+\esssup_{\tau \in {\cal T}_\theta}  \tilde{\zeta_\tau}^2| \mathcal{F}_\theta]+ 
 \eta E[\int_\theta^Te^{\beta (s-\theta)}(\tilde f_s)^2 ds|\mathcal{F}_\theta]\quad {\rm a.s.}
\end{equation*}
This result holds for all Lipschitz drivers $f^1$ and $f^2$  satisfying Assumption \ref{Royer}.\\
{\bf Step 2}: 
Note that $(Y^2, Z^2, k^2)$ is the solution the DRBSDE associated with barriers $\xi ^{2},\zeta ^{2}$ and 
driver $f^1(t,y,z,\mathpzc{k}) + \delta f_t$.  By applying the result of Step 1 to the driver $f^1(t,y,z,\mathpzc{k})$ and the driver $f^1(t,y,z,\mathpzc{k}) + \delta f_t$ (instead of $f^2$), we get the desired result.
\fproof

\begin{remark}
The previous two corollaries show  the relevance of the characterization  of the solution of the (non-linear) DRBSDE with irregular obstacles as the value of an "extended" $\mathcal{E}^f$-Dynkin game, as established in Theorem \ref{Thm_extended_game}. In particular, this characterization, 
together with  the existence of $\varepsilon$-saddle points (of stopping systems) for the "extended" $\mathcal{E}^f$-Dynkin game, allows us to provide estimates with universal constants which, it seems, cannot be obtained by using Gal'chouk-Lenglart's formula. 
Indeed, up to now in the literature,  It\^o-type techniques have not proved useful for showing estimates with universal constants,  even  in the  simplest case of continuous barriers and Brownian filtration (cf. Remark 4.5 in \cite{DQS2} for details).
\end{remark}

\section{Application to the pricing of game options with irregular payoffs}\label{gameoptions}
In this paragraph, we illustrate how the results of 
Section \ref{subsect_right} can be applied to the problem of pricing of game options (with non-right-continuous pay-offs) in a class of  market models with imperfections. \\
We set $E:= \R$, $\nu(de) := \lambda \delta_1(de)$, where $\lambda$ is a positive constant, and where $\delta_1$ denotes the Dirac measure at $1$. 
The process $N_t := N([0,t] \times \{1\})$ is then a Poisson process with parameter $\lambda$, and we have
$\tN_t := \tN([0,t] \times \{1\})  = N_t - \lambda t.$ 

We  assume that the filtration is the natural filtration associated with $W$ and $N$.

 We consider a financial market which consists of one risk-free asset, whose price process $S^0$ satisfies $dS_t^{0}=S_t^{0} r_tdt$ with $S_0^0=1$, and two risky assets with price processes $S^{1},S^{2}$ satisfying:
 $$dS_t^{1}=S_{t^-}^{1}[\mu_t^1dt +  \sigma^1_t dW_t+ \beta_t^1d\tilde N_t];\quad dS_t^{2}=S_{t^-}^{2} [\mu^2_tdt+\sigma^2_tdW_t + \beta_t^2d\tilde N_t],$$
 with $S_0^1>0$ and $S_0^2>0$. 
We suppose that the processes $\sigma^1,\sigma^2,$ $\beta^1,\beta^2,$ $r, \mu^1,\mu^2$ are 
predictable and bounded, with $\beta_t^i>-1$ for $i=1,2$. 
Let $\mu_t:= (\mu^1,\mu^2)'$ and let $\Sigma_t:= (\sigma_t, \beta_t)$ be the $2\times 2$-matrix with first column 
$\sigma_t:=(\sigma_t^1,\sigma_t^2)'$ and second column  $\beta_t:=(\beta_t^1,\beta_t^2)'$. 
We suppose that $\Sigma_t$ is invertible and that the coefficients of $\Sigma_t^{-1}$ are bounded.

We  consider an agent who can invest his/her initial wealth $x\in\R$ in the three assets. 

For $i=1,2$, we denote by $\varphi_t^i$ the amount invested in the $i^{\textit{th}}$ risky asset. 
A process $\varphi= (\varphi^1, \varphi^2)'$ belonging to ${\mathbb H}^2 \times  {\mathbb H}^2_{\nu}$ will be called a \emph{portfolio strategy}.

The value of the associated portfolio (or {\em wealth}) at time $t$ is denoted  by $X^{x, \varphi}_t$ (or simply by $X_t$). 
In the case of a perfect market, we have
\begin{align*}\label{portfolio}
dX_t & = (r_t (X_t- \varphi_t^1 - \varphi_t^2)+\varphi_t^1 \mu^1_t +\varphi_t^2 \mu^2_t  ) dt +
(\varphi_t^1 \sigma^1_t + \varphi_t^2 \sigma^2_t) dW_t +
(\varphi_t^1 \beta^1_t + \varphi_t^2 \beta^2_t) d\tilde N_t \\
& = (r_t X_t+\varphi_t' (\mu_t - r_t{\bf 1}) ) dt+ \varphi_t' \sigma_t dW_t + 
\varphi_t' \beta_t  d\tilde N_t,
\end{align*}
where ${\bf 1}=(1,1)'$.
%
%
%
More generally, we will suppose that there may be some \emph{imperfections} in the market,  taken into account via 
the {\em nonlinearity} of the
dynamics of the wealth and encoded in a Lipschitz driver $f$ satisfying Assumption~\ref{Royer} (cf. the example hereafter; cf. also \cite{EQ96} or \cite{DQS3} for other examples of imperfections). More precisely, 
we suppose that  the  wealth process  $X^{x, \varphi}_t$ (also $X_t$)
satisfies  the forward differential equation: 
\begin{equation}\label{riche}
-dX_t= f(t,X_t, {\varphi_t}' \sigma_t, {\varphi_t}' \beta_t) dt - {\varphi_t}' \sigma_t dW_t-{\varphi_t}' \beta_t d\tilde N_t, \;  ; \;  X_0=x,
\end{equation}
 or,
equivalently, setting $Z_t= {\varphi_t}' { \sigma}_t$ and
  $k_t= {\varphi_t}' \beta_t $,
 \begin{equation}\label{wea}
-dX_t= f(t,X_t, Z_t,k_t ) dt -  Z_t dW_t- k_t d\tilde N_t ; \;  X_0=x.
\end{equation}
Note that $(Z_t, k_t)= {\varphi_t}' { \Sigma}_t$, which  is equivalent to ${\varphi_t}' = (Z_t, k_t)\, { \Sigma}_t^{-1}$.

\begin{remark}\label{richessemartingale}
Note that the wealth process $X^{x, \varphi}$ is an $\mathcal{E}^f$-martingale, since
$X^{x, \varphi}$ is the solution of the 
BSDE with driver $f$, terminal time $T$ and terminal condition $X_T^{x, \varphi}$. 
\end{remark}

\begin{remark}
This model includes the case of a perfect market, for which $f$ is a linear driver given by
$
f(t,y,z,k)= -r_t y- (z, k)\, { \Sigma}_t^{-1} (\mu_t - r_t{\bf 1}).
$

Another example is given by the case of a  borrowing rate $R_t$ different from the lending rate $r_t$, satisfying generally $R_t \geq r_t$ (cf. e.g. \cite{Korn} and \cite{C1}). In this case, we have:
\begin{align*}
dX_t & = (r_t (X_t- \varphi_t^1 - \varphi_t^2)^+- R_t (X_t- \varphi_t^1 - \varphi_t^2)^-
+\varphi_t^1 \mu^1_t +\varphi_t^2 \mu^2_t  ) dt + \varphi_t' \sigma_t dW_t + 
\varphi_t' \beta_t  d\tilde N_t \\
& = (r_t X_t+\varphi_t' (\mu_t - r_t{\bf 1}) - (R_t -r_t)(X_t- \varphi_t^1 - \varphi_t^2)^-) dt+ \varphi_t' \sigma_t dW_t + 
\varphi_t' \beta_t  d\tilde N_t.
\end{align*}
The wealth process $X$ thus satisfies the dynamics \eqref{riche} with the driver $f$ given by
$f(t,x, {\varphi}' \sigma_t, {\varphi}' \beta_t)= -r_t x-\varphi' (\mu_t - r_t{\bf 1}) + (R_t -r_t)(x- \varphi' {\bf 1})^-$.

Our model also includes the case of a repo market on which 
 the  risky assets are supposed to be traded (cf. \cite{Brigo}). For $i=1,2$, we denote by $b^i_t$ (resp. $l^i_t$) the borrowing (resp. lending) repo rate
 \footnote{also called repo rate for long (resp. short) position}  at time $t$ for the $i^{\textit{th}}$ risky asset. The processes $b^i_t$, $l^i_t$, for $i=1,2$, are supposed to be bounded predictable processes. In this case, we have
 \footnote{$\varphi^i_t >0$ means that we need some risky asset number $i$, so we borrow it, while if 
  $\varphi^i_t <0$, we lend it.}
$$dX_t=(r_t X_t+\varphi_t' (\mu_t - r_t{\bf 1}) 
+ l^1_t (\varphi^1_t)^-  - b^1_t (\varphi^1_t)^+  + l^2_t (\varphi^2_t)^-  - b^2_t (\varphi^2_t)^+) dt + \varphi_t' \sigma_t dW_t + 
\varphi_t' \beta_t  d\tilde N_t.
$$
Hence,  the wealth process $X$ satisfies  \eqref{riche} with the driver $f$ given by\\
$f(t,x, {\varphi}' \sigma_t, {\varphi}' \beta_t)= -r_t x-\varphi' (\mu_t - r_t{\bf 1}) - l^1_t (\varphi^1)^-  + b^1_t (\varphi^1)^+ -  l^2_t (\varphi^2)^-  + b^2_t (\varphi^2)^+$.

An other example is given by the case of a large seller  (see e.g. \cite{GDQS} Section 3.5 or \cite{DQS3} for details).
\end{remark}

Let $T>0$ be a given terminal time. Let $(\xi, \zeta)$ be  an admissible pair of processes.\\
 We recall that a \textit{game option} is a financial instrument which gives the buyer the right  \textit{to exercise} at any  stopping  time 
$\tau \in  {\cal T}$ and the seller the right  \textit{to cancel} at any stopping time  $\sigma \in  {\cal T}$. If the buyer exercises at time $\tau$ before the seller cancels, then the seller pays to the buyer the amount $ \xi_{\tau}$; if the seller cancels at time $\sigma$ before the buyer exercises,  the seller pays to the buyer   
 the amount $ \zeta_{\sigma}$ at the cancellation time $\sigma$. The difference 
$\zeta - \xi \geq 0$ corresponds to a penalty which the seller pays to the buyer in the case of an early cancellation of the contract. 
Thus, if the seller chooses a cancellation time $\sigma$ and the buyer chooses an exercise time $\tau$,  the former pays to the latter 
the payoff $I(\tau, \sigma)$ (defined in \eqref{eq0}) 
at time 
$\tau \wedge \sigma.$ 

\begin{definition}
For an initial wealth $x\in\R$, a {\em super-hedge} (for the seller) against the game option is a pair $(\sigma, \varphi)$ of a stopping time $\sigma \in {\cal T}$ and a portfolio strategy $ \varphi$ $\in$  ${\mathbb H}^2 \times  {\mathbb H}^2_{\nu}$ such that \footnote{Note that condition \eqref{condA} is equivalent to 
$X^{x, \varphi}_{t \wedge \sigma} \geq I(t,\sigma), \;\; 0 \leq t \leq T \quad {\rm a.s.}$}
\begin{equation}\label{condA}
X^{x, \varphi}_{t } \geq \xi_t,  \; 0\leq t \leq \sigma \; \text{ a.s. and } X^{x, \varphi}_{\sigma } \geq \zeta_{\sigma}
 \text{ a.s.} 
\end{equation}
We denote by ${\cal S} (x)$ the  set of all super-hedges associated with the initial wealth $x$. The  {\em seller's price} $u_0$ of the game option (at time $0$) is defined by 
\begin{equation} \label{seller's price}
 u_0:= \inf \{x \in \R,\,\, \exists  (\sigma, \varphi) \in {\cal S} (x) \}.
 \end{equation}
\end{definition}

 With the help of Theorem \ref{caracterisation}, Proposition \ref{lalabis} (first assertion) and Proposition \ref{sautsYA} (iii), we give  a dual formulation for the seller's   price $u_0$ via the 
 $\mathbfcal{E}^f$-Dynkin game, and we show the existence of a super-hedge for the seller under a left regularity assumption on the "cancellation payoff" $\zeta$.  We also obtain a characterization of $u_0$ in terms of the DRBSDE \eqref{def_solution_DRBSDE}.

\begin{proposition}[Seller's  price and super-hedge of the game option] \label{superhedging} 
Suppose that $(\xi, \zeta)$ satisfies  Mokobodzki's condition, that $\xi$ is right-u.s.c. and $\zeta$ is right-l.s.c.
Assume moreover that  $\zeta$ is left lower-semicontinuous along stopping times. The seller's  price  $u_0$ satisfies:
\begin{equation}\label{u0}
u_0=\inf_{\sigma \in \mathcal{T} }  \sup_{\tau \in \mathcal{T}} \mathbfcal{E}^f_{0,\tau \wedge \sigma}[I(\tau, \sigma)]= 
\sup_{\tau\in \mathcal{T} } \inf_{\sigma \in \mathcal{T}} 
\mathbfcal{E}^f_{0,\tau \wedge \sigma}[I(\tau, \sigma)]=Y_0,
\end{equation}
where $(Y, Z,k, h, A,A',C,C')$ is the  solution of the 
DRBSDE  
associated with driver $f$ and barriers $\xi$ and $ \zeta$.

Let  $\sigma^*:= \inf \{ t \geq 0,\,\, Y_t = \zeta_t\}$ and $\overline\sigma:= \inf \{ t \geq 0,\,\, A'_t > 0\,\,{\rm or}\,\, C'_{t^-}>0\}$. The pairs $(\sigma^*, \varphi^*)$ and  $(\overline\sigma, \varphi^*)$ are super-hedges  associated with the initial amount $u_0$.
\end{proposition}

\begin{remark}
This result generalizes Theorem 3.12 in \cite{DQS3} which concerns the case when the payoff processes $\xi$ and $\zeta$ are right-continuous. 
 \end{remark}

\begin{remark}
In the  special case of a perfect market model, our result gives that $u_0$ is equal to the 
  value  of a classical Dynkin game problem, which generalizes the results previously
 shown in the literature (cf., e.g., \cite{Kifer,H}) to the case of non-right-continuous payoff processes. 
  \end{remark}

\begin{proof}  
The proof 
relies on Theorem \ref{caracterisation}, Proposition \ref{lalabis}, and on similar arguments to those in  \cite{DQS3} (where the case of game options with RCLL payoffs in a market with default is treated). \\
By Theorem \ref{caracterisation}, we have $Y_0=\inf_{\sigma \in \mathcal{T} }  \sup_{\tau \in \mathcal{T}} \mathbfcal{E}^f_{0,\tau \wedge \sigma}[I(\tau, \sigma)]= 
\sup_{\tau\in \mathcal{T} } \inf_{\sigma \in \mathcal{T}} 
\mathbfcal{E}^f_{0,\tau \wedge \sigma}[I(\tau, \sigma)].$ 
Hence,  in order to prove the other assertions of the   theorem,  it  remains to show that $u_0=Y_0$, $( \sigma^*,\varphi^*) \in {\cal S}(Y_0)$, and $(\overline\sigma, \varphi^*)\in {\cal S}(Y_0)$.\\
%
%
Let us first show that
$( \sigma^*,\varphi^*)$ and $(\overline\sigma, \varphi^*)$ belong to ${\cal S}(Y_0)$.
 Recall that $ \sigma^* \leq \overline\sigma$ a.s. (cf. Remark \ref{etoilebar}).
 Since $\zeta$ is left-l.s.c. along stopping times, 
 the process  
$A'$ is continuous (cf.  Proposition \ref{continuous} (iii)). 
By definition of $\overline\sigma$, we derive that the processes $(A'_t)$ and $(C'_{t-})$ are  a.s. constant on $[0, \overline\sigma)$, and hence on $[0, \overline\sigma]$ since they are left-continuous.
 For almost every $\omega$, for all $t \in [0, \overline\sigma(\omega)]$, we thus have
\begin{align}\label{forwardb}
Y_t(\omega)=Y_0-\int_0^t  f(s,\omega, Y_s(\omega),Z_s(\omega),k_s(\omega))ds+M_t(\omega)-A_t(\omega)-C_{t-}(\omega), 
\end{align}
where $M_t:= \int_0^t Z_sdW_s+
\int_0^t k_s d\tilde N_s$.
Now, for almost every $\omega$,  the wealth $X_.^{Y_0, \varphi^* }(\omega)$ (associated with the initial capital $Y_0$ and the portfolio strategy $\varphi^*$) satisfies  the deterministic forward  differential equation:
\begin{align}\label{ri}
X_t^{Y_0, \varphi^*}(\omega) = Y_0-\int_0^t f(s,X_s^{Y_0, \varphi^*}(\omega),Z_s(\omega),k_s(\omega))ds + M_t(\omega)
, \,\, 
0 \leq t \leq T.
 \end{align}
%
Since $(A_t)$ and $(C_{t-})$ are non-decreasing, by applying the classical comparison result on $[0, \overline\sigma(\omega)]$  for the two forward differential equations \eqref{forwardb} and \eqref{ri}, 
we derive that, for almost every $\omega$, we have
 \begin{equation}\label{ooo}
X_t^{Y_0, \varphi^*}(\omega) \geq Y_t (\omega), \quad 0 \leq t \leq \overline\sigma (\omega).
\end{equation}
From this  inequality and from the fact that $Y$ is greater than or equal to the lower barrier $\xi$,   we deduce $X_t^{Y_0, \varphi^*} \geq Y_t\geq \xi_t,$  $0 \leq t \leq \overline\sigma \quad {\rm a.s.}$ Since $A'$ is continuous, by Proposition \ref{lalabis} (first assertion), we have 
 $Y_{\sigma^*} =\zeta_{\sigma^*}$ and $Y_{\overline\sigma} =\zeta_{\overline\sigma}$ a.s.
Hence, from the inequality \eqref{ooo}, we  get    
$  X_{\sigma^*}^{Y_0, \varphi^*} \geq \zeta_{\sigma^*} \quad {\rm a.s.}$ and
$  X_{\overline\sigma}^{Y_0, \varphi^*} \geq \zeta_{\overline\sigma} \quad {\rm a.s.}$
 We conclude that $(\sigma^*,\varphi^*)$ and $(\overline\sigma, \varphi^*)$ belong to ${\cal S}(Y_0)$. 
 From this and the definition of $u_0$, we immediately derive the inequality $Y_0 \geq u_0$. 
   
To complete the proof of the theorem, it remains to show the converse inequality $u_0\geq Y_0$.
Let $x \in \R$ be such that there exists $( \tilde\sigma,\tilde\varphi) \in {\cal S}(x)$. We show that $x \geq Y_0$. 
 Since $( \tilde\sigma,\tilde\varphi) \in {\cal S}(x)$, we have
$X^{x, \tilde\varphi}_{t } \geq \xi_t$, $0\leq t \leq \tilde\sigma$ a.s. and 
  $X^{x, \tilde\varphi}_{\tilde\sigma } \geq \zeta_{\tilde \sigma}$ a.s. 
For each $\tau \in \mathcal{T}$, we thus get the inequality
$X^{x,\tilde \varphi}_{\tau \wedge \tilde\sigma} \geq I(\tau, \tilde\sigma) \quad {\rm a.s.}$.
By the non-decreasing property 
of $\mathcal{E}^f$, together with the $\mathcal{E}^f$-martingale property of $X^{x, \tilde\varphi}$ (cf. Remark \ref{richessemartingale}),
we thus obtain 
 $x= \mathcal{E}^f _{0,\tau \wedge \tilde\sigma}[X^{x, \tilde\varphi}_{\tau \wedge \tilde\sigma}] \geq 
\mathcal{E}_{0, \tau \wedge \tilde \sigma}^f [ I(\tau, \tilde\sigma)]$, for each $\tau \in \mathcal{T}$. By taking 
the supremum over $\tau \in \mathcal{T}$, we derive 
\begin{align*}
x  \geq 
 \sup_{\tau \in \mathcal{T}}
\mathcal{E}_{0, \tau \wedge \tilde \sigma}^f [ I(\tau, \tilde \sigma)]\geq \inf_{\sigma \in \mathcal{T}} \sup_{\tau \in \mathcal{T}}
\mathcal{E}_{0, \tau \wedge \sigma}^f [ I(\tau, \sigma)],
\end{align*}
where the last inequality is obvious. Hence, 
$x  \geq 
 \inf_{\sigma \in \mathcal{T}} \sup_{\tau \in \mathcal{T}}
\mathcal{E}_{0, \tau \wedge \sigma}^f [ I(\tau, \sigma)]=Y_0$ (we recall that the equality is due to Theorem \ref{caracterisation}). As $x$ is an arbitrary initial capital for which there exists a super-hedge and as $u_0$ is defined as the infimum of such $x$'s, we get $u_0\geq Y_0$, which is the desired inequality. 
The proof is thus complete. 
\end{proof}

\begin{remark}
By similar arguments as those used in the proof of Theorem 3.11 in \cite{DQS3}, we can show that the equalities \eqref{u0} still hold even if $\zeta$ is not left lower-semicontinuous along stopping times. However, in that case, there does not necessarily exist a super-hedge for the seller of the game option.
\end{remark}

We now give some examples of game options, also called {\em callable American options}, with non right-continuous pay-offs. 
\begin{Example} 
We consider a callable American {\em basket} option 
 whose payoff processes $\xi$ and $\zeta$  are of the form: 
 $\xi_t:=g(S^1_t)h(S^2_t)$ and $\zeta_t:=\delta g(S^1_t)$ for $t\in[0,T]$, where $\delta$ is a positive constant, where $h:\R\rightarrow [0,\delta]$ is an u.s.c. function, and where $g:\R\rightarrow \R^+$ is a continuous function such that $(g(S_t^1))\in\mathcal{S}^2$. Note that the process 
 $(g(S_t^1))$ is right-continuous and left-continuous along stopping times (since the process $(S_t^1)$ is right-continuous and left-continuous along stopping times).
  Note that since an u.s.c. function can be written as the limit of a (non increasing) sequence of continuous functions,
 the process  $(h(S_t^2))$ is optional, which implies that $\xi$ is optional. Moreover, since $h$ is u.s.c.\,and since the process $(S_t^2)$ is right-continuous, the process $(h(S_t^2))$, and hence $(\xi_t)$, is right-u.s.c. 
Note that we have $\xi \leq \zeta$. Suppose that
the pair $(\xi, \zeta)$ 
 satisfies  Mokobodzki's condition. By Proposition \ref{superhedging}, the seller's price of the game option is equal to the value function of the ${\cal E}^f$-Dynkin game, and is also
characterized as the solution of the DRBSDE associated with driver $f$ and
 barriers $\xi$ and $\zeta$. \\
 An example of an game option of this type is given by $g(x):=(x-K)^+$ (or $(K-x)^+$) and $h(x):=
 {\bf 1}_{[K, + \infty)}(x)$. Note that $h$ is u.s.c. \footnote{Note that the process $\xi$  is right-u.s.c.\,, but, in general neither left-limited nor right-limited. }
Moreover, by Tanaka's formula, the process $\zeta_t = \delta g(S_t^1)$ is a semimartingale. Hence, by Lemma \ref{Moki_lem} and Remark \ref{Moki_rem} in the Appendix, Mokobodzki's condition holds.

 \end{Example}

 \begin{Example}
 We now consider a cancellable American option of {\em Lookback} type, whose payoff processes are of the form: 
 $\xi_t:=g(S^1_t)h(\inf_{0\leq s\leq t} S^1_s)$ and $\zeta_t:=\delta g(S^1_t)$ for $t\in[0,T]$, where the functions $g$ and $h$ satisfy the same assumptions as in the previous example.
By the same arguments as above, we get the characterizations of the seller's price of the game option  given in Proposition \ref{superhedging}.\\
For example, {\em cancellable  American call options with lower barrier} enter in this framework. In this case, 
 the payoff processes are of the form:
 $\xi_t:= ( S^1_t - K)^+{\bf 1}_{ \inf_{0\leq s\leq t} S^1_s \geq L}$ and 
 $\zeta_t:= \delta  ( S^1_t - K)^+$, where $\delta$ is a constant greater than or equal to $1$, and $K$ and $L$ are positive constants such that $K >L$. Note that  the "cancellation" payoff process $\zeta$ can also be chosen equal to $( S^1_t - K)^+ + \delta$, where $\delta$ is a positive constant.
\end{Example}

\paragraph{Conclusion}
1. In this paper, we have thus formulated a notion of doubly reflected BSDE in the case of completely irregular barriers and general filtration. \\
2. We have first shown that the existence of a  solution is equivalent to 
the so-called Mokobodzki's condition.\\
3. We have then shown that, if a  solution exists, it is unique.\\
4. In the case where the barriers $\xi$ and $- \zeta$ are right-u.s.c. (and satisfy  Mokobodzki's condition), 
we have shown that the unique solution of the doubly reflected BSDE is characterized as the value of an 
$\mathbfcal{E}^{f}$-{\em Dynkin game} {\bf over stopping times} (for a general filtration). Moreover, when 
$ \xi $ and $-\zeta$ are left u.s.c. along stopping times, there exists  
saddle point for this game.\\
$\bullet$ We have given an application of this characterization to the pricing of game options (whose payoffs $\xi$ and $- \zeta$ are  right-u.s.c.) in an imperfect market.\\
5. In the case where the barriers $\xi$ and $\zeta$ are completely irregular (and satisfying Mokobodzki's condition), the unique solution of the doubly reflected BSDE is characterized as the value of an "extended"
$\mathbfcal{E}^{f}$-{\em Dynkin game} {\bf over stopping systems}. 
Moreover,  this "extended"
$\mathbfcal{E}^{f}$-{\em Dynkin game} admits $\varepsilon$-saddle point.
\\
6. Using the  results of the previous point, 
we have shown a comparison theorem, as well as a priori estimates with universal constants for doubly reflected BSDEs with completely irregular barriers.

It is still an open question to establish, in the Markovian case, some links between DRBSDEs with non-right-continuous barriers and related obstacle problems. The results provided in the present paper, in particular the comparison theorem and the a priori estimates with universal constants, will be valuable tools in the study of this problem. 

It is also interesting to establish  analogous results for other types of DRBSDEs, such as, for example, 
DRBSDEs driven by a random measure.

\section{Appendix}

\subsection{Reflected BSDE with driver $0$ and irregular obstacle }\label{sec3}

Let $T>0$ be a fixed terminal time. 
Let $\xi= (\xi_t)_{t\in[0,T]}$ be a 
 process in ${\cal S}^2$,
 called  \emph{barrier}, or  \emph{obstacle}.

 The following result has been proved in  \cite{MG3}  (cf. Theorem 3.1):
\begin{proposition} \label{rexiuni}
Let  $\xi$  be a process in $\mathcal{S}^2$.
There exists a unique solution of the reflected BSDE with driver equal to $0$ and obstacle $\xi$, that is a unique process $(X,\pi, l, h,A,C)\in {\cal S}^2 \times \H^2 \times \H^2_\nu \times {\cal M}^{2, \bot} \times{\cal S}^2\times {\cal S}^2  $ such that
\begin{align}\label{RBSDE}
&X_t=\xi_T-\int_{t}^T  \pi_s dW_s-\int_{t}^T \int_E l_s(e) \tilde N(ds,de) - (h_T-h_t)+A_T-A_t+C_{T-} -C_{t-},
0\leq t \leq T \text{ a.s.,}\\
&  X_t \geq \xi_t   \,\, \text{ for all } t\in[0,T]  \text{ a.s.,}\nonumber\\
& A \text{ is a nondecreasing right-continuous predictable process 
with } A_0= 0 \text{ and such that } \nonumber\\
& \int_0^T {\bf 1}_{\{X_{t^-} > \overline{\xi}_t\}} dA_t = 0 \text{ a.s.  } 
\nonumber\\  
& C \text{ is a nondecreasing right-continuous adapted purely discontinuous process with } C_{0-}= 0  \nonumber\\
& \text{ and such that }
(X_{\tau}-\xi_{\tau})(C_{\tau}-C_{\tau-})=0 \text{ a.s. for all }\tau\in{\cal T}_0. \nonumber
\end{align}
\end{proposition}

We introduce the following operator:
\begin{definition}[Operator induced by an RBSDE with driver $0$]\label{definitionRef}
For a  process $(\xi_t)$ $\in$  $\mathcal{S}^2$,  we denote by $\Ref[\xi]$ the first component of the solution to the Reflected BSDE with (lower) barrier $\xi$ and with driver $0$. 
\end{definition}
\begin{remark} Note that by Proposition \ref{rexiuni} the  operator $\Ref: \xi \mapsto 
\Ref[\xi]$ 
 is well-defined on $\mathcal{S}^{2}$.
 \end{remark}

 We give some useful properties of the operator $\Ref$ in the following two lemmas. 

\begin{lemma}\label{stopping_result}
The operator $\Ref: S^{2}\rightarrow S^{2}$ satisfies the following properties:
\begin{enumerate}
\item
The operator $\Ref$ is nondecreasing,  
that is, for  $\xi, \xi'$ $\in S^{2}$ such that $\xi \leq  \xi'$ 
we have $\Ref[\xi]\leq \Ref[\xi']$.  
\item
If $\xi$ $\in S^2$ is a strong  supermartingale, then $\Ref[\xi]=\xi.$ 
\item
For each $\xi$ $\in S^{2}$,  $\Ref[\xi]$ is a strong supermartingale and satisfies $\Ref[\xi]\geq \xi$.
\end{enumerate}
\end{lemma}
\dproof 
By definition, we have $\Ref[\xi]=X$, where $X= ( X_t)_{t\in[0,T]}$ is the first coordinate of the solution of the reflected BSDE \eqref{RBSDE}.
Now, by   Theorem 3.1 
 in \cite{MG3},
the process $X$ is equal to the value function of the classical optimal stopping problem with payoff $\xi$, that is for each stopping time
$\theta$, we have
$$
  X_\theta  = \esssup_{\tau \in \stops} E[ \xi_{\tau}  \mid \Fc_\theta]\quad {\rm a.s.}
$$
Hence, by classical results of Optimal Stopping Theory, the process $ \Ref[\xi]=X$ is equal to
the {\em Snell envelope} of the process $\xi$, that is, the smallest strong supermartingale
 greater than or equal to $\xi$. Using this property, we easily derive the three assertions of the lemma.
 \fproof
 
 \begin{remark} \label{Rk_monotone} We recall that the nondecreasing limit of a sequence of strong supermartingales is a strong supermartingale
(which can be easily shown by the Lebesgue theorem for  conditional expectations).
 \end{remark}

 We now show a monotone convergence result for  the operator $\Ref$.
\begin{lemma} \label{monotoneconvergence}
Let $(\xi^n)$ be a sequence of processes belonging to $S^{2}$, supposed to be  nondecreasing, i.e., such that 
for each $n\in \N$, $\xi^n \leq \xi^{n+1}$. Let $\xi:= \lim_{n \rightarrow + \infty} \xi^n$. 
If $\xi \in {\cal S}^2$, then  $\Ref[ \xi] = \lim_{n \rightarrow + \infty} \Ref[ \xi^n]$.
\end{lemma}

\dproof
As the operator $\Ref$ is nondecreasing, the sequence  $(\Ref[ \xi ^n])$ is nondecreasing.
Let $ X:= \lim_{n \rightarrow + \infty} \Ref[ \xi ^n]$. 
Again, due to the nondecreasingness of   the operator $\Ref$,  we have $\Ref[\xi ^0] \leq \Ref[\xi ^n] \leq \Ref[\xi]$, for all $n \in \N$. 
By letting $n$ go to $+\infty$, we get $\Ref[\xi ^0] \leq X$ and 
\begin{equation}\label{leq}
X \leq \Ref[\xi].
\end{equation} 
In particular, we have $X \in {\cal S}^2$. Let us now show that $X \geq \Ref[\xi].$
By definition of $
\Ref[\xi ^n]$ as the solution of the reflected BSDE with obstacle $\xi ^n$, we have $\Ref[\xi ^n] \geq \xi ^n$, for all $n \in \N$. By letting $n$ go to $+\infty$, we get $X \geq \xi$. Hence, 
\begin{equation} \label{inter}
\Ref[X] \geq \Ref[\xi].
\end{equation}
We note now that for each $n \in \N$, $\Ref[\xi ^n]$ is a strong supermartingale (cf. Lemma \ref{stopping_result}). It follows that
$X$ is a strong supermartingale as the nondecreasing limit of a sequence of strong supermartingales (cf. Remark \ref{Rk_monotone}).
Hence,  $X= \Ref[X]$ (cf. Lemma \ref{stopping_result}, second 
assertion). By \eqref{inter}, we thus have $X \geq \Ref[\xi]$, which, using \eqref{leq}, implies
$X = \Ref[\xi]$.
\fproof

\subsection{Proofs}

\paragraph{Proof of Proposition \ref{seq}}  The proof relies on Lemmas \ref{stopping_result} and 
\ref{monotoneconvergence}.
We first show that  ${\cal X}^n \geq 0$ and $ {\cal X}^{'n}\geq 0$, for all $n\in\N$. By definition, ${\cal X}^{'n}_T= {\cal X}^{n}_T =0$. 
Since $\Tilde{\xi}_T^{f}= \Tilde{\zeta}_T^{f}=0$, it follows that 
$({\cal X}^{' n}+\Tilde{\xi}^{f}){\bf 1}_{[0,T)}= {\cal X}^{' n}+\Tilde{\xi}^{f}$ and 
$({\cal X}^{ n}-\Tilde{\zeta}^{f}){\bf 1}_{[0,T)}= {\cal X}^{ n}-\Tilde{\zeta}^{f}$. Moreover, since ${\cal X}^{n}$ is 
a strong supermartingale, we have  $ {\cal X}^{n}_\theta \geq 
E[ {\cal X}^{n}_T\, \vert \Fc_\theta]$ $=0$ a.s.  for all $\theta\in \T_0$,  which implies that ${\cal X}^n \geq 0$. \footnote{Recall that, by a result of the general theory of processes, if $\phi \in {\cal S}^2$ and $\phi' \in {\cal S}^2$  are such that 
  $\phi_\theta \leq \phi'_\theta$ a.s. for all $\theta\in\T_0,$ then $\phi \leq \phi'$.}
Similarly, we see that  $ {\cal X}^{'n}\geq 0$.

We prove recursively that $({\cal X}^{n})_{n \in {\mathbb N}}$ and $({\cal X}^{'n})_{n \in {\mathbb N}}$ are nondecreasing sequences of processes. We 
have ${\cal X}^{1} \geq 0={\cal X}^{0}$ and ${\cal X}^{'1} \geq 0={\cal X}^{'0}$. Suppose that 
${\cal X}^{n} \geq {\cal X}^{n-1}$ and ${\cal X}^{'n} \geq {\cal X}^{'n-1}$. The induction hypothesis and the nondecreasingness of the operator 
$\Ref$ (cf. Lemma \ref{stopping_result}) give 
\begin{align}\label{classique}
\Ref[{\cal X}'^{n}+\Tilde{\xi}^{f}]\geq \Ref[{\cal X}^{'n-1}+\Tilde{\xi}^{f}] \quad ; \quad 
\Ref[{\cal X}^{n}-\Tilde{\zeta}^{f}] \geq \Ref[{\cal X}^{n-1}-\Tilde{\zeta}^{f}].
\end{align}
Hence,  ${\cal X}^{n+1} \geq {\cal X}^{n}$ and ${\cal X}^{'n+1} \geq {\cal X}^{'n}$, which is the desired result.

We now  define two processes $H^{f}$ and $H^{'f}$ as follows: 
$$
H_t^{f}:=H_t+E[\xi_T^-|\cal{F}_t]+E[\int_t^Tf^-(s)ds|\cal{F}_t]; \quad
H_t^{'f}:=H'_t +E[\xi_T^+|\cal{F}_t]+E[\int_t^Tf^+(s)ds|\cal{F}_t],
$$
where $H$ and $H'$ come from Mokobodzki's condition for $(\xi,\zeta)$ (cf. Eq. \eqref{Moki}). 
We note that $H^{f}$ and $H^{'f}$  are nonnegative strong supermartingales in  ${\cal S}^2$. 
From Mokobodzki's condition, we get
\begin{equation}\label{Mokodeux}
  \tilde \xi^{f} \leq H^{f} -H^{'f} \leq \tilde \zeta^{f}.
\end{equation}

We prove recursively that  ${\cal X}^n \leq H^{f}$ and ${\cal X}^{'n}\leq H^{'f}$, for all $n \in {\mathbb N}$.  
Note first that ${\cal X}^0 =0\leq H^{f}$ and ${\cal X}^{'0}=0\leq H^{'f}$. Suppose now that ${\cal X}^n \leq H^{f}$ and 
${\cal X}^{'n}\leq H^{'f}$. 
From this, together with \eqref{Mokodeux}, we get ${\cal X}^{'n}\leq H^{'f} \leq H^{f}- \Tilde{\xi}^{f}$, which implies 
${\cal X}^{'n} + \Tilde{\xi}^{f} \leq H^{f}$. 
Since 
the operator 
$\Ref$ is non decreasing, we derive  ${\cal X}^{n+1}= \Ref [{\cal X}^{'n} + \Tilde{\xi}^{f} ] \leq \Ref[H^{f}]$. Since $H^{f}$ is a strong supermartingale, the second  assertion of Lemma \ref{stopping_result} gives
$\Ref[H^{f}]=H^{f}$. 
Hence, ${\cal X}^{n+1} \leq H^{f}$. Similarly, we show ${\cal X}^{'n+1}\leq H^{'f}$. The desired conclusion follows.

By definition, we have ${\cal X}^{f} =\lim \uparrow 
{\cal X}^{n}$ and ${\cal X}^{'f}=\lim \uparrow 
{\cal X}^{'n}$.
The processes ${\cal X}^{f}$ and ${\cal X}^{'f}$ are optional (valued in $[0, + \infty]$) as the limit of sequences of  optional (nonnegative)
processes.
 Since for all $n \in {\mathbb N}$, ${\cal X}^{n}_T= 
{\cal X}^{'n}_T =0$  a.s.\,, we have ${\cal X}^{f}_T={\cal X}^{'f}_T =0$ a.s.\,
Moreover, since  for all $n \in {\mathbb N}$, $0 \leq {\cal X}^n \leq H^{f}$ and $0 \leq {\cal X}^{'n}\leq H^{'f}$, we obtain
$0 \leq {\cal X}^{f} \leq H^{f}$ and $0 \leq {\cal X}^{'f} \leq H^{'f}$. As $H^{f}, H^{'f}$ $\in$ ${\cal S}^2$, 
it follows that ${\cal X}^{f}$ and ${\cal X}'^{f}$ belong to ${\cal S}^2$.

Moreover, ${\cal X}^{f}$ and ${\cal X}^{'f}$ are strong supermartingales as limits of nondecreasing sequences of strong supermartingales (cf. Remark \ref{Rk_monotone}).

%
%

It remains to show that ${\cal X}^{f}$ and ${\cal X}^{'f}$ are solutions of the system \eqref{system}. 
Recall that, since ${\cal X}^{'n}_T=\Tilde{\xi}_T^{f}=0$, by \eqref{sys1}, we have  for all  $n \in \N$,
\begin{equation}\label{sys1bis}
{\cal X}^{ n+1}=\Ref[{\cal X}^{' n}+\Tilde{\xi}^{f}]. 
\end{equation}
Note  that the sequence $({\cal X}^{'n}+\Tilde{\xi}^{f})_{ n \in \N } $  is non decreasing and converges to 
${\cal X}^{'f}+\Tilde{\xi}^{f}$. 
By Lemma \ref{monotoneconvergence}, we thus derive that 
 $\lim_{n \rightarrow + \infty} \Ref[{\cal X}^{'n}+\Tilde{\xi}^{f}] = \Ref[{\cal X}^{' f}+\Tilde{\xi}^{f}]$.
 Hence, by letting $n$ tend to $+ \infty$ in \eqref{sys1bis}, we get
 ${\cal X}^{f} = \Ref[{\cal X}^{'f}+\Tilde\xi^{f}]$. Similarly, it can be shown that ${\cal X}^{'f}=\Ref[{\cal X}^{f}-\Tilde \zeta^{f}]$. Since 
  ${\cal X}^{f}_T={\cal X}^{'f}_T =0$ a.s.\,, it follows that ${\cal X}^{f}$ and ${\cal X}^{'f}$ are solutions of the system \eqref{system}. 
  
   Note now that ${\cal X}^{f}$, ${\cal X}^{'f}$ satisfy the inequalities $ \tilde \xi^{f} \leq {\cal X}^{f} -{\cal X}^{'f} \leq \tilde \zeta^{f}.$ Moreover, they are the minimal nonnegative strong supermartingales  in $ \mathcal{S}^2$ satisfying
these inequalities. Indeed, if $J,J^{'}$ are nonnegative strong supermartingales  in $ \mathcal{S}^2$ satisfying $ \tilde \xi^{f} \leq J -J^{'} \leq \tilde \zeta^{f}, $ then, using the same arguments as above, we derive that ${\cal X}^{f}\leq J$ and ${\cal X}^{'f}\leq J^{'}$.

From this minimality property, it follows that  $({\cal X}^{f}, {\cal X}^{'f})$ is also characterized as the minimal solution of the system \eqref{system} of coupled RBSDEs. 
\fproof

\paragraph{Proof of Lemma \ref{Lemma_estimate}}
Let $\beta>0$ and  $\varepsilon>0$ be such that $\beta\geq\frac 1 {\varepsilon^2}$.
We  set $\tilde Y:=Y-\bar Y$, $\tilde Z:=Z-\bar Z$, $\tilde A:=A-\bar A$, $\tilde A':=A'-\bar A'$, $\tilde C:=C-\bar C$, $\tilde C':=C'-\bar C'$, $\tilde k:=k-\bar k$, 
$\tilde h:=h-\bar h$,
 and $\tilde f(\omega, t):=f(\omega, t)-\bar f(\omega, t)$.
We note that $\tilde Y_T=\xi_T-\xi_T=0;$ moreover, 
$$-d\tilde Y_t=\tilde f_t dt +d\tilde A_t - d\tilde A'_t +d\tilde C_{t-}- d\tilde C'_{t-}-\tilde Z_t dW_t- \int_{E} \tilde k_t(e) \tilde N(dt,de)- d\tilde h_t,\quad t\in[0,T].$$
Thus we see that $\tilde Y$ is an {\em optional strong  semimartingale} in the vocabulary of \cite{Galchouk} 
%
 with decomposition 
$\tilde Y_t= \tilde Y_0+ M_t+\alpha_t+\gamma_t$, where
$M_t:=\int_0^t \tilde Z_s dW_s+\int_0^t \int_{E} \tilde k_s(e) \tilde N(ds,de) + \tilde h_t$, $\alpha_t:= -\int_0^t  \tilde f_s ds- \tilde A_t+\tilde A'_t$ 
and 
$\gamma_t:=-\tilde C_{t-}+\tilde C'_{t-}$ (cf., e.g.,  Theorem A.3.  and Corollary A.2 in \cite{MG}).
Applying Gal'chouk-Lenglart's formula (cf.  Corollary A.2 in \cite{MG}) to $\e^{\beta t}\tilde Y_t^2$, and using the property  
$\langle \tilde h^c, W \rangle =0$, 
we get: almost surely, for all $t\in[0,T]$,
\begin{equation*}
\begin{aligned}
\e^{\beta T}\tilde Y_T^2&=\e^{\beta t}\tilde Y_t^2+\int_{]t,T]}\beta\e^{\beta s} (\tilde Y_{s})^2 ds
-2\int_{]t,T]} \e^{\beta s}\tilde Y_{s-}\tilde f_s ds -2\int_{]t,T]} \e^{\beta s}\tilde Y_{s-}d\tilde A_s \\
&+2\int_{]t,T]} \e^{\beta s}\tilde Y_{s-}d\tilde A'_s-
\int_{[t,T[} 2\e^{\beta s}\tilde Y_{s}d(\tilde C)_{s+}+
\int_{[t,T[} 2\e^{\beta s}\tilde Y_{s}d(\tilde C')_{s+}\\
&+ \sum_{t<s\leq T}\e^{\beta s}(\tilde Y_s-\tilde Y_{s-})^2+\sum_{t\leq s<T}\e^{\beta s}(\tilde Y_{s+}-\tilde Y_{s})^2 +\\
&+\int_{]t,T]} \e^{\beta s} \tilde Z_s^2 ds+
 \int_{]t,T]} \e^{\beta s} d\langle \tilde h^c \rangle _s   + \tilde M_T-\tilde M_t,
\end{aligned}
\end{equation*}
where 
\begin{equation}\label{defM}
\tilde M_t:= 2\int_{]0,t]} \e^{\beta s}\tilde Y_{s-}\tilde Z_s d W_s+2
\int_{]0,t]} \e^{\beta s}\int_E \tilde Y_{s-}\tilde k_s(e) \tilde N(ds,de)
  +2\int_{]0,t]} \e^{\beta s}\tilde Y_{s-}d\tilde h_s.
 \end{equation}

Thus, we get (recall that $\tilde Y_T=0$): almost surely, for all $t\in[0,T]$,
\begin{equation}\label{eq0_lemma_estimate}
\begin{aligned}
\e^{\beta t}\tilde Y_t^2+ \int_{]t,T]} \e^{\beta s} \tilde Z_s^2 ds +
 \int_{]t,T]} \e^{\beta s} d\langle \tilde h^c \rangle _s 
 &=  -\int_{]t,T]}\beta\e^{\beta s} (\tilde Y_{s})^2 ds+
2\int_{]t,T]} \e^{\beta s}\tilde Y_{s}\tilde f_s ds+2\int_{]t,T]} \e^{\beta s}\tilde Y_{s-}d\tilde A_s\\
&-2\int_{]t,T]} \e^{\beta s}\tilde Y_{s-}d\tilde A'_s +2\int_{[t,T[} \e^{\beta s}\tilde Y_{s}d\tilde C_{s}-2\int_{[t,T[} \e^{\beta s}\tilde Y_{s}d\tilde C'_{s}\\
& -\sum_{t<s\leq T}\e^{\beta s}(\tilde Y_s-\tilde Y_{s-})^2 
-\sum_{t\leq s<T}\e^{\beta s}(\tilde Y_{s+}-\tilde Y_{s})^2 -( \tilde M_T-\tilde M_t).
\end{aligned}
\end{equation}
We give hereafter an upper bound for some of the terms appearing on the right-hand side (r.h.s. for short) of the above equality.

Let us first consider the sum of \textbf{the first and the  second term } on the r.h.s. of equality  \eqref{eq0_lemma_estimate}. By applying the inequality $2ab\leq (\frac a \varepsilon)^2+\varepsilon^2 b^2$, valid for all $(a,b)\in\R^2$, we get: a.s. for all $t\in[0,T]$,
\begin{equation*}
\begin{aligned}
-\int_{]t,T]}\beta\e^{\beta s} (\tilde Y_{s})^2 ds+2\int_{]t,T]}  \e^{\beta s}\tilde Y_{s}\tilde f_s ds
&\leq(\frac 1 {\varepsilon^2}-\beta)\int_{]t,T]}\e^{\beta s} (\tilde Y_{s})^2 ds+\varepsilon^2 \int_{]t,T]}  \e^{\beta s}\tilde f^2(s) ds.
\end{aligned}
\end{equation*}
As $\beta\geq\frac 1 {\varepsilon^2}$, we have $(\frac 1 {\varepsilon^2}-\beta)\int_{]t,T]}\e^{\beta s} (\tilde Y_{s})^2 ds \leq 0$, for all $t\in[0,T]$ a.s. 

For \textbf{the third  term} (resp. \textbf{the fourth  term}) on the r.h.s. of \eqref{eq0_lemma_estimate} it can be shown that, a.s. for all $t\in[0,T]$,  $+2\int_{]t,T]} \e^{\beta s}\tilde Y_{s-}d\tilde A_s\leq 0$ (resp. $-2\int_{]t,T]} \e^{\beta s}\tilde Y_{s-}d\tilde A'_s\leq 0$) The proof uses property \eqref{RBSDE_A} of the definition of the DRBSDE and the properties  $Y\geq \xi$, $\bar Y\geq \xi$ (resp. $Y\leq \zeta $, $\bar Y\leq \zeta$) ; the details are similar to those in the case of  RBSDE (with one lower obstacle) (cf., for instance,  the proof of Lemma 3.2 in \cite{MG}).

For the \textbf{fifth} and \textbf{sixth} terms on the r.h.s. of \eqref{eq0_lemma_estimate} we show that,  a.s. for all $t\in[0,T]$, $+2\int_{[t,T[} \e^{\beta s}\tilde Y_{s}d\tilde C_{s}\leq 0$  and  $-2\int_{[t,T[} \e^{\beta s}\tilde Y_{s}d\tilde C'_{s}\leq 0$. 
These inequalities are  based on 
property \eqref{RBSDE_C} of the DRBSDE, on  the non-decreasingness of (almost all trajectories of) $C$, $\bar C$, $C'$ and $\bar C'$,  and on the inequalities $Y\geq \xi$, $\bar Y\geq \xi$, $Y\leq \zeta$, $\bar Y\leq \zeta$. The details, which are similar to those of the proof of Lemma 3.2 in \cite{MG}, are left to the reader.  
The above observations, together with equation \eqref{eq0_lemma_estimate}, lead to the following inequality: a.s., for all $t\in[0,T]$,  
\begin{equation}\label{eq7_lemma_estimate}
\begin{aligned}
\e^{\beta t}\tilde Y_t^2+\int_{]t,T]} \e^{\beta s} \tilde Z_s^2 ds
 +
 \int_{]t,T]} \e^{\beta s} d\langle \tilde h^c \rangle _s 
 &\leq \varepsilon^2 \int_{]t,T]}  \e^{\beta s}\tilde f^2(s) ds-\sum_{t<s\leq T}\e^{\beta s}(\tilde Y_s-\tilde Y_{s-})^2\\
 &-( \tilde M_T-\tilde M_t).
\end{aligned}
\end{equation}

from which we derive  estimates for $\|\tilde Z\|^2_\beta$,
$\|\tilde k\|^2_{\nu,\beta}$, $\|\tilde h\|^2_{\beta, \mathcal{M}^2}$, and then an estimate for $\vvertiii{\tilde Y}^2_\beta.$
\paragraph*{ Estimate for $\|\tilde Z\|^2_\beta$, 
$\|\tilde k\|^2_{\nu,\beta}$ and $\|\tilde h \|^2_{\beta, {\cal M}^2}$}
Note first that we have: 
\begin{equation*}
\begin{aligned}
  \sum_{t< s\leq T}\e^{\beta s}(\Delta \tilde h_s) ^2&+\int_{]t,T]} \e^{\beta s} ||\tilde k_s||_{\nu}^2 ds-\sum_{t<s\leq T}\e^{\beta s}(\Delta \tilde Y_s)^2
= -\sum_{t<s\leq T}\e^{\beta s}(\Delta \tilde A_s- \Delta \tilde A'_s) ^2\\
&-\int_{]t,T]} \e^{\beta s}\int_E \tilde k_s^2(e)  \tilde N(ds,de)
  -2 \sum_{t<s\leq T}\e^{\beta s}  (\Delta \tilde A_s- \Delta \tilde A'_s) \Delta \tilde h_s -2 \sum_{t<s\leq T}\e^{\beta s}\tilde k_s(p_s) \Delta \tilde h_s ,
\end{aligned}
\end{equation*}
where, we have used the fact that $N(\cdot,de)$ "does not have jumps in common" with the predictable processes $A$ and $A'$, since $N(\cdot,de)$ jumps only at totally inaccessible stopping times.\\ 
By adding the term $\int_{]t,T]} \e^{\beta s}|| \tilde k_s||_{\nu}^2 ds$ $+   \sum_{t< s\leq T}\e^{\beta s}(\Delta \tilde h_s) ^2$ on both sides of inequality \eqref{eq7_lemma_estimate}, by using the above computation   and the well-known equality $[\tilde h]_t= \langle \tilde h^c \rangle_t+ \sum (\Delta \tilde h)_s^2$, we get 
\begin{equation}\label{eq_last_000}
\begin{aligned}
 \e^{\beta t}\tilde Y_t^2+\int_{]t,T]} \e^{\beta s} \tilde Z_s^2 ds &+\int_{]t,T]} \e^{\beta s}|| \tilde k_s||_{\nu}^2 ds  +  \int_{]t,T]} \e^{\beta s} d[ \tilde h] _s 
\leq \varepsilon^2 \int_{]t,T]}  \e^{\beta s}\tilde f^2(s) ds - (M'_T- M'_t)\\
& -2 \sum_{t<s\leq T}\e^{\beta s} (\Delta \tilde A_s- \Delta \tilde A'_s) \Delta \tilde h_s -2 \int_t^T d[\tilde h\,, \int_0^\cdot \int_E \e^{\beta s}\tilde k_s(e) \tilde N(ds,de)\,]_s,
\end{aligned}
\end{equation}
with $M'_t= \tilde M_t + \int_{]t,T]} \e^{\beta s}\int_E  \tilde k_s^2(e)  \tilde N(ds,de)$ (where $\tilde M$ is given by \eqref{defM}).\\
  By classical arguments, which  use Burkholder-Davis-Gundy inequalities, we can show that the local martingale $M'$ is a martingale.
Moreover, since $\tilde h$  $\in$ 
 ${\cal M}^{2,\bot}$, by Remark \ref{2i}, 
  we derive that  the expectation of the last term of the above inequality \eqref{eq_last_000} is equal to $0$.
Furthermore, since $\tilde h$ is a martingale, for each predictable stopping time $\tau$, we have 
$E[\Delta \tilde h_\tau / {\cal F}_{\tau-}]=0$ (cf., e.g., Chapter I, Lemma (1.21) in \cite{J}). Moreover, 
since $\tilde A$ and $\tilde A'$ are predictable, we get that $\Delta \tilde A_\tau- \Delta \tilde A'_\tau$ is ${\cal F}_{\tau-}$-measurable (cf., e.g., Chap I (1.40)-(1.42) in \cite{J}), which implies that
$E[(\Delta \tilde A_\tau- \Delta \tilde A'_\tau) \Delta \tilde h_\tau / {\cal F}_{\tau-}]=
(\Delta \tilde A_\tau- \Delta \tilde A'_\tau) E[ \Delta \tilde h_\tau / {\cal F}_{\tau-}]=0$. We thus get 
$E[\sum_{0< s \leq T}\e^{\beta s}(\Delta \tilde A_s - \Delta \tilde A'_s) \Delta \tilde h_s]=0.$

By applying \eqref{eq_last_000} with $t=0$, and by taking expectations on both sides of the resulting inequality, we obtain
$\tilde Y_0^2+ \|\tilde Z\|^2_\beta+\|\tilde k\|^2_{\nu,\beta}+  \|\tilde h \|^2_{\beta, {\cal M}^2} \leq \varepsilon^2\|\tilde f\|^2_\beta.$        
We deduce that $\|\tilde Z\|^2_\beta\leq \varepsilon^2\|\tilde f\|^2_\beta$, $ \|\tilde k\|^2_{\nu,\beta}\leq \varepsilon^2\|\tilde f\|^2_\beta  $ and $\|\tilde h \|^2_{\beta, {\cal M}^2} \leq \varepsilon^2\|\tilde f\|^2_\beta$, which are the desired estimates \eqref{eq_initial_Lemma_estimate}.


The proof of the estimate for $\vvertiii{\tilde Y}^2_\beta$ follows from inequality \eqref{eq7_lemma_estimate} and the estimates for $\|\tilde Z\|^2_\beta$, 
$\|\tilde k\|^2_{\nu,\beta}$ and $\|\tilde h \|^2_{\beta, {\cal M}^2}$, by using exactly the same arguments 
as in the proof of the a priori estimates for reflected BSDEs given  in \cite{MG3} (cf. the Appendix). The proof is thus complete.
\fproof

\paragraph{Proof of Theorem \ref{exiuni}}
 For each $\beta>0$, we denote by $\mathcal{B}_\beta^2$ the space $\mathcal{S}^2 \times \H^2\times \H^2_\nu$ which we  equip with the norm $\|(\cdot,\cdot,\cdot) \|_{\mathcal{B}_\beta^2}$ defined by 
$\| (Y, Z, k)\|_{\mathcal{B}_\beta^2}^2:=\vvertiii {Y}_{\beta}^2  +   \| Z\|_{\beta}^2+\|k \|_{\nu,\beta}^2, $ for  $(Y,Z,k)\in \mathcal{S}^2 \times \H^2\times \H^2_\nu.$ Since $(\H^2, \|\cdot\|_\beta)$, $(\H^2_\nu,\|\cdot\|_{\nu,\beta})$, and  $(\mathcal{S}^2, 
\vvertiii {\cdot}_{\beta})$ are Banach spaces, it follows that
$(\mathcal{B}_\beta^2, \|\cdot \|_{\mathcal{B}_\beta})$ is a Banach space.

We define a mapping $\Phi$ from $\mathcal{B}_\beta^2$ into
itself as follows: for a given $(y, z, l) \in \mathcal{B}_\beta^2$, we set $\Phi (y, z, l):= (Y, Z, k)$, where $Y, Z, k$ 
are
the first three components of the solution $(Y, Z,k, h,A,C,A',C')$ to the DRBSDE  associated with driver $f_s:= 
f(s, y_s, z_s, l_s)$ and with the pair of admissible barriers $(\xi,\zeta)$. 
The mapping $\Phi$ is well-defined by Theorem  \ref{OY}.  

%
  
%
%
Using the estimates provided in Lemma \ref{Lemma_estimate} and following the same arguments as in the proof of Theorem 3.4 in \cite{MG}, 
we derive that for a suitable choice of the parameter $\beta>0$ the mapping  $\Phi$ is a contraction from the Banach space $\mathcal{B}_\beta^2$ into itself.  

By the Banach fixed-point theorem, we get that $\Phi$ has a unique fixed point in $\mathcal{B}_\beta^2$, denoted by $(Y,Z,k)$, that is, such that $(Y,Z,k)= \phi (Y,Z,k)$. By definition of the mapping $\Phi$, the process $(Y,Z,k)$ is thus equal to the first three components of the solution $(Y, Z,k, h, A, C, A', C' )$ to the DRBSDE associated with the driver process $g(\omega,t):=f(\omega, t,Y_t(\omega),Z_t(\omega),k_t(\omega))$ and with the pair of barriers $(\xi,\zeta)$. This property first implies that $(Y,Z,k,h,A,C,A',C')$ is the unique solution 
to the DRBSDE with parameters $(f,\xi,\zeta)$. 
%
%
\fproof

\paragraph{Proof of the statement of  Proposition \ref{Rk_identification}}
 Let $(A,C)$ (resp. $(A',C')$) be the {\em Mertens process} associated with the strong supermartingale ${\cal X}^{f}$ (resp. ${\cal X}^{'f}$), that is satisfying 
$${\cal X}^{f}_t  = E[A_T-A_t +C_{T-} -C_{t-} \,\vert \, {\cal F}_t ]\, (\text{resp. }
 {\cal X}^{'f}_t = E[A'_T-A'_t +C'_{T-} -C'_{t-} \,\vert \, {\cal F}_t ]), \text{ for all } t\in[0,T].$$
 We have to show that $A, C, A'$ and $C'$ are equal to the four last coordinates of the solution of the DRBSDE associated with  parameters  $(\xi,\zeta,f)$.
To this purpose, we apply the same arguments as those used in 
the proof of Proposition \ref{critereexistence} to $X={\cal X}^{f}$ and $X'= {\cal X}^{'f}$. Let $B$, $D$, $B'$ and $D'$ be defined as in this proof. Set $H_t:=E[B_T-B_t +D_{T-} -D_{t-} \,\vert \, {\cal F}_t ]$ and  
$H'_t:=E[B'_T-B'_t +D'_{T-} -D'_{t-}\,\vert \, {\cal F}_t ]$. Since $dB_t << dA_t$, $dB'_t << dA'_t$, $dD_t << dC_t$ and $dD'_t << dC'_t$, we have $H \leq {\cal X}^{f}$ and $H' \leq {\cal X}^{'f}$. Moreover, $H-H'={\cal X}^{f}-{\cal X}^{'f}$, which yields that $\tilde{\xi}^{f} \leq H-H' \leq \tilde{\zeta}^{f}$. By the minimality property of  $({\cal X}^{f}, {\cal X}^{'f})$ (cf. the last assertion of Proposition \ref{seq}), we derive that $H = {\cal X}^{f}$ and $H'= {\cal X}^{'f}$. Hence, $B = A$, $B' = A'$, $D = C$ and $D'= C'$. 
By the properties of $B, B', D,$ and $D'$, we thus get $dA_t \perp dA'_t$ and  $dC_t \perp dC'_t$.
Let now $Y$ be defined by \eqref{oY} with $X={\cal X}^{f}$ and $X'= {\cal X}^{'f}$, and let $( Z,k,h)$ be defined as in the proof of Proposition \ref{critereexistence}. The process
 $(Y,Z,k,h,A,C,A',C')$ is then the solution of the doubly reflected BSDE with parameters $(f,\xi,\zeta)$. The proof is thus complete. 
\fproof

 \paragraph{Proof of Corollary \ref{optimal}} 
 Let $\theta \in \T_0$. By Theorem \ref{caracterisation}, we have $Y_{\theta}=\overline{V}(\theta)=\underline{V}(\theta) $ a.s. 
Moreover, by Proposition \ref{continuous} (iii), since  $ \xi $ (resp. $\zeta$) is left u.s.c.(resp. left l.s.c.) along stopping times, it follows that the nondecreasing process $A$ (resp. $A'$) is continuous.
Since $\sigma_{\theta}^*\leq \overline{\sigma}_{\theta}$ a.s. (cf. Remark \ref{etoilebar}),  by Proposition \ref{lalabis} (first assertion), the process $(Y_t, \,\theta \leq t \leq
\tau \wedge \sigma_{\theta}^*)$ 
is a strong $\mathbfcal{E}^{^{f}}$-supermartingale. We thus get
\begin{equation}\label{supermart}
Y_{\theta} \, \geq \, \mathbfcal{E}^{^{f}}_{_{\theta,\tau \wedge \sigma_{\theta}^*}}[Y_{\tau \wedge \sigma_{\theta}^*}] \quad \text{ a.s. }
\end{equation}
Since $Y \geq \xi$ and $Y_ { \sigma_{\theta}^* } = \zeta_{  \sigma_{\theta}^*}$ a.s.\,(by Proposition \ref{lalabis}),
we also have
$$ Y_{\tau \wedge {\sigma_{\theta}^*}}= Y_{\tau}\textbf{1}_{\tau \leq \sigma_{\theta}^*}+Y_{\sigma_{\theta}^*}\textbf{1}_{\sigma_{\theta}^* < \tau} \, \geq \, \xi_{\tau}\textbf{1}_{\tau \leq \sigma_{\theta}^*}+\zeta_{\sigma_{\theta}^*}\textbf{1}_{\sigma_{\theta}^* < \tau}=I(\tau, \sigma_{\theta}^*)  \quad \text{ a.s. }$$
By inequality \eqref{supermart} and the non decreasing property of $\mathbfcal{E}^{^{f}}$, we get
$Y_{\theta} \, \geq \, \mathbfcal{E}^{^{f}}_{_{\theta,\tau \wedge \sigma_{\theta}^*}}[I (\tau, \sigma_{\theta}^*)]$ a.s.\,
Similarly, one can show that for each $\sigma \in \T_{\theta}$, we have:
$
Y_{\theta} \,  \leq \,  \mathbfcal{E}^{^{f}}_{_{\theta,\tau^*_{\theta} \wedge \sigma}}[I (\tau^*_{\theta}, \sigma)]$ a.s.\,It follows that $( \tau^*_{\theta},  \sigma^*_{\theta})$ is a saddle point at time $\theta$. 
Similarly, using Proposition \ref{lalabis}, it can be shown that $(\overline{\tau}_{\theta}, \overline{\sigma}_{\theta})$ is a   saddle point at time $\theta$, which ends the proof.
\fproof

\begin{lemma}\label{Moki_lem} Let $\xi$ be an optional process which can be written 
$\xi_t := M_t + \alpha_t + \gamma_{t^-}$, where $M$  a square integrable martingale, $\alpha$ and $\gamma$ are
RCLL adapted processes with $\alpha_0=\gamma_{0^-}=0$, and with
 square integrable  total variation that is, $E( \vert \alpha \vert_T^2) < \infty$ and $E( \vert \gamma \vert_T^2) < \infty$.
Then, the process $\xi$ can be written as the difference of two non negative square integrable strong supermartingales.
\end{lemma}
\dproof
By the above Proposition A.7 in \cite{DQS2}, there exists an unique  pair $(A, A')$ of non decreasing RCLL adapted processes such that  $A_T$ and $A'_T$ are square integrable, and
 $\alpha=A'-A$ with $dA_t \perp dA'_t$. 
 Similarly, there exists an unique  pair $(C, C')$ of non decreasing RCLL adapted processes such that  $C_T$ and $C'_T$ are square integrable, and
 $\gamma=C'-C$ with $dC_t \perp dC'_t$. 
 
The processes $H$ and $H'$ defined by $H_t:=E[\xi^+_T+A_T-A_t + C_{T^-} - C_{t^-} | \mathcal{F}_t]$ and 
$H'_t:=E[\xi^{-}_T+A'_T - A'_t + C'_{T^-} - C'_{t^-}| \mathcal{F}_t]$ are non negative strong supermartingales belonging to $ \mathcal{S}^2$. 
Moreover, we have $\xi_t$ $=$ $E[ \xi_T+\alpha_T - \alpha_t + \gamma_{T^-} - \gamma_{t^-} | \mathcal{F}_t]$ $=$ $ H_t - H'_t$, which gives the desired result. 
\fproof

\begin{remark}\label{Moki_rem} From this property, we derive that if $\xi$ and $\zeta$ 
are optional processes in ${\cal S}^2$ with $\zeta_T= \xi_T$ and $\xi \leq \zeta$, and such that $\xi$ (or $\zeta$) satisfies the assumptions of Lemma \ref{Moki_lem}, then the pair $(\xi,\zeta)$ satisfies Mokobodzki's condition  \eqref{Moki}.
\end{remark}

 \begin{example}\label{contre-exemple}
 We give a simple example of classical Dynkin game in a deterministic framework which does not admit a value. 
 Let us suppose that $T=1$. Let $\xi$ and $\zeta$ be two maps defined on $[0,1]$ by 
 $\xi(t):= \frac{1}{2} {\bf 1}_{\mathbb{Q}} (t)$ for all $t \in [0,1]$,
 $\zeta(s):= 1 +  {\bf 1}_{\mathbb{R} \backslash \mathbb{Q}} (s)$ for all $s \in [0,1)$, and 
 with $\zeta (1)=  \frac{1}{2}$. \footnote{Actually, the value of $\zeta (1)$ does not intervene in the criterion.}
 Note that the pair  $(\xi, \zeta)$  is admissible and satisfies Mokobodzki's condition. Moreover, $\xi$ and $-\zeta$ are not right-u.s.c.  
We have: \\
 For each $s \in {\mathbb{Q}} \cap [0,1]$ (resp. $s \in (\mathbb{R} \backslash \mathbb{Q}) \cap [0,1]$), we have $\sup_{t \in[0,1]} I(t,s)= 1$ (resp. $= 2$). 
 Hence,  $\inf_{s \in[0,1]} \sup_{t \in[0,1]} I(t,s)=1$.\\
 On the other hand, for each $t \in {\mathbb{Q}} \cap [0,1]$ (resp. $t \in (\mathbb{R} \backslash \mathbb{Q}) \cap [0,1]$), we have 
 $\inf_{s \in[0,1]}  I(t,s)= 1/2$ (resp. $= 0$). Hence, 
$\sup_{t \in[0,1]} \inf_{s \in[0,1]}  I(t,s)= 1/2$, from which we derive that there does not exist a value for this deterministic Dynkin game.
 \end{example}



\end{document}